\documentclass{article}
\usepackage[text={155mm,220mm},centering]{geometry}
\usepackage{algorithm2e}
\usepackage[utf8]{inputenc}

\usepackage{hyperref} 
\hypersetup{	colorlinks=true,	linkcolor=blue,	filecolor=magenta,	urlcolor=blue,	allcolors=blue,}

\usepackage{yun}

\usepackage{graphicx}
\usepackage{graphics}
\graphicspath{{graphs/}} 
\usepackage{pinlabel}

\def \P{P}

\def \mytheta{\theta}
\def \myf{\phi}
\def \myF{\Phi}

\def \Mul{\operatorname{Mul}}
\def \ker{\operatorname{ker}}

\def \supp{\operatorname{supp}}
\def \KS{D_{\operatorname{KS}}}
\def \MMD{D_{\operatorname{MMD}}}

\def \Child{\operatorname{Child}}
\def \Desc{\operatorname{Desc}}

\renewcommand{\myred}[1]{{#1}}

\title{Minimum $\Phi$-distance estimators for finite mixing measures}

\author{
    Yun Wei$^1$,
    Sayan Mukherjee, 
    XuanLong Nguyen$^2$  
}
\date{%
$^1$Department of Mathematical Science, The University of Texas at Dallas, TX\\
    $^2$ Department of Statistics, University of Michigan, Ann Arbor, MI 
     \\%
    ~\\
}

\begin{document}

\maketitle

\begin{abstract} 
   Finite mixture models have long been used across a variety of fields in engineering and sciences. Recently there has been a great deal of interest in quantifying the convergence behavior of the \emph{mixing measure}, a fundamental object that encapsulates all unknown parameters in a mixture distribution. In this paper we propose a general framework for estimating the mixing measure arising in finite mixture models, which we term minimum $\Phi$-distance estimators. We establish a general theory for the minimum $\Phi$-distance estimator, where sharp probability bounds are obtained on the estimation error for the mixing measures in terms of the suprema of the associated empirical processes for a suitably chosen function class $\Phi$. 
    Our framework includes several existing and seemingly distinct estimation methods as special cases using a weakened identifiability condition, but also motivates new estimators. For instance, it extends the minimum Kolmogorov-Smirnov distance estimator to the multivariate setting, and it extends the method of moments to cover a broader family of probability kernels beyond the Gaussian. Moreover, it also includes methods that are applicable to complex (e.g., non-Euclidean) observation domains, using tools from reproducing kernel Hilbert spaces. It will be shown that under general conditions the methods achieve optimal rates of estimation under Wasserstein metrics in either minimax or pointwise sense of convergence; the latter case can be achieved when no upper bound on the finite number of components is given. 
    %
    %
    \myred{Also of interest is a sharp inequality that captures the local information geometry for general mixture models precisely in terms of moment differences between mixing measures.} 
    
\end{abstract}


\tableofcontents

\section{Introduction}

Since the early work of \cite{pearson1894contributions}, finite mixture models have long been used as a modeling tool across a variety of fields in engineering and sciences \cite{peel2000finite}. 
They are deployed in clustering analysis \cite{bishop2006pattern}, as well as modeling heterogeneous data distributions, e.g.,  \cite{van1996rates,genovese2000rates,ghosal2001entropies,wei2022convergence}. Recently there has been a great deal of interest in quantifying the convergence behavior of the \emph{mixing measure}, a fundamental object that encapsulates all unknown parameters in a mixture distribution \cite{chen1995optimal, nguyen2013convergence, ho2019singularity,heinrich2018strong, wu2020optimal}. In this paper we propose a general framework for estimating the mixing measure arising in finite mixture models. This framework not only includes many existing estimation methods as special cases but also motivates new ones of interest.

The estimation framework that we study involves a general notion of distance on the space of latent mixing measures, which requires the evaluation of probability measures using a suitable class of test functions. The function class will be generically named $\Phi$ in this paper, and the corresponding distance the $\Phi$-distance. Specializing the function class $\Phi$ to concrete instances leads to well-known estimation methods.
For instance, a special case of the $\Phi$-distance is the Kolmogorov–Smirnov (KS) distance for univariate distributions, which results in the minimum KS distance method \cite{deely1968construction, chen1995optimal, Chen1996-ia, heinrich2018strong}. Another special case is the $\ell_\infty$ distance of moment vectors for finite mixture distributions, which yields the so-called denoised method of moments estimator \cite{wu2020optimal}. It is worth noting that the minimax optimality analysis of both methods has only been established recently \cite{heinrich2018strong,wu2020optimal}.

The minimum $\Phi$-distance estimator studied in this paper is considerably more general and more broadly applicable than the aforementioned works. It can be applied to both multivariate parameter spaces and multivariate domains of observed data, as well as to general families of probability kernels for modeling the mixture components. Moreover, the minimum $\Phi$-distance estimation framework leads to methods that are applicable to complicated (e.g., non-Euclidean) observation domains. In particular, we study a minimum distance estimator based on the maximum mean discrepancy (MMD), a particular $\Phi$-distance that arises in a different context (of learning with reproducing kernel Hilbert spaces) \cite{gretton2012kernel}. \myred{Minimum MMD distance estimators have been studied in \cite{briol2019statistical, cherief2022finite} where they focus on density estimation rate, not parameter convergence rates as in our paper. }

A general theory for the general minimum $\Phi$-distance estimator introduced in this paper involves obtaining sharp probability bounds on the estimation error for the mixing measures in terms of the behavior of the suprema of the empirical processes associated with the function class $\Phi$. As a direct consequence of this general theory, we are able to obtain the optimal rates of estimation for all three specific estimators mentioned above. Notably, optimality is established in both a uniform convergence (minimax) sense and a pointwise sense. For instance, for the minimum KS-distance estimator, we generalize the existing results of \cite{heinrich2018strong} from univariate to multivariate scenario, while relaxing some of the assumptions of their theorems. For the denoised moment method of \cite{wu2020optimal}, under our general theory it becomes possible to extend  
to replicate the results for Gaussian mixtures to broader families of kernels, namely, the natural exponential families with quadratic variance functions (NEF-QVF) \cite{morris1982natural}. \myred{Some other relevant papers on moment methods in mixture models are \cite{Kalai2010-jn, Anandkumar2012-nq, Hardt2015-rg}. We also apply our general theories to establish a convergence rate for multi-dimensional Gaussian mixture models for estimators based on moment tensors, which have been studied previously in \cite{pereira2022tensor} but lack theoretical convergence rate results. } 

An interesting aspect about estimation in finite mixture models, in the setting where the number of parameters is unknown, is that the optimal minimax rate is typically far slower than the pointwise optimal rate of parameter estimation. To achieve the pointwise optimal rate of estimation for the mixing measure, we study a plug-in estimator which consists of two steps: first, obtain a consistent estimate of the number of mixture components, and second, estimate the mixing measure based on the former
estimate. Both steps make essential use of the chosen $\Phi$ function classes. It will be shown that under quite general conditions, the pointwise rate of convergence of the proposed estimator is indeed the optimal $n^{-\frac{1}{2}}$ under the $\ell_1$ Wasserstein metric, which is much faster than the minimax (optimal) rate of $n^{-\frac{1}{2(2d_1-1)}}$, where $d_1$ is the effective degree of freedom representing the amount of overfitting by the mixture model. Such a phenomenon was established for the minimum KS-distance estimator in \cite{heinrich2018strong} and will be proved here for the general minimum $\Phi$-distance estimator. 

Another noteworthy, perhaps deeper aspect, of our general theory for the minimum $\Phi$-distance estimator is illuminated by the dual and separate roles that the class $\Phi$ of test functions plays. On the one hand, $\Phi$ has to be sufficiently rich to enable the identification of the mixing measure --- this is related to the condition of strong identifiability employed in the existing literature \cite{chen1995optimal,Rousseau-Mengersen-11, nguyen2013convergence,ho2019singularity,heinrich2018strong}. On the other hand, since the mixing measure is not observed directly --- only samples of the mixture distributions are given from which the elements in $\Phi$ can be estimated. Thus $\Phi$ has to be sufficiently small if the distance is to be evaluated efficiently from the empirical data. Exploiting the balance between these two forces --- one is theoretical and another computational --- allows one to design suitable $\Phi$ classes, as well as obtain a sharp analysis of the corresponding estimator under intrinsic identification conditions. In fact, such conditions are shown to be weaker than the more standard strong identifiability conditions considered in the literature. \myred{To the best of our knowledge, we are the first to weaken the standard strong identifiability condition in the literature. This relaxation is particularly useful when the function class $\Phi$ is chosen to be of finite cardinality. Usefulness of this relaxation is demonstrated by method of moments, and also by studying the mixture of Bernoulli/multinomial distributions, which outperforms the existing results in the literature \cite{manole2021estimating}.}

\myred{One key technical inequality we established is the following: for any two mixing measures $G,H$ in some suitable space,  
$$W_{2k-1}^{2k-1}(G,H)\leq C_1 \mbf_{2k-1}(G,H) \leq C_2 \sup_{\myf\in \myF} \left|\int \myf dG - \int \myf dH  \right| \leq C_3 \mbf_{2k-1}(G,H).  $$
By choosing a suitable function class $\myF$, $\sup_{\myf\in \myF} \left|\int \myf dG - \int \myf dH  \right|$ represents the distance between two mixture densities, and thus the above characterizes precisely the information geometry of mixture models in terms of the moment difference between the corresponding mixing measures, $\mbf_{2k-1}(G,H)$. Moreover, the moment difference is further bounded below by the Wasserstein distance. We note that related inequalities are obtained for Gaussian mixture models \cite[Theorem 4.2]{doss2020optimal}, but such inequalities for general mixture models are new, to the best of our knowledge. } \myred{The new inequality above also proposes to use moment difference as a candidate to measure the mixing measure estimation errors in the sense that moment difference captures the local information of finite mixture models, and it also yields the popular Wasserstein distance errors  in the literature.}

Finally, we note several other related strands of recent work regarding mixture model estimation. The theoretical analysis of pointwise convergence behavior in finite mixture models has been explored extensively in under complex model settings \cite{Rousseau-Mengersen-11, nguyen2013convergence,ho2016strong,ho2016convergence,wei2022convergence}. The development of methods for estimating the unknown number of mixture components continues to be of interest, as shown in \cite{guha2021posterior,manole2021estimating,cai2020power,cai2021finite} and the references therein.
   Other researchers \cite{hettmansperger2000almost, elmore2004estimating,cruz2004semiparametric,Jochmans1996Nonparametric,vandermeulen2019operator,ritchie2020consistent,aragam2021uniform,bing2022sketched} studied nonparametric mixtures, i.e., no parametric forms for the probability kernels for a given component are assumed, and the focus is on the problem of density estimation due to the nonparametric setup. 
By contrast, in this paper we study mixtures  with the parametric form of component distribution imposed, since in practice prior knowledge on the component distributions might be available. Moreover, we investigate the convergence behavior of parameter estimates, which are generally more challenging to address than that of the mixture density function, as pointed out in \cite{nguyen2013convergence, ho2016strong, ho2016convergence, heinrich2018strong, ho2019singularity,wei2022convergence,doss2020optimal,ashtiani2018nearly}.

The rest of the paper will proceed as follows. Section~\ref{sec:unifiedframework} presents the minimum $\Phi$-distance estimation framework and develop a general theory of the analysis of uniform convergence (minimax) rates for this class of estimators. Section~\ref{sec:examples} presents three specific instances of the minimum $\Phi$-distance estimators, including the minimum KS-distance method, the denoised moment method \myred{(using one-dimensional and higher-dimensional moment tensors)}, and a novel estimator based on the MMD distance. In Section~\ref{sec:pointwiserate} we obtain the pointwise rate of convergence for the mixing measures, by studying an estimation method based on the minimum $\Phi$-distance estimates. Section \ref{sec:discussion} outlines several open questions, as well as related results of potential interest. All proofs are given in the Appendix.

\subsection{Notation} 
Denote the set of natural numbers by $\Nb=\{0,1,\ldots\}$ 
the set $[k]:=\{1,2,\ldots,k\}$.  $\Nb_+$ denotes the positive natural numbers. The maximum between two numbers is denoted by $a\vee b$ or $\max\{a,b\}$. The minimum between two numbers is denoted by $a\wedge b$ or $\min\{a,b\}$. $\Gamma(x)$ denotes the Gamma function.
$\tilde{\Theta}^\circ$ is the interior of a set $\tilde{\Theta}$. The complement of a set $A$ is denoted by $A^c$. For a finite set $A$, $|A|$ denotes its cardinality. $1_A(x)$ for a set $A$ is the indicator function taking the value $1$ when $x\in A$ and $0$ otherwise. $1_{p\geq a}$ for a logical statement like $p\geq a$ is  $1$ if the statement is true and $0$ otherwise. 

The vector of all zeros is denoted as
$\bm{0}$ (in bold). Any vector $x\in \Rb^d$ is a column vector with its $i$-th coordinate denoted by $x^{(i)}$. The span of a vector is denoted $\operatorname{span}(v)=\{a v| a\in \Rb  \}$. The inner product between two vectors $a$ and $b$ is denoted by $a^\top b$ or $\langle a, b\rangle$.
The multi-index notation for $\alpha\in \Nb^q$
imposes the following 
$$
|\alpha|:=\sum_{i\in [q]} |\alpha^{(i)}|, \quad \alpha! := \prod_{i\in [q]} \alpha^{(i)}!, \quad \theta^\alpha :=  \prod_{i\in [q]} \left(\theta^{(i)}\right)^{\alpha^{(i)}},
$$
\myred{where $\theta\in \Rb^q$.} Denote $\Ic_k:=\{\alpha\in \Nb^q \mid  |\alpha|\leq k \}$. For two multi-indices $\alpha,\gamma\in \Nb^q$, $\alpha\leq \gamma$ if and only if $\alpha^{(i)}\leq \gamma^{(i)}$ for any $i\in [q]$. For a multi-index $\alpha$, the operator $D^{\alpha}$ means partial derivative of order $\alpha^{(i)}$ to the $i$-th coordinate. Note in this paper that the partial derivative is always with respect to $\theta$, i.e. $D^\alpha p(x \mid \theta)=\frac{\partial^\alpha }{\partial \theta^{\alpha}} p(x \mid \theta)$.  

For any probability measure $\P$ and $Q$ on measure space  $(\Xf,\Xc)$ with densities respectively $p$ and $q$  with respect to    some base measure $\lambda$, the variational distance between them is
	$V(\P,Q) = \sup_{A\in \Xc} |\P(A)-Q(A)| = \frac{1}{2} \int_{\Xf}  |p(x)-q(x)|d\lambda$. 

We denote the Dirac measure at $\theta$ as $\delta_{\theta}$. For a finite signed (discrete) measure $G=\sum_{i\in [k]} {p_i}\delta_{\theta_i}$ on $\Rb^q$, its $\alpha$-th moment is $m_\alpha(G)= \int \theta^{\alpha} \textrm{d}G(\theta) = \sum_{i\in [k]} {p_i}\theta_i^\alpha\in \Rb^q$. Denote by $\mbf_k(G) :=(m_\alpha(G))_{\alpha\in \Ic_k}\in \Rb^{|\Ic_k|}$ the vector of all $\alpha$-th moments of $G$ for $\alpha\in \Ic_k$. We also write $m_\alpha(Z)=m_\alpha(G)$ or $\mbf_k(Z)=\mbf_k(G)$ when $Z\sim G$, i.e., $Z$ is a random variable drawn from probability distribution $G$. In general, for a measurable function $\myf$ defined on $\Theta$, its integral w.r.t. a distribution $G=\sum_{i=1}^k p_i \delta_{\theta_i} \in \Ec_k(\Theta)$ is denoted by $G\myf:=\int \myf dG = \sum_{i=1}^k p_i \myf(\theta_i)$. The notation $G\myf$ is used to emphasize that $G$ can be viewed as an linear operator on measurable functions on $\Theta$. 

Denote by $G-\theta:=\sum_{i\in [k]} {p_i}\delta_{\theta_i-\theta}$  the signed measure obtained by shifting the support points of $G$ by $-\theta$. Denote $S_\epsilon G:=\sum_{i\in [k]} {p_i}\delta_{\epsilon \theta_i}$ to be the signed measure obtained by scaling the support points of $G$ by $\epsilon$, where $S_{\epsilon}$ is viewed as an operator on signed measures. 

Denote by $C(\cdot)$ or $c(\cdot)$  a positive finite constant depending only on its parameters and the probability kernel $\{\Pb_\theta\}_{\theta\in \Theta}$. In the presentation of inequality bounds and proofs,  they may differ from line to line. 



\section{A general framework for estimation}
\label{sec:unifiedframework}

Consider a family of probability distributions  $\{\Pb_{\theta}\}_{\theta\in \Theta}$ on measurable space $(\Xf,\Xc)$, 
where $\theta$ are the parameters of the family and $\Theta\subset \Rb^q$ is the parameter space. Throughout this paper it is assumed that the map $\theta\mapsto \Pb_{\theta}$ is injective.  The space of all discrete probability distributions with exactly (or at most) $k$ distinct atoms on $\Theta$ is denoted by $\Ec_k(\Theta)$ (respectively, $\Gc_k(\Theta)$). It is clear that $\Gc_k(\Theta)=\cup_{\ell\in [k]} \Ec_\ell(\Theta)$. Given a finite discrete probability measure $G=\sum_{i=1}^{k} p_i \delta_{\theta_i} \in \Ec_{k}(\Theta)$, the mixture distribution on $(\Xf,\Xc)$ induced by $G$ is given by 
$\Pb_G(\textrm{d}x) = \sum_{i=1}^{k} p_i \Pb_{\theta_i}(\textrm{d}x)$. $G$ is called the mixing measure corresponding to the mixture distribution $\Pb_G$. Given i.i.d. observed samples $X_1,\ldots,X_n\overset{\iidtext}{\sim} \Pb_{G^*}$ for some fixed but unknown mixing measure $G^*$, the goal is to estimate \myred{$G^*=\sum_{i\in [k^*]}p^*_i\delta_{\theta^*_i}\in \Ec_{k^*}(\Theta)$}, which contains all the parameters of interest $k^*$, $p^*_i$, $\theta^*_i$. To be clear, in this paper we will not assume the number of mixture components $k^*$ is known. We will construct estimators on $\Gc_k(\Theta)$ for some $k$ and assume $G^*\in \Gc_k(\Theta)$ so $k$ is a known upper bound for $k^*$; there is one general result where $G^*\not\in \Gc_k(\Theta)$ is not required and we will point this out. Our task in this paper is to study the convergence rate of estimating the mixing measure $G^*$. In particular we will pay attention to the dependence on the upper bound $k$.

In order to quantify the convergence of mixing measures in mixture models, a useful device is a suitably defined optimal transport distance~\cite{nguyen2013convergence,villani2003topics}.
Consider the Wasserstein-$\ell$ distance with respect to (w.r.t.) the Euclidean distance on $\Theta$: for all $ G=\sum_{i=1}^{k}p_i\delta_{\theta_i}, G'=\sum_{i=1}^{k'}p'_i\delta_{\theta'_i} 
$, we define
\begin{equation}
W_\ell(G,G') = \left(\min_{\bm{q}} \sum_{i=1}^{k}\sum_{j=1}^{k'} q_{ij}\|\theta_i-\theta'_j\|_2^\ell\right)^{1/\ell}, \label{eqn:Wpdef}
\end{equation}
where the infimum is taken over all joint probability distributions $\bm{q}$ on $[k]\times [k']$ such that, when expressing $\bm{q}$ as a $k\times k'$ matrix, the marginal constraints hold: $\sum_{j=1}^{k'} q_{ij} = p_i$ and $\sum_{i=1}^{k} q_{ij} = p'_j$. 
 We state $G_{n}\overset{W_\ell}{\to}G$ if $G_{n}$ converges to $G$ under the $W_\ell$ distance.

\subsection{Minimax and pointwise convergence bounds} 

A standard way for characterizing the difficulty of an estimation problem is via minimax 
lower bounds for the quantity of interest. An estimation procedure is then evaluated against this metric of performance; the procedure is considered optimal in the \emph{minimax sense} if the corresponding minimax estimation upper bound guarantee matches the minimax lower bound under the same setting. It must be noted that for mixture models, the optimal minimax estimation rate is typically much slower than the optimal \emph{pointwise} estimation rate for the mixing measure. Thus, in theory a ``minimax optimal" procedure is not necessarily optimal in the sense of pointwise convergence, and vice versa. To fully assess the quality of a proposed estimation procedure, in this paper we will characterize the proposed estimation procedure using both types of convergence bounds.

For finite mixture models, the optimal pointwise convergence rate for the Wasserstein metrics are the parametric $n^{-1/2}$ under quite general settings. Various estimation methods have been shown to achieve this rate of pointwise convergence (possibly up to a logarithm factor) \cite{heinrich2018strong,ho2020robust,guha2021posterior}. In this paper, the analysis of pointwise convergence will be deferred to Section~\ref{sec:pointwiserate}.
On the other hand, a precise minimax bound for overfitted finite mixture models may vary with the model setting. The first such example \myred{for general mixture models} was established by \cite{heinrich2018strong} in the univariate parameter setting, i.e., when $q=1$ and the mixture distribution is on $\Xf=\Rb$. \myred{Prior to \cite{heinrich2018strong}, there are also related minimax results \cite{Hardt2015-rg, Kalai2010-jn, chen2003tests} for Gaussian mixture models.  }
Here, we shall present a general and somewhat stronger result (comparing to \cite{heinrich2018strong}) that is moreover applicable to the $q \in \mathbb{N}$ setting, and that relies on a weaker assumption on the kernel $\Pb_{\theta}$. 
Within this subsection, assume that $\{\Pb_{\theta}\}_{\theta\in \Theta}$ has density $\{p(x \mid \theta)\}_{\theta\in \Theta}$ w.r.t. a dominating measure $\lambda$ on $(\Xf,\Xc)$.  The following technical assumption imposes a  regularity of the density family $\{p(x \mid \theta)\}_{\theta\in \Theta}$. It restricts the partial derivatives of members in this family. The assumption is quite mild compared to those considered in the existing literature, which will be discussed shortly. Note in this paper the partial derivative is always with respect to $\theta$, i.e., $D^\alpha p(x \mid \theta)=\frac{\partial^\alpha }{\partial \theta^{\alpha}} p(x \mid \theta)$.  

\begin{assump}
We say that the probability kernel $\{p(x \mid \theta)\}_{\theta\in \Theta}$ satisfies Assumption $A(\theta_0,m)$ if 
 1) there exists $b>0$ such that for $\lambda$-a.e. $x$, $p(x \mid \theta)$ is $m
 $-th order continuously differentiable w.r.t. $\theta$ in 
$\{\theta\in \Theta: \|\theta-\theta_0\|_2< b\}$; and 2) there exists a unit vector $\psi\in \Rb^q$ such that 
\begin{equation}
A:= \max_{|\alpha|=m}\ \sup_{\substack{\theta'\in\sp(\psi) \\ \|\theta'\|_2\leq b}}\ \sup_{t\in [0,1]}\int \frac{\left(D^\alpha p(x \mid \theta_0+t\theta')\right)^2}{p(x \mid \theta_0+\theta')}  d\lambda <\infty. \label{eqn:minimaxhypothesis2}
\end{equation}
\end{assump}


Let $\Pc(\Theta)$ be the space of all probability measure on $\Theta$ endowed with Borel sigma algebra. Denote by $\Ef_n$ the set of all estimators (measurable random elements) taking values in $\Pc(\Theta)$ 
based on $\iidtext$ samples $X_1,\ldots,X_n$ from the mixture distribution $\Pb_{G^*}$.  In the following $\Eb_{G^*} f(X_1,\ldots,X_n)$ denotes the expectation when $\{X_i\}_{i\in [n]}\overset{\iidtext}{\sim} \Pb_{G^*}$. 

\begin{thm}[Minimax lower bound]
\label{thm:minimax}
\begin{enumerate}[label=(\alph*)]
\item \label{item:thm:minimaxb}
Suppose that the probability kernel $\{p(x \mid \theta)\}_{\theta\in \Theta}$ satisfies Assumption $A(\theta_0,2k-1)$ for some $\theta_0\in \Theta$. Then for any $n\geq 1$, 
\begin{equation*}
\inf_{\hat{G}_n\in \Ef_n}\ \sup_{ \substack{G^*\in \Gc_k(\Theta)}} \Eb_{G^*} W_1(\hat{G}_n,G^*) \geq C(A,q,k) n^{-\frac{1}{4k-2}}. 
\end{equation*}

\item \label{item:thm:minimaxa}
Consider any $k_0\leq k$ and fix a $G_0\in \Ec_{k_0}(\Theta)$.  Set $\epsilon_n=n^{-\frac{1}{4d_1-2}}$, where $d_1=k-k_0+1$. Suppose that there exists  a support point $\theta_0$ of $G_0$ such that the probability kernel $\{p(x \mid \theta)\}_{\theta\in \Theta}$ satisfies Assumption $A(\theta_0,2d_1-1)$. Then for any $a>0$, for any $n\geq 1$, 
\begin{equation}
\inf_{\hat{G}_n\in \Ef_n}\ \sup_{ \substack{G^*\in \Gc_k(\Theta)\\ W_1(G^*,G_0)<a \epsilon_n }  } \Eb_{G^*} W_1(\hat{G}_n,G^*) \geq C(A,q,d_1) n^{-\frac{1}{4d_1-2}}.  \label{eqn:uniformintegraluppbou}
\end{equation}
\end{enumerate}
\end{thm}

    Part \ref{item:thm:minimaxb} follows directly from part \ref{item:thm:minimaxa} with $G_0=\delta_{\theta_0}$, and the proof of part \ref{item:thm:minimaxa} is in Section \ref{sec:proofofminimaxtheorem}. Part \ref{item:thm:minimaxa} is known as a local minimax lower bound since the true mixing measure $G^*$ is within a shrinking neighborhood of some $G_0\in \Ec_{k_0}(\Theta)$, which can be thought of as prior information that one believes the true mixing measure $G^*$ to lie in. Since $W_1(G^*,G_0)<a \epsilon_n$ and $G^*\in \Gc_k(\Theta)$, for large $n$ we must have $k^*\in [k_0,k]$. Hence the quantity $d_1=k-k_0+1$ is termed an \emph{overfitted index}.  
    The local minimax lower bound (when ignoring the constant multiplier independent of $n$) $n^{-\frac{1}{4d_1-2}}$ depends on the overfitted index: 
    the more accurate the prior information $G_0$ is, the less overfit, the smaller $d_1$, and the smaller the local minimax lower bound. 
    In particular, the slowest local minimax lower bound happens when $k_0=1$, that is when $G_0=\delta_{\theta_0}$, which is also the (global) minimax lower bound $n^\frac{1}{4k-2}$ in Part \ref{item:thm:minimaxb}. Several specific estimators will be shown to have uniform convergence rates matching the minimax lower bounds in Theorem \ref{thm:minimax} up to a constant multiplier so the exponents of $n$ can not be improved. 

\begin{rem}
A similar minimax result for the case $q=1$ is established \cite[Theorem 3.2]{heinrich2018strong}. Theorem \ref{thm:minimax} has several notable improvements. Firstly, Theorem \ref{thm:minimax} works for multivariate parameter spaces. Secondly, our $\epsilon_n$ is smaller than $n^{-\frac{1}{4d_1-2}+\kappa}$ for some $\kappa>0$ as in \cite[Theorem 3.2]{heinrich2018strong} and thus Theorem \ref{thm:minimax} is more general; moreover, the technical assumptions in Theorem \ref{thm:minimax} are also weaker (see the next remark for details). Finally, the proof of Theorem \ref{thm:minimax} appears to be simpler since it does not rely on the local asymptotic normality 
argument as \cite[Theorem 3.2]{heinrich2018strong}. \myeoe
\end{rem}

\begin{rem}
Note that in the univariate parameter setting $q=1$, a similar assumption called  $(p,\alpha)$-smooth in \cite[Definition 2.1]{heinrich2018strong} with $p\in [m]$ and $\alpha=2$ implies  our weaker assumption $A(\theta_0,m)$ for any $m$. In fact \cite[Theorem 3.2]{heinrich2018strong} requires more: $p\in [2k+2]$ (which roughly means more differentiability than $2k-1$ in our results) and $\alpha\in [4]$ (which roughly means stronger integrability condition than our square integrability condition).   \myeoe
\end{rem} 

\subsection{Inverse bounds and implications on convergence rates}
\label{sec:invbou}
In this subsection we introduce a general distance to measure the deviation between two mixing distributions. Then, we introduce inverse bounds, a collection of inequalities which relate this new distance to the Wasserstein distance. Such inverse bounds are useful to derive convergence rates for an estimation procedure.


Consider a family $\myF$ of real-valued functions defined on $\Theta$. Roughly speaking, $\myF$ is a collection of test functions such that for each $\myf\in \myF$, $G\myf$ can be relatively easy to estimate based on data samples from $\Pb_{G}$ (recall the notation  $G\myf:=\int \myf dG$). We shall be more precise about this when we introduce our estimator in Section \ref{sec:mindisest}. Given $\myF$ we use $\sup_{\myf\in \myF} |G \myf-H \myf|$ to measure the deviation between two mixing measures $G$ and $H$. A natural requirement for the test functions is the following property. 
\begin{definition}
\label{def:distin}
	$\Gc_k(\Theta)$ is \textit{distinguishable} by $\myF$ if for any $G \neq H\in \Gc_k(\Theta)$, $\sup_{\myf\in \myF} |G \myf-H \myf|>0$. 
\end{definition}
If $\Gc_k(\Theta)$ is distinguishable by $\myF$, then  $\sup_{\myf\in \myF} |G \myf-H \myf|$ is a distance on $\Gc_k(\Theta)$.

\begin{exa}[Moment deviations between mixing distributions] \label{exa:moment}
Suppose $q=1$ for simplicity in this example. We will consider the general $q$ setting in detail in Section \ref{sec:momentdeviation}. Consider $\myF_2=\{(\theta-\theta_0)^j\}_{j\in [2d_1-1]}$ to be a finite collection of polynomials with $\theta_0$ a fixed constant. 
Then 
\begin{equation}
\sup_{\myf\in \myF} |G \myf-H \myf| = \sup_{j\in [2d_1-1]} |m_j(G-\theta_0)-m_j(H-\theta_0)|=\|\mbf_{2d_1-1}(G-\theta_0)-\mbf_{2d_1-1}(H-\theta_0)\|_{\infty}, \label{eqn:momentdistance}
\end{equation}
which is the maximum deviation of the first $2d_1-1$ moments of $G-\theta_0$ and $H-\theta_0$. Note that one may also include the index $j=0$ in the definition of $\myF_2$.  \myeoe
\end{exa}

\begin{exa}[Integral probability metrics]
\label{exa:divationbetweenmixturedistributions}
Consider $\myF=\left\{ \theta \mapsto \int f_1(x) \Pb_\theta(dx) | f_1 \in \Fc_1 \right\}$ where $\Fc_1$ is some subset of $\Mc$, the space of all measurable functions on $(\Xf,\Xc)$. Note that each $f_1 \in \Fc_1$ defines a function of $\theta$, $\theta \mapsto \int f_1(x) \Pb_\theta(dx)$. Then
\begin{equation}
\sup_{\myf\in \myF} |G \myf-H \myf| = \sup_{f_1 \in \Fc_1} \left|\int f_1 d\Pb_G-\int f_1 d\Pb_H\right|, \label{eqn:integraldistance}
\end{equation}
which is the integral probability metrics (IPM) \cite{muller1997integral, sriperumbudur2012empirical}  between mixture distributions induced respectively by the mixing distributions $G$ and $H$. 

When $\Fc_1=\left\{x \mapsto 1_{B}(x)|B\in \Xc \right\} $, $\eqref{eqn:integraldistance}$ represents the total variation distance $V(\Pb_G,\Pb_H)$.  When the underlying space $(\Xf,\Xc)=(\Rb,\Bc(\Rb))$, the real line endowed with the Borel sigma algebra, and $\Fc_1= \left\{x \mapsto 1_{(-\infty,a]}(x)|a\in \Rb \right\}$,  $\eqref{eqn:integraldistance}$ represents the Kolmogorov-Smirnov (KS) distance $\KS(\Pb_G,\Pb_H)$, which is the maximum deviation of the cumulative distribution functions (CDF) of the mixture distributions. We refer to the $\myF$ in the previous case as $\myF_0$ and will discuss it in detail in Section \ref{sec:minimumdistance estimator}. As we can see, with different choices of $\Fc_1$, we are able to obtain different IPMs. Other IPMs of interest include  Wasserstein-1 distance, Dudley's metric \cite[Chapter 11]{dudley2018real} and maximum mean discrepancy (MMD) \cite{gretton2012kernel}. \myeoe
\end{exa}

A powerful property for $\myF$ to possess, under suitable identification conditions that will be introduced, is a global inverse bound relating $\sup_{\myf\in \myF} |G \myf-H \myf|$ to a Wasserstein distance: 
\begin{equation}
	\inf_{  G\neq H\in \Gc_k(\Theta) } \frac{\sup_{\myf\in \myF} |G \myf-H \myf|}{W_{2k-1}^{2k-1}(G,H)} >0. \label{eqn:globalinversebound}
\end{equation}
It is clear that $\Gc_k(\Theta)$ is distinguishable by $\myF$ is a necessary condition for  \eqref{eqn:globalinversebound} to hold. To establish a uniform convergence rate around a neighborhood of some $G_0\in \Ec_{k_0}(\Theta)$, we need a local version of \eqref{eqn:globalinversebound}.  The local inverse bound relating $\sup_{\myf\in \myF} |G \myf-H \myf|$ to a Wasserstein distance is: 
\begin{equation}
\liminf_{ \substack{G,H\overset{W_1}{\to} G_0\\ G\neq H\in \Gc_k(\Theta) }} \frac{\sup_{\myf\in \myF} |G \myf-H \myf|}{W_{2d_1-1}^{2d_1-1}(G,H)} >0. \label{eqn:localinversebound}
\end{equation}
In the above inequality $d_1$ is a function of $G_0$: each $G_0$ has a unique number of atoms $k_0$, and thus has a unique overfit index $d_1=k-k_0+1$. The local inverse bound \eqref{eqn:localinversebound} and the global inverse bound \eqref{eqn:globalinversebound} are related by the following lemma.

\begin{lem} \label{lem:localtoglobal}
Suppose that $\Theta$ is compact. 
	If  \eqref{eqn:localinversebound} holds for any $G_0\in \Gc_k(\Theta)$ and  $\Gc_k(\Theta)$ is distinguishable by $\myF$, then  \eqref{eqn:globalinversebound} holds.
\end{lem}

\myred{
The following lemma states some equivalent formulations of inverse bounds. 

\begin{lem}[Equivalent versions of inverse bounds]
\label{lem:equinvbou}
 \begin{enumerate}[label=(\alph*)]
\item \label{lem:equinvboua}
\eqref{eqn:globalinversebound} is equivalent to
\begin{align*}
	 W_{2k-1}^{2k-1}(G,H) \leq C' \sup_{\myf\in \myF} |G \myf-H \myf|, \quad \forall G, H\in \Gc_k(\Theta)
\end{align*}
for some constant $C'$ (that possibly depends on the model). 

\item  \label{lem:equinvboub}
 Fix $G_0\in \Ec_{k_0}(\Theta)$. \eqref{eqn:localinversebound} is equivalent to the following:
 there exist $r(G_0)$ and  $C(G_0)$, where their dependence on $\myF, \Theta, k_0,k$ are suppressed, such that for any $G,H\in \Gc_k(\Theta)$ satisfying $W_1(G_0,G)<r(G_0)$ and $W_1(G_0,H)<r(G_0)$, 
\begin{align*}
 W_{2d_1-1}^{2d_1-1}(G,H)  
\leq &  C(G_0)   \sup_{\myf\in \myF} |G \myf-H \myf|.
\end{align*}

\item \label{lem:equinvbouc}
Suppose that $\Theta$ is compact and that $\Gc_k(\Theta)$ is distinguishable by $\myF$.  Fix $G_0\in \Ec_{k_0}(\Theta)$.  \eqref{eqn:localinversebound} is equivalent to the following: there exist $r(G_0)$ and  $C(G_0)$, where their dependence on $\myF, \Theta, k_0,k$ are suppressed, such that for any $ G, H\in \Gc_k(\Theta)$ satisfying $W_1(G_0,H)<r(G_0)$, 
\begin{align*}
 W_{2d_1-1}^{2d_1-1}(G,H)
\leq &  C(G_0)     \sup_{\myf\in \myF} |G \myf-H \myf|.
\end{align*}
\end{enumerate}
\end{lem}
}

To appreciate the fundamental roles of the inverse bounds in our framework, we state the following result on uniform convergence rates of any estimators in the following lemma. \myred{The proof is straightforward based on Lemma \ref{lem:equinvbou} and thus is omitted.}

\myred{
\begin{lem}[Consequences of inverse bounds] \label{lem:convergencerate}
Suppose that $\Theta$ is compact. Let $\hat{G}_n\in \Ef_n$ be any estimator.
 \begin{enumerate}[label=(\alph*)]
\item \label{itema:lem:convergencerate}
Suppose that \eqref{eqn:globalinversebound} holds. Then for any $G\in G_k(\Theta)$, any $t>0$
\begin{align}
	 \left\{W_{2k-1}^{2k-1}(G,\hat{G}_n)\geq t \right\} \cap \{ \hat{G}_n \in \Gc_k(\Theta) \} 
	\subset\left\{ \sup_{\myf\in \myF} |G \myf-\hat{G}_n \myf|\geq C(\myF,\Theta,k) t \right\} \cap \{ \hat{G}_n \in \Gc_k(\Theta) \} . \label{eqn:globalconvergencerategenerallem}
\end{align}

\item \label{itemb:lem:convergencerate}
Fix $G_0\in \Ec_{k_0}(\Theta)$. Suppose that \eqref{eqn:localinversebound} holds and that $\Gc_k(\Theta)$ is distinguishable by $\myF$. Then there exist $r(G_0)$ and  $C(G_0)$, where their dependence on $\myF, \Theta, k_0,k$ are suppressed, such that for any $G\in \Gc_k(\Theta)$ satisfying $W_1(G_0,G)<r(G_0)$, any $t >0$,
\begin{align}
\left\{ W_{2d_1-1}^{2d_1-1}(G,\hat{G}_n)\geq t \right\} \cap \{ \hat{G}_n \in \Gc_k(\Theta) \} 
\subset &    \left\{ \sup_{\myf\in \myF} |G \myf-\hat{G}_n \myf|\geq C(G_0)  t  \right\} \cap \{ \hat{G}_n \in \Gc_k(\Theta) \} . \label{eqn:convergencerategenerallem}
\end{align}
\end{enumerate}
\end{lem}
}
\begin{rem}
We provide several interpretations for \eqref{eqn:globalconvergencerategenerallem}, and omit the interpretations for
\eqref{eqn:convergencerategenerallem} due to the similarity.
Equation \eqref{eqn:globalconvergencerategenerallem} states that the event on $(\Xf,\Xc)$ defined in terms of a Wasserstein distance between mixing measures is a subset of the event defined in terms of $\sup_{\myf\in \myF} |G \myf-\hat{G}_n \myf|$. To quantify the convergence rate in the Wasserstein distance, it suffices to have a control on $\sup_{\myf\in \myF} |G \myf-\hat{G}_n \myf|$ for an estimator $\hat{G}_n$. 

Since \eqref{eqn:globalconvergencerategenerallem} is a relationship of events on the probability space $(\Xf,\Xc)$, we may evaluate the events under any probability measure to obtain a upper bound of the tail probability $\Pb\left(W^{2k-1}_{2k-1}(G,\hat{G}_n)\geq t\right)$. In this paper, the natural probability measure is the true $\Pb_{G^*}$ under which one obtain the observed i.i.d. samples $X_1,\ldots,X_n\overset{\iidtext}{\sim} \Pb_{G^*}$.
Another example for the probability measure, not covered in this work, is a posterior distribution $\Pi(\cdot|X_1,\cdots,X_n)$, which is derived via Bayes' formula from some prior distribution $\Pi$ 
on the space of mixing measures (see \cite{nguyen2013convergence,wei2022convergence}).   

It is noted that Lemma \ref{lem:convergencerate} is quite general and that it does not require  $G^*\in \Gc_k(\Theta)$. Note also that $G$ can be any mixing measure in $\Gc_k(\Theta)$, not necessarily the true mixing measure $G^*$. Thus, in the scenario of model misspecification, i.e., when $G^*\not \in \Gc_k(\Theta) $, one may take $G\in \argmin_{G'\in \Gc_k(\Theta)} \sup_{\myf\in \myF}|G'\myf - G^*\myf|$ be a projection of $G^*$ onto $\Gc_k(\Theta)$. We leave such directions to interested readers (see also \cite{guha2021posterior}). In the remainder of the paper we will work with the well-specified setting, i.e., $G^*\in \Gc_k(\Theta)$.
\myeoe
\end{rem}

\subsection{Minimum $\myF$-distance estimators and uniform convergence rates}
\label{sec:mindisest}
We now present a general estimator called the minimum $\myF$-distance estimator, which controls  $\sup_{\myf\in \myF} |G \myf-\hat{G}_n \myf|$.
In Section \ref{sec:invbou} we have stated that ``$\myF$ is a collection of functions such that for each $\myf\in \myF$, $G\myf$ can be relatively easy to estimate based on data samples from $\Pb_G$''. The precise assumption is as follows. 

\begin{definition}
\label{def:estimatable}
The family $\myF$ is said to be \emph{estimatable} on $\Gc_{k}(\Theta)$ if
for each $\myf\in \myF$, there exists a measurable function $t_\phi$ defined on $\Xf$ such that $G \myf=\Eb_{G}t_\myf(X_1)$ for any $G\in \Gc_k(\Theta)$. In other words, $t_\myf(X_1)$ is an unbiased estimate for $G\myf$.    
\end{definition}

If $\myF$ is estimatable 
then $G^* \myf = \Eb_{G^*} t_\myf(X_1)$, a quantity which may be estimated by its empirical analog $\frac{1}{n} \sum_{i\in [n]} t_\myf(X_i) $. We say the finite mixture model $\Pb_G$ is \emph{identifiable} on $G_k(\Theta)$ if $G\mapsto \Pb_{G}$ is injective on $\Gc_k(\Theta)$. The next lemma is a straightforward result connecting several definitions, the proof is omitted.

\begin{lem}
    If $\myF$ is estimatable on $\Gc_{k}(\Theta)$ and $\Gc_k(\Theta)$ is distinguishable by $\myF$, then the mixture model is identifiable on $\Gc_k(\Theta)$. \myeoe
\end{lem}

Suppose that $\Theta$ is compact and suppose that $\myF$ is estimatable on $\Gc_{k}(\Theta)$.
Define
$$
\hat{G}_n(\ell)\in \argmin_{G'\in \Gc_\ell(\Theta)} \sup_{\myf\in\myF}\left|G' \myf-\frac{1}{n}\sum_{i\in [n]}t_\myf(X_i)\right|, \quad \forall \ell\in \Nb_+. 
$$ 
\begin{figure}[ht!]
        \labellist
        \small \hair 2pt
        \pinlabel {$\Pb_n$} at 290 494 
        \pinlabel {$\Pb_{\hat{G}_n}$} at 235 277 
        \pinlabel {$\Pb_{G^*}$} at 455 217 
        \pinlabel {$n^{-\frac{1}{2}}$} at 425 397 
        \pinlabel {$\{\Pb_{G}: G\in \Gc_k(\Theta)\}$} at 117 127 
        \endlabellist
        \centering
        \includegraphics[width=0.3\linewidth]{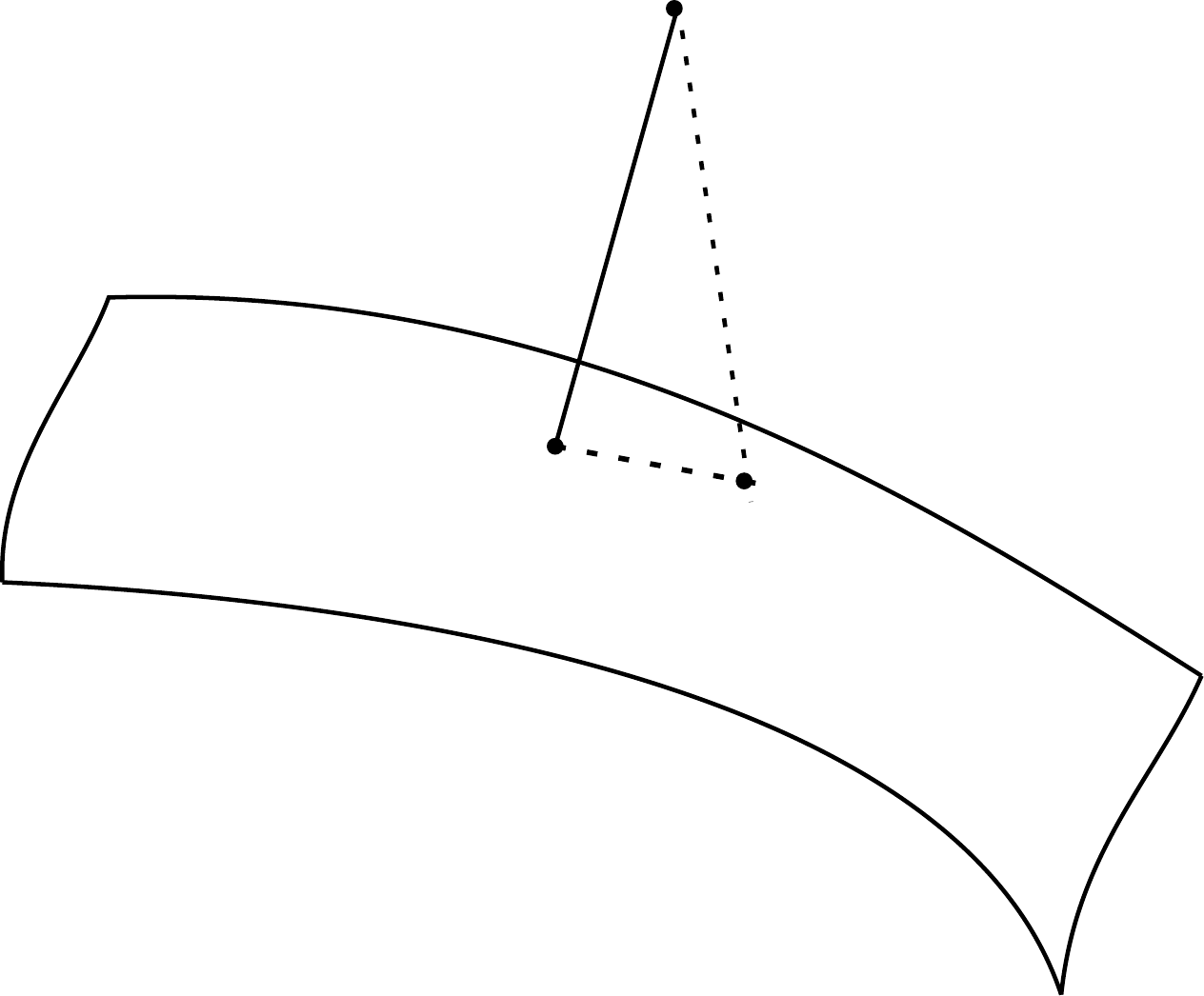}
        \caption{\myred{Minimum distance estimators: The set 
        $\Ac:=\{\Pb_{G}: G\in \Gc_k(\Theta)\}$, depicted by the surface in the plot, is the space of all mixture probability distributions. $\Pb_{G^*}$ is the true mixture distribution, which is an element on $\Ac$. Denote $\Pb_n$ to be the empirical measure based on $X_1,\ldots,X_n\sim \Pb_{G^*}$.  $\Pb_n$ is typically not on $\Ac$. Here we project $\Pb_n$ to $\Ac$ by finding the element $\Pb_{G'}$ that has the smallest distance to $\Pb_n$, where $\sup_{\myf\in\myF}\left|G' \myf-\frac{1}{n}\sum_{i\in [n]}t_\myf(X_i)\right|$ is a distance between $\Pb_{G'}$ to $\Pb_n$ (see Section \ref{sec:integralprobabilitymetric} ahead for more details).}   }
        \label{fig:exaofdef}
    \end{figure}
 Note that $\hat{G}_n(\ell)$ is well-defined since $\sup_{\myf\in\myF}\left|G' \myf-\frac{1}{n}\sum_{i\in [n]}t_\myf(X_i)\right|$ is non-negative 
 and lower semicontinuous w.r.t. $G'$, which implies that its minimum is attained on the compact space $\Gc_\ell(\Theta)$. Our estimator would be $\hat{G}_n=\hat{G}_n(k)$, which is termed a \emph{minimum $\myF$-distance estimator}. 

 Intuitively, since $\myF$ is estimatable on $\Gc_k(\Theta)$,  when $n$ is large, for any $G' \in \Gc_l(\Theta)$, one expects
 $$
 \sup_{\myf\in\myF}\left|G' \myf-\frac{1}{n}\sum_{i\in [n]}t_\myf(X_i)\right| \approx \sup_{\myf\in\myF}\left|G' \myf-G^* \myf\right|,
 $$
 and hence if $\Gc_k(\Theta)$ is distinguishable by $\myF$, one expects $\hat{G}_n$ to be close to $G^*$.  A summary of the estimation procedure is stated in Algorithm \ref{alg:unified}.

\RestyleAlgo{ruled}
\begin{algorithm}
\caption{Minimum $\myF$-distance estimators} \label{alg:unified}
\KwData{$X_1,\ldots, X_n\overset{\iidtext}{\sim} \Pb_{G^*}$}
\KwResult{$\hat{G}_n$}
$\bar{t}_\myf \gets \frac{1}{n}\sum_{i\in [n]}t_\myf(X_i)$, for each $\myf\in \myF$\;
$\hat{G}_n\in \argmin_{G'\in \Gc_k(\Theta)} \sup_{\myf\in\myF}\left|G' \myf-\bar{t}_\myf\right|$
\end{algorithm}

It follows that for any minimum $\myF$-distance estimator $\hat{G}_n$ and for any $G\in \Gc_k(\Theta)$, by the triangle inequality, 
\begin{equation}
\sup_{\myf\in \myF} |\hat{G}_n \myf- G \myf |\leq 2\sup_{\myf\in \myF} \left|\frac{1}{n}\sum_{i\in [n]}t_\myf(X_i)- G \myf \right|. \label{eqn:Gnproperty}
\end{equation}
\myred{As we shall see in some specific instances in Section \ref{sec:examples}, the $\myF$ distance can often serve as distance between mixture densities. Thus the above can be seen as density convergence rates are upper bounded by the supremum of some empirical process.}

By Lemma \ref{lem:convergencerate}\myred{, which lower bounds mixture densities distances by Wasserstein distance or moment differences between mixing distributions,} with $G=G^*$ and \eqref{eqn:Gnproperty} we immediately have the following uniform convergence rates \myred{of the mixing measures.} 

\begin{thm}[Uniform convergence rate] \label{thm:convergencerate}
Suppose that $\Theta$ is compact and suppose that $\myF$ is estimatable on $\Gc_{k}(\Theta)$. Let $\hat{G}_n$ be a minimum $\myF$-distance estimator. 
 \begin{enumerate}[label=(\alph*)]
\item 
Suppose that \eqref{eqn:globalinversebound} holds. Then there is a positive constant  $C$, where its dependence on $ \Theta, k, \myF$ and the probability kernel $\{\Pb_\theta\}$ is suppressed, such that for any $G^*\in \Gc_k(\Theta)$, any $t>0$, \myred{and for any $D\in \{W_{2k-1}^{2k-1}, \mbf_{2k-1} \} $, }\footnote{Throughout this paper if any of quantity is not measurable, the probability and the expectation should be understood as outer probability and outer expectation \cite[Section 1.2]{wellner2013weak}.}
\begin{align}
	\Pb_{G^*} \left(D(G^*,\hat{G}_n)\geq t \right) 
	\leq  \Pb_{G^*} \left( \sup_{\myf\in \myF} \left|\frac{1}{n}\sum_{i\in [n]}t_\myf(X_i)- G^* \myf \right| \geq C  t  \right), \label{eqn:globalconvergencerategeneralpro}
\end{align} 
and
\begin{align*}
	\Eb_{G^*} D(G^*,\hat{G}_n)
	\leq C \, \Eb_{G^*}  \sup_{\myf\in \myF} \left|\frac{1}{n}\sum_{i\in [n]}t_\myf(X_i)- G^* \myf \right| . 
\end{align*}

\item 
Fix $G_0\in \Ec_{k_0}(\Theta)$ for some $k_0\in [k]$. Suppose that \eqref{eqn:localinversebound} holds and that $\Gc_k(\Theta)$ is distinguishable by $\myF$. Then there exists $r(G_0)$, $C(G_0)$ and $c(G_0)$, where their dependence on $ \Theta, k_0,k,\myF$ and the probability kernel $\{\Pb_\theta\}$ are suppressed, such that for any $G^*\in \Gc_k(\Theta)$ satisfying $W_1(G_0,{G^*})<r(G_0)$, \myred{and for any $D\in \{W_{2d_1-1}^{2d_1-1}, \mbf_{2d_1-1} \} $, }
\begin{align}
\Pb_{G^*} \left(D(G^*,\hat{G}_n)\geq t \right) 
\leq   \Pb_{G^*} \left( \sup_{\myf\in \myF} \left|\frac{1}{n}\sum_{i\in [n]}t_\myf(X_i)- G^* \myf \right| \geq C(G_0)  t  \right), \label{eqn:convergencerategeneral1}
\end{align}
and
\begin{align*}
\Eb_{G^*} D(\hat{G}_n,G^*) 
\leq  C(G_0)  \Eb_{G^*} \sup_{\myf\in \myF} \left|\frac{1}{n}\sum_{i\in [n]}t_\myf(X_i)- G^* \myf \right|. 
\end{align*}

\end{enumerate}
\end{thm}

\begin{rem}
\label{rem:2.15}
%
By checking the proof \myred{in the appendix}, the property that $\myF$ is estimatable on $\Gc_k(\Theta)$ is actually not used to develop the above theorem; indeed any class of functions $\{\tilde{t}_\myf\}_{\myf\in \myF}$ on $\Xf$, not necessarily the particular class $\{t_\myf\}_{\myf\in \myF}$ in the Definition \ref{def:estimatable}, can be used in the definition of minimum $\myF$-distance estimators and does not change the conclusions of the above theorem. It is the performance of the estimator $\hat{G}_n$, which is governed by $\sup_{\myf\in \myF} \left|\frac{1}{n}\sum_{i\in [n]}\tilde{t}_\myf(X_i)- G^* \myf \right|$, that is affected by the choice of $\{\tilde{t}_\myf\}_{\myf\in \myF}$. Indeed, if we choose the particular class $\{t_\myf\}_{\myf\in \myF}$ in Definition \ref{def:estimatable}, in lieu of the intuition discussed after the definition of minimum $\myF$-distance estimators, we expect that $\Eb_{G^*}\sup_{\myf\in \myF} \left|\frac{1}{n}\sum_{i\in [n]}t_\myf(X_i)- G^* \myf \right|$ to be small; in fact if $\myF$ is estimatable on $\Gc_k(\Theta)$ it is the suprema of an empirical process: 
$$
\sup_{\myf\in \myF} \left|\frac{1}{n}\sum_{i\in [n]}t_\myf(X_i)- G^* \myf \right|=\sup_{\myf\in \myF} \left|\frac{1}{n}\sum_{i\in [n]}\left(t_\myf(X_i)- \Eb_{G^*} t_\myf(X_i)\right) \right|.
$$ 
Once we specialize the general framework to a specific example where $\myF$ and $t_\myf$ are concrete, we can say more about the empirical process and hence obtain a concrete convergence rate in terms of $n$ for the right hand sides of  \eqref{eqn:globalconvergencerategeneralpro} and \eqref{eqn:convergencerategeneral1}. Such examples will be provided in Section \ref{sec:examples}. 

Theorem \ref{thm:convergencerate} demonstrates a trade-off in choosing the (estimatable) function class $\myF$. On the one hand, $\myF$ has to be rich enough so that the inverse bounds \eqref{eqn:globalinversebound} and \eqref{eqn:localinversebound} hold. On the other hand, $\myF$ has to be small enough so that the governing empirical process $\sup_{\myf\in \myF} \left|\frac{1}{n}\sum_{i\in [n]}t_\myf(X_i)- G^* \myf \right|$ is well behaved to obtain sharp rates of convergence.
\myeoe
\end{rem}

\myred{
\begin{rem}
The proof in the appendix is only presented for $D$ to be Wasserstein distance. But in view of Theorem \ref{thm:inverseboundmoment} ahead, the proof can be trivially adapted to moment difference.  \myeoe
\end{rem}
}


\subsection{Sufficient identification conditions for local inverse bounds}


At the core of our methodological and theoretical framework is the precise connection between the choice of function class $\Phi$ and the convergence rates for the mixing distribution that this choice affects. In particular, $\Phi$ needs to be sufficiently rich so that the inverse bounds hold.
To apply Theorem \ref{thm:convergencerate}, we need to verify that \eqref{eqn:localinversebound} holds for any $G_0\in \Gc_k(\Theta)$, in lieu of  Lemma \ref{lem:localtoglobal}. In this subsection we provide sufficient conditions to establish that the local inverse bound \eqref{eqn:localinversebound} holds for any $G_0\in \Gc_k(\Theta)$.

\begin{definition} \label{def:linearindependentdomain}
The family $\myF$ is said to be a \textit{$(m,k_0,k)$ linear independent domain} if the following hold: 1) Each $\myf\in \myF$ is $m$-th order continuously differentiable on $\Theta$; \footnote{To make sense of the differentiability at the boundary of $\Theta$, it suffices to treat $\myf\in \myF$ as functions defined on a larger domain $\tilde{\Theta}$ and our prior parameter space $\Theta\subset \tilde{\Theta}^\circ$. 
} and 2) Consider any integer $\ell\in [k_0,2k-k_0]$, and any vector $(m_1,m_2,\ldots,m_\ell)$ such that $1\leq m_i\leq m+1$ for $i\in [\ell]$ and $\sum_{i=1}^\ell m_i\in [2k_0, 2k]$, then
for any distinct $\{\theta_i\}_{i\in [\ell] }\subset \Theta$, the  operators $\{D^\alpha|_{\theta=\theta_i}\}_{ 0\leq |\alpha|< m_i, i\in [\ell]}$  on $\myF$ are linear independent, i.e.,  
\begin{subequations}
\begin{align}
\sum_{i=1}^{\ell}\  \sum_{|\alpha|\leq m_i-1} a_{i\alpha} D^\alpha \myf(\theta_i)  = &0,  \quad \forall \myf\in \myF \label{eqn:linearinddomaina}\\
\sum_{i\in [\ell]}   a_{i\bm{0}}   = &0, \label{eqn:linearinddomainb}
\end{align}
\end{subequations}
if and only if 
$$
a_{i\alpha}=0, \quad \forall \ 0\leq |\alpha|< m_i, \ i\in [\ell].
$$
\end{definition}

\begin{rem} \label{rem:linearinddomainconstraint}
The equation \eqref{eqn:linearinddomainb} can be seen as \eqref{eqn:linearinddomaina} with $\myf \equiv 1_\Theta$, the constant function $1$ on $\Theta$. So a slightly more accurate terminology should be ``the  operators $\{D^\alpha|_{\theta=\theta_i}\}_{ 0\leq |\alpha|< m_i, i\in [\ell]}$  on $\myF\cup \{1_\Theta\}$ are linear independent''. 

It is clear that if a subset of $\myF$ is a $(m,k_0,k)$ linear independent domain, then so is $\myF$. Another observation is that if $\myF$ is a $(m,k_0,k)$ linear independent domain then $\myF$ is a $(m',k',k)$ linear independent domain for any $m'\leq m$ and $k'\geq k_0$.
\myeoe
\end{rem}

A related and somewhat more standard notion of strong identifiability has been widely studied in the previous work \cite{chen1995optimal,nguyen2013convergence,ho2016strong,heinrich2018strong,ho2019singularity}, which are roughly the linear independence between the mixture kernel density (or CDF) and its derivatives. In the next definition we generalize the concept to our general framework based on test function $\myF$, which will recover the existing definitions once a suitable $\myF$ is chosen. 

\begin{definition}[$m$-strong identifiability] 
\label{def:mstrong}
A family $\myF$ of functions of $\theta$ is $m$-strongly identifiable if each $\myf\in\myF$ is $m$-order continuously differentiable; and for any finite set of $\ell$ distinct points $\theta_i\in \Theta$, 
$$
\sum_{i=1}^{\ell}\  \sum_{|\alpha|\leq m} a_{i\alpha} D^\alpha \myf(\theta_i)  = 0,  \quad \forall \myf\in \myF
$$    
if and only if 
$$
a_{i\alpha}=0, \quad \forall \ 0\leq |\alpha|\leq m, \ i\in [\ell].
$$
\end{definition}
\myred{ 
\begin{rem}
\label{rem:lindomvsmstr}
The advantage of the definition of $m$-strong identifiability is that it is simpler and more straightforward to verify. But it is clear that $\myF$ is $m$-strongly identifiable implies that $\myF$ is a $(m,k_0,k)$ linear independent domain for any $k_0\leq k$. That $\myF$ is a $(m,k_0,k)$ linear independent domain is an improvement over $m$-strong identifiability due to the reduced number of equations required, ones that arise from a careful consideration of possible allocations of atoms of $k$-component mixing measures converging to a fixed $k_0$-component mixing measure. The relaxation from $m$-strong identifiability to our definition of linear independent domain while maintaining the guaranteed inverse bounds (Theorem \ref{thm:inversebound}) is one of the key contributions in this paper. In fact, when $\myF$ is of finite cardinality, the  linear system in the definition of linear independent domain is much better behaved than that in the definition of strong identifiability, since the former has less variables while the number of equations remain the same. In particular, there are some important examples, e.g., family of monomials $\myF_2$ (see Section \ref{sec:momentdeviation}), that are $(m,k_0,k)$ linear independent domain but not $m$-strongly identifiable. Moreover, we show in Example \ref{exa:mixmuldis} for mixture of multinomial distributions, by using our weaker condition of linear independent domain, the inverse bounds hold if and only if $N\geq 2k-1$ which improves the previous results \cite[Proposition 1 and Corollary 1]{manole2021estimating}.  \myeoe
\end{rem}
}

The following general theorem establishes that $(m,k_0,k)$ linear independent domains are sufficient conditions for establishing fundamental local inverse bounds.


\begin{thm} 
\label{thm:inversebound}
If $\Theta\subset \Rb^q$ is compact.
\begin{enumerate}[label=(\alph*)]

\item \label{thm:inversebounda}
If $\myF$ is a $(2d_1-1,k_0,k)$ linear independent domain, then \eqref{eqn:localinversebound} holds for any $G_0\in \Ec_{k_0}(\Theta)$.
\item \label{thm:inverseboundb}
If $\myF$ is a $(2k-1,1,k)$ linear independent domain, then \eqref{eqn:localinversebound} holds for any $G_0\in \Gc_{k_0}(\Theta)$ for any $k_0\in [k]$. 
\end{enumerate}
\end{thm}

By Remark \ref{rem:linearinddomainconstraint} if $\myF$ is a $(2k-1,1,k)$ linear independent domain this implies that $\myF$ is a $(2d_1-1,k_0,k)$ linear independent domain for any $k_0\in [k]$, hence in Theorem \ref{thm:inversebound} part \ref{thm:inverseboundb} immediately follows from part \ref{thm:inversebounda}.
It is also noted that
in Section \ref{sec:optimality} we show that the exponent $2d_1-1$ of the denominator in \eqref{eqn:localinversebound} is optimal. In general, the compactness assumption is necessary for inverse bounds to hold; see relevant discussions in Remark \ref{rem:compactness}. 
\begin{rem}
Notice that $\sup_{\myf\in \myF} |G \myf-H \myf| = \sup_{\myf\in \myF\cup \{1_\Theta\}} |G \myf-H \myf| $ since $G 1_\Theta - H 1_\Theta = 0$. So we may always assume that $1_\Theta\in \myF$ without affecting \eqref{eqn:localinversebound}. For $\myf=1_\Theta$, the corresponding $t_\myf=1_\Xf$. Hence the minimum $\myF$-distance estimator $\hat{G}_n$ also remains unchanged by replacing $\myF$ with $\myF\cup \{1_\Theta\}$. See Remark \ref{rem:linearinddomainconstraint} for a related discussion.
\myeoe
\end{rem}

\begin{rem}
The proof of Theorem \ref{thm:inversebound} part \ref{thm:inversebounda} follows a structure similar to that of the proof of \cite[Theorem 6.3]{heinrich2018strong}, which is based on their original construction of a coarse-grained tree for the space of supporting atoms. While \cite[Theorem 6.3]{heinrich2018strong} only deals with the special case that $\myF=\myF_0$ (cf. Example \ref{exa:divationbetweenmixturedistributions}), Theorem \ref{thm:inversebound} is more general and can be applied to other function class $\myF$. Even for the special case $\myF=\myF_0$, Theorem \ref{thm:inversebound} has several improvements: it applies to the multivariate distribution while \cite[Theorem 6.3]{heinrich2018strong} only considers the univariate mixture distribution; moreover, in the technical sense the assumptions needed are also relaxed. We will revisit case $\myF=\myF_0$ in Section \ref{sec:minimumdistance estimator} and discuss the comparisons in more detail.    \myeoe
\end{rem}

\myred{

\subsection{Information geometry of finite mixture models}
\label{sec:infgeo}
In this section we present an improvement to the inverse bounds obtained in Theorem \ref{thm:inversebound}, by relating the $\Phi$-distance of mixture densities to the moment difference distance of mixing measures. We establish not only an lower bound, but also an upper bound and therefore characterizes the local geometry of mixture distributions in terms of the moment difference. 

\begin{thm} 
\label{thm:inverseboundmoment}
If $\Theta\subset \Rb^q$ is compact.
\begin{enumerate}[label=(\alph*)]

\item \label{itema:thm:inverseboundmoment}
If $\myF$ is a $(2d_1-1,k_0,k)$ linear independent domain, then for any $G_0\in \Ec_{k_0}(\Theta)$, it holds:
\begin{equation}
\liminf_{ \substack{G,H\overset{W_1}{\to} G_0\\ G\neq H\in \Gc_k(\Theta) }} \frac{\sup_{\myf\in \myF} |G \myf-H \myf|}{\mbf_{2d_1-1}(G-\theta_0,H-\theta_0)} >0, \label{eqn:localinverseboundmoment}
\end{equation}
where $\theta_0$ is a arbitrary element in $\Rb^q$.
\item \label{itemb:thm:inverseboundmoment}
If $\myF$ is a $(2k-1,1,k)$ linear independent domain, then \eqref{eqn:localinverseboundmoment} holds for any $G_0\in \Gc_{k_0}(\Theta)$ for any $k_0\in [k]$. 
\item \label{itemc:thm:inverseboundmoment}
Suppose 
$$
\sup_{|\alpha|\leq 2d_1-1} \sup_{\theta\in \Theta} \sup_{\myf\in\myF} \left| D^\alpha \myf(\theta) \right|<\infty
$$
and that
there is a uniform continuity modulus $w(\cdot)$ such that: for any $\alpha$ with $|\alpha|=m$,  
$$\sup_{\myf\in \myF}|D^\alpha \myf(\theta) - D^\alpha \myf(\theta')|\leq w(\theta-\theta')$$ 
with $\lim_{h\to 0}  w(h) =0$. Then
\begin{equation}
\limsup_{ \substack{G,H\overset{W_1}{\to} G_0\\ G\neq H\in \Gc_k(\Theta) }} \frac{\sup_{\myf\in \myF} |G \myf-H \myf|}{\mbf_{2d_1-1}(G-\theta_0,H-\theta_0)} <\infty. 
\label{eqn:localinverseboundmomentupp}
\end{equation}
where $\theta_0$ is a arbitrary element in $\Rb^q$. 
\end{enumerate}
\end{thm}

The above theorem basically states that  $\sup_{\myf\in \myF} |G \myf-H \myf|$, the distance between mixture densities, which will be discussed in detailed in Section \ref{sec:examples}, is roughly the same as $\mbf_{2d_1-1}(G-\theta_0,H-\theta_0)$, the moment difference between the corresponding mixing measures. Specifically, given some regularity condition specified in parts \ref{itema:thm:inverseboundmoment} and \ref{itemc:thm:inverseboundmoment}, it then holds that in a small neighborhood around $G_0$: for some constant $c,C$,
\begin{equation}
c\mbf_{2d_1-1}(G-\theta_0,H-\theta_0)  \leq \sup_{\myf\in \myF} |G \myf-H \myf| \leq C \mbf_{2d_1-1}(G-\theta_0,H-\theta_0) \label{eqn:locinfgeo}
\end{equation}
for any $G$ and $H$ in the neighborhood. Recall $d_1=2(k-k_0+1)-1$ depends on the number of component $k_0$ of $G_0$. If more differentiability conditions as in part \ref{itemb:thm:inverseboundmoment} and $\Gc_k(\Theta)$ is distinguishable by $\myF$, then it also holds that 
\begin{equation}
c\mbf_{2k-1}(G-\theta_0,H-\theta_0)  \leq \sup_{\myf\in \myF} |G \myf-H \myf| \leq C \mbf_{2k-1}(G-\theta_0,H-\theta_0)  \label{eqn:glabalequimoment}
\end{equation}
for any $G,H\in \Gc_k(\Theta)$.

\begin{rem}
\myred{The inequalities in Eq.\eqref{eqn:glabalequimoment} (Eq. \eqref{eqn:locinfgeo}) demonstrate that the moment difference of mixing measures is an intrinsic metric that captures precisely the global (local) geometry of data population densities, under a (relaxed) condition of strongly identifiable finite mixtures,  
even though the
Wasserstein distance arguably continues to be an appealing surrogate for measuring the quality of the model's parameter estimates. }
In fact, by Lemma \ref{lem:momentbound} we know that the moment difference is an upper bound for Wasserstein distance between mixing measures, and consequently, upper bounds for estimation errors in terms of moment difference will automatically yield the same upper bounds for estimation errors in terms of Wasserstein distance. Moreover, since Theorem \ref{thm:inverseboundmoment} provides an inverse bound in terms of the moment difference, as we have done in Theorem \ref{thm:convergencerate} we can restate all the remaining examples in Section \ref{sec:examples} in terms of estimation errors defined via moment difference, which are strictly speaking stronger results as discussed in the previous sentence. 
\myeoe
\end{rem}

\begin{rem}
\label{rem:infgeocom}
    There is some related work on quantifying the distances between mixing measures. For location Gaussian mixtures, \cite[Theorem 4.2]{doss2020optimal} establishes a similar inequality \textcolor{blue}{to \eqref{eqn:glabalequimoment}} for squared Hellinger distance, KL divergence, and $\chi^2$-divergence. It is noteworthy that their results, by leveraging special properties of the Gaussian distribution, quantify how the constant coefficients depend on $q$ and $k$. However, it is not straightforward to see how their proof can be generalized beyond location Gaussian mixtures. Another related work is \cite{gassiat2014local} where they quantify the local geometry with respect to Hellinger distance of general location mixtures by a pseudo metric. 
The strength of Theorem \ref{thm:inverseboundmoment} is that it is applicable to general mixture models and a general integral probability metric, which extend beyond either location Gaussian mixtures or location mixtures considered in such prior work. \myred{For instance, one can easily obtain sufficient conditions for general mixture models to obtain Eq.\eqref{eqn:glabalequimoment} (Eq. \eqref{eqn:locinfgeo}) type inequality with $\sup_{\myf\in \myF} |G \myf-H \myf|$ specialized to KS distance (see Section \ref{sec:minimumdistance estimator}), MMD (see Section \ref{sec:MMD}), total variational distance (see Section \ref{sec:beyondsup}), and beyond; we leave the details to interested readers. Furthermore, the local type inequality \eqref{eqn:locinfgeo} appears to be novel, to the best of our knowledge. }
\myeoe
\end{rem}

}

\section{Instances of the minimum $\myF$-distance estimators}
\label{sec:examples}
In this section we shall see that specializing the function class $\Phi$ leads to existing estimation methods as well as new ones.

\subsection{Minimum IPM estimators}
\label{sec:integralprobabilitymetric}

First, we specialize the general results in Section \ref{sec:unifiedframework} to the case where $\myF$ takes the form in Example \ref{exa:divationbetweenmixturedistributions}. Consider $\myF=\left\{ \theta \mapsto \int f_1(x) \Pb_\theta(dx) \mid f_1 \in \Fc_1 \right\}$ where $\Fc_1$ is some subset of $\Mc$, the space of all measurable functions on $(\Xf,\Xc)$. 
Then
\begin{equation*}
\sup_{\myf\in \myF} |G \myf-H \myf| = \sup_{f_1 \in \Fc_1} \left|\int f_1 d\Pb_G-\int f_1 d\Pb_H\right|.
\end{equation*}
\myred{See Table \ref{tab:ipm_examples} for a list of function classes $\Fc_1$ and the associated IPM distances. }

\begin{table}[h!]
\centering
\begin{tabular}{|c|c|}
\hline
\textbf{Function Class} $\boldsymbol{\mathcal{F}_1}$ & \textbf{Distance Metric} \\
\hline
$\left\{1_{(-\infty,a]}(x)\,:\, a \in \mathbb{R} \right\}$ & Kolmogorov–Smirnov distance \\
\hline
$\left\{f : \|f\|_\infty \leq 1 \right\}$ & Total variation distance \\
\hline
$\left\{f : \|f\|_L \leq 1 \right\}$ & Wasserstein-1 distance \\
\hline
\{f : \text{unit ball in RKHS} \} & Maximum mean discrepancy (MMD) \\
\hline
$\left\{f : \|f\|_L + \|f\|_\infty \leq 1 \right\}$ & Dudley's metric \\
\hline
$\ldots$ & $\ldots$ \\
\hline
\end{tabular}
\caption{Examples of function classes $\mathcal{F}_1$ and their associated IPM distances}
\label{tab:ipm_examples}
\end{table}

Such a function class $\myF$ is automatically estimatable, since for each $f_1\in \Fc_1$, or equivalently for each $\myf\in \myF$, there exists a function 
$
t_\phi=f_1
$ 
defined on $\Xf$ such that $ \Eb_{G}t_\myf(X) = \int f_1 d\Pb_G= \int \int f_1 d\Pb_\theta dG = G \myf  $ holds for any probability measure $G$, including $G\in \Gc_k(\Theta)$. 
In the remainder of this subsection we will take $t_\myf$ in the form presented in the previous sentence.
In this case the corresponding minimum $\myF$-distance estimator becomes 
\begin{align*}
\hat{G}_n\in & \argmin_{G'\in \Gc_k(\Theta)} \sup_{f_1 \in \Fc_1}\left|\int f_1 d\Pb_{G'}-\frac{1}{n}\sum_{i\in [n]}f_1(X_i)\right| \\
= & \argmin_{G'\in \Gc_k(\Theta)} \sup_{f_1 \in \Fc_1}\left|\int f_1 d\Pb_{G'}-\int f_1 d\hat{\Pb}_n\right|,
\end{align*}
where $\hat{\Pb}_n:=\frac{1}{n}\sum_{i\in [n]}\delta_{X_i}$ denotes the empirical measures. We refer to the minimum $\myF$-distance estimators in this case as \emph{minimum IPM estimators}.

Theorem \ref{thm:convergencerate} can be applied in this case, with the governing empirical process 
\begin{equation}
\sup_{\myf\in \myF} \left|\frac{1}{n}\sum_{i\in [n]}t_\myf(X_i)- G^* \myf \right|=\sup_{f_1 \in \Fc_1}\left|\frac{1}{n}\sum_{i\in [n]}f_1(X_i)-\int f_1 d\Pb_{G^*}\right|.  \label{eqn:empiricalprocessIPM}
\end{equation}


\myred{
In fact, one may view minimum $\Phi$ distance estimators as minimum IPM estimators. 
\begin{rem}[minimum $\Phi$-distance estimators are also minimum IPM estimators]
\label{rem:GMMisIPM}
Although we presented  minimum IPM estimators as a specific instance in the bigger category of minimum $\myF$-distance estimator, note that they are actually equivalent: if $\myF$ is estimatable, then \eqref{eqn:momentmethod} can be written as
$$
\hat{G}_n\in \argmin_{G'\in \Gc_k(\Theta)} \sup_{\myf\in\myF}\left|G' \myf-\frac{1}{n}\sum_{i\in [n]}t_\myf(X_i)\right|= \argmin_{G'\in \Gc_k(\Theta)} \sup_{\myf\in\myF}\left|\int t_\myf(x)d\Pb_{G'}-\int t_\myf(x)d\Pb_n\right|. 
$$
which is a minimum IPM estimator with $\Fc_1= \{t_\myf: \myf\in \myF \}$. \myeoe 
\end{rem}
}

For minimum IPM estimators,  the freedom lies in the class $\Fc_1$ of functions on $\Xf$. Two specific choices of $\Fc_1$ with  concrete uniform convergence rates are provided in the following.

\subsubsection{Minimum KS-distance estimators}

\label{sec:minimumdistance estimator}

In this subsection consider $(\Xf,\Xc)=(\Rb^d,\Bc(\Rb^d))$. Take $\Fc_1 := \left\{1_{(-\infty,x]}(\cdot) \mid x\in \Rb^d \right\}$ where $(-\infty,x]:=(-\infty,x^{(1)}]\times \ldots \times(-\infty,x^{(d)}]$. 
For such $\Fc_1$, we have  $\myF = \myF_0 = \{ \theta \mapsto F(x \mid \theta) \mid x\in \Rb^d  \}$, where $F(x \mid \theta)$ as a function of $x\in \Rb^d$ is the CDF corresponding to $\Pb_\theta$. This has been presented briefly in Example \ref{exa:divationbetweenmixturedistributions} with $d=1$ and we now elaborate it more in this subsection. For $\myF=\myF_0$, 
$$
\sup_{\myf\in \myF_0} |G \myf-H \myf|=\sup_{x \in \Rb^d} \left|\int_{(-\infty,x]}  d\Pb_G-\int_{(-\infty,x]} d\Pb_H\right|  =: \KS(\Pb_G,\Pb_H), 
$$
where $\KS(\cdot,\cdot)$ stands for Kolmogorov–Smirnov distance. 
Let $F_n(x):= \frac{1}{n}\sum_{i\in [n]}1_{(-\infty,x]}(X_i)$ be the empirical CDF and  $\hat{\Pb}_n=\frac{1}{n}\sum_{i\in[n]}\delta_{X_i}$ the empirical measure. Then a minimum IPM estimator is $\hat{G}_n \in \argmin_{G'\in \Gc_k}\KS(\Pb_{G'},\hat{\Pb}_n)$, which is historically known as a minimum distance estimator \cite{deely1968construction, chen1995optimal, heinrich2018strong}. To avoid any confusion under our general framework, we shall call this a \emph{minimum KS-distance estimator}.
\begin{algorithm}
\caption{Minimum KS-distance estimators} 
\KwData{$X_1,\ldots, X_n\overset{\iidtext}{\sim} \Pb_{G^*}$}
\KwResult{$\hat{G}_n$}
$\hat{F}_n(x) \gets \frac{1}{n}\sum_{i\in [n]} 1_{(-\infty,x)}(X_i)$, for each $x\in \Rb^d$\;
$\hat{G}_n\in \argmin_{G'\in \Gc_k(\Theta)} \sup_{x\in\Rb^d}\left| \int_{\Theta} F(x|\theta) dG' -\hat{F}_n(x) \right|$
\end{algorithm}

We now refine the assumptions in the general framework to the specific case $\myF_0$. Firstly, it is clear that $\Gc_k(\Theta)$ is distinguishable by $\myF=\myF_0$ if and only if the mixture model is identifiable on $G_k(\Theta)$. Secondly, the family $\myF_0$ is a $(m,k_0,k)$ linear independent domain if the following two conditions hold: 1) For each $x \in \Rb^d$, $F(x|\theta)$ is $m$-th order continuously differentiable on $\Theta$;  2) Consider any integer $\ell\in [k_0,2k-k_0]$, and any vector $(m_1,m_2,\ldots,m_\ell)$ such that $1\leq m_i\leq m+1$ for $i\in [\ell]$ and $\sum_{i=1}^\ell m_i\in [2k_0, 2k]$.
For any distinct $\{\theta_i\}_{i\in [\ell] }\subset \Theta$, the functions (as functions of $x\in \Rb^d$) $\{D^\alpha F  (x \mid \theta_i)\}_{ 0\leq |\alpha|< m_i, i\in [\ell]}$  are linear independent, i.e.,  
\begin{subequations}
\begin{align}
\sum_{i=1}^{\ell}\  \sum_{|\alpha|\leq m_i-1} a_{i\alpha} D^\alpha F(x \mid \theta_i)  = &0,  \quad \forall x\in \Rb^d \label{eqn:linearinddomainCDFa}\\
\sum_{i\in [\ell]}   a_{i\bm{0}}   = &0, \label{eqn:linearinddomainCDFb}
\end{align}
\end{subequations}
if and only if 
$$
a_{i\alpha}=0, \quad \forall \ 0\leq |\alpha|< m_i, \ i\in [\ell].
$$
The condition that $\myF_0$ is a $(m,k_0,k)$ linear independent domain is indeed a condition on the probability kernels $\Pb_\theta$ or its corresponding CDF $F(x \mid \theta)$. 
One can similarly specialize the definition of $m$-strong identifiability in this case, which is the same as \cite[Definition 2.2]{heinrich2018strong} (with an additional equality constraint \eqref{eqn:linearinddomainCDFb}). Note that in this case, $0$-strong identifiability implies the mixture model is identifiable; that is why the definition is termed strongly identifiable even though the definition is on linear independence between the functions and their derivatives. The strong identifiability condition \cite[Definition 2.2]{heinrich2018strong} is a stronger assumption in the sense that it implies that $\myF_0$ is a $(m,k_0,k)$ linear independent domain and the mixture model is identifiable. 

  To obtain a concrete convergence rate for the minimum KS-distance estimators by applying Theorem \ref{thm:convergencerate}, it remains to control the governing empirical process from \eqref{eqn:empiricalprocessIPM}:
 \begin{equation}
\sup_{f_1 \in \Fc_1}\left|\int f_1 d\Pb_{G^*}-\frac{1}{n}\sum_{i\in [n]}f_1(X_i)\right| = \KS(\Pb_{G^*},\hat{\Pb}_n) .  \label{eqn:empiricalprocessKS}
\end{equation}
The next lemma is such a result for any probability measure $\Pb$, not necessarily the mixture probability measures $\Pb_G$.

 \begin{lem}
 \label{lem:KSempiricalprocess}
    Let $X_1,X_2,\ldots,X_n$ be i.i.d. samples from the probability measure $\Pb$ on $\Rb^d$. Then
    $$
    \Eb \KS(\hat{\Pb}_n,\Pb) \leq C \sqrt{\frac{d}{n}},
    $$
    where $C$ is independent of $\Pb$.
\end{lem} 
\begin{proof}
    It follows directly from \cite[Corollary 7.18]{van2014probability} or \cite[Theorem 8.3.26]{vershynin2018high}.
\end{proof}
 
 The next theorem is an immediate consequence of Theorem \ref{thm:convergencerate}  
 combined with Theorem \ref{thm:inversebound} and Lemma \ref{lem:KSempiricalprocess}. 

\begin{thm} \label{thm:convergencerateminimumdistanceestimator} 
Suppose that $\Theta$ is compact, and that the mixture model is identifiable on $\Gc_k(\Theta)$. Let $\hat{G}_n$ be a minimum KS-distance estimator.
 \begin{enumerate}[label=(\alph*)] 
\item \label{itema:thm:convergencerateminimumdistanceestimator}
Assume that $\myF_0$ is a $(2k-1,1,k)$ linear independent domain. There exists $C$ where its dependence on $ \Theta, k$ and the probability kernel $\{\Pb_\theta\}$ are suppressed, such that \myred{for $D\in \{W_{2k-1}^{2k-1}, \mbf_{2k-1} \} $},
 \begin{align*}
 \sup_{G^*\in \Gc_k(\Theta)}	\Eb_{G^*} D(G^*,\hat{G}_n )
 	\leq   C n^{-\frac{1}{2}}. 
 \end{align*}
\item \label{itemb:thm:convergencerateminimumdistanceestimator}
Fix $G_0\in \Ec_{k_0}(\Theta)$. Assume that $\myF_0$ is a $(2d_1-1,k_0,k)$ linear independent domain. Then there exist positive constants $r(G_0)$, $C(G_0)$ and $c(G_0)$, where their dependence on $ \Theta, k_0,k$ and the probability kernel $\{\Pb_\theta\}$ are suppressed, such that 
\myred{for $D\in \{W_{2d_1-1}^{2d_1-1}, \mbf_{2d_1-1} \}$},
     \begin{align*}
     \sup_{\substack{G^*\in \Gc_k(\Theta):W_1(G_0,G^*)<r(G_0) }}	\Eb_{G^*} D(G^*,\hat{G}_n )
     	\leq   c(G_0) n^{-\frac{1}{2}}. 
     \end{align*}
\end{enumerate}
\end{thm}

\begin{rem} \label{rem:minimaxrateKS}
\myred{Since $W_\ell(P,Q)$ is increasing in $\ell$, for $\ell\in [2k-1]$,
    \begin{multline}
     \Eb_{G^*} W_{1}(G^*,\hat{G}_n ) \leq \Eb_{G^*}W_\ell(G^*,\hat{G}_n ) \leq \Eb_{G^*} W_{2k-1}(G^*,\hat{G}_n ) \\ \leq \left(\Eb_{G^*} W^{2k-1}_{2k-1}(G^*,\hat{G}_n )\right)^{\frac{1}{2k-1}}\leq C_k \left(\Eb_{G^*}\mbf_{2k-1}(G^*,\hat{G}_n)\right)^{\frac{1}{2k-1}}, 
     \end{multline}
    where the last inequality is due to Jensen's inequality. The above inequality, Theorem \ref{thm:convergencerateminimumdistanceestimator} \ref{itema:thm:convergencerateminimumdistanceestimator} and Theorem \ref{thm:minimax} then imply that the minimax optimal rate for any $W_{\ell}(G^*,\hat{G}_n )$ and $\mbf_{2k-1}(G^*,\hat{G}_n)$ is $n^{-\frac{1}{2(2k-1)}}$ for any $\ell\in [2k-1]$ under the setting of Theorem \ref{thm:convergencerateminimumdistanceestimator} \ref{itema:thm:convergencerateminimumdistanceestimator}.} 
    By similar arguments we can also obtain the minimax optimal rate for any $W_{\ell}(G^*,\hat{G}_n )$ and $\mbf_{2d_1-1}(G^*,\hat{G}_n )$ is $n^{-\frac{1}{2(2d_1-1)}}$ for any $\ell\in [2d_1-1]$ under the setting of Theorem \ref{thm:convergencerateminimumdistanceestimator} \ref{itemb:thm:convergencerateminimumdistanceestimator}.  So the uniform convergence rate over the whole mixing measure space $\Gc_k(\Theta)$ for minimum KS-distance estimators is $n^{-\frac{1}{4k -2}}$. If some prior knowledge that $G^*$ is in a small neighborhood of a discrete distribution $G_0\in \Ec_{k_0}(\Theta)$, then the (now local) uniform convergence rate improves (decreases) to  $n^{-\frac{1}{4d_1 -2}}$ with $d_1=k-k_0+1$. 
    \myeoe
\end{rem}

\begin{rem}
Similar uniform convergence rates for the minimum KS-distance estimator in the case $d=1$ and $q=1$ were established by \cite[Theorem 3.3]{heinrich2018strong} who corrected an earlier study of \cite{chen1995optimal}. Here, we extend these results to any finite $d$ and any finite $q$ by specializing our general framework to the choice $\myF=\myF_0$. Note that even in the case $d=q=1$, Theorem \ref{thm:inversebound} improves \cite[Theorem 3.3]{heinrich2018strong} in the technical sense that less assumptions are imposed: our assumption that $\myF_0$ is a $(m,k_0,k)$ linear independent domain is weaker than the strong identifiable assumption \cite[Proposition 2.3]{heinrich2018strong} in that less differentiability and only a constrained linear independence are assumed; more importantly, our theorem does not require the uniform continuity modulus assumption \cite[Assumption B(k)]{heinrich2018strong}. 
\myeoe
\end{rem}

\subsubsection{Minimum MMD estimators}
\label{sec:MMD}

In this subsection we present another example of minimum IPM estimators,
which we call the minimum MMD estimator. MMD stands for maximum mean discrepancy, a metric that arises from a particular choice of $\Phi$ using reproducing kernel Hilbert spaces (RKHS) \cite{gretton2012kernel}. Unlike the minimum KS-distance estimators, the minimum MMD estimator seems novel in the literature of mixture models to the best of our knowledge. 
We emphasize that while the minimum KS estimator may be difficult to apply to non-Euclidean or high-dimensional or complex structured data domain $\Xf$ (since the Kolmogorov-Smirnov distance evaluation involves finding the supremum over $\Xf$), the minimum MMD may be more applicable in such settings thanks to the powerful machinery of the RKHS. As the minimum MMD estimators represent a novel instance of our general framework, they shall be treated in considerable detail. 

\paragraph{Maximum mean discrepancy}
 First, we recall some basic background of the RKHS and the associated MMD metric. In this section, $\Xf$ is assumed to be any topological space endowed with the $\sigma$-algebra $\Xc=\Bc(\Xf)$, the Borel measurable sets.  
 Consider a real-valued symmetric and positive semidefinite kernel function $\ker(\cdot,\cdot)$ on the measurable space $(\Xf,\Xc)$.
Let $\Hc$ denote the reproducing kernel Hilbert space (RKHS) associated with the reproducing kernel $\ker(\cdot,\cdot)$  with its inner product $\langle \cdot, \cdot \rangle_\Hc$, i.e., $\Hc$ is a Hilbert space of functions on $\Xf$ which satisfies the reproducing property: $h(x) = \langle \ker(\cdot, x), h \rangle_\Hc$ for all $h \in \Hc$ and $x \in \Xf$; for more details for RKHS please refer to \cite[Chapter 4]{steinwart2008support}. Moreover, assume $\ker(\cdot,\cdot)$ is measurable, i.e.,  $\ker(\cdot,x)$ is a measurable function for each $x\in \Xf$, then each member $h$ of $\Hc$ is a measurable function on $\Xf$ \cite[Lemma 4.24]{steinwart2008support}.

Denote by $\Mc_b(\Xf,\Xc)$ the space of all finite signed measures on $(\Xf,\Xc)$. Each $\Pb\in \Mc_b(\Xf,\Xc)$ defines  a linear map $h\mapsto \int_\Xf h d\Pb $ on $\Hc$. Suppose  $\ker(\cdot,\cdot)$ is bounded hereafter, i.e., $\|\ker\|_\infty:=\sup_{x\in \Xf} \sqrt{\ker(x,x)} <\infty $, and then the above linear map is bounded and hence $\Pb$ can be identified as a member $\mu(\Pb)$ in $\Hc$ by Riesz Representation Theorem \cite[Lemma 26]{sriperumbudur2011universality}, given as below:
\begin{equation}
\label{eq:mu}
\mu(\Pb)(\cdot) = \int \ker(\cdot,x) d\Pb(x)\in \Hc, \quad \forall \Pb\in \Mc_b(\Xf,\Xc),
\end{equation}
and satisfies 
$$
\langle  \mu(\Pb), h\rangle_\Hc = \int_\Xf h d\Pb = \int \langle \ker(\cdot,x), h \rangle_\Hc d\Pb(x), \quad \forall h\in \Hc.
$$

Denote by $\Pc(\Xf,\Xc)$ the space of all probability measures on $(\Xf,\Xc)$. 
Then, the \emph{maximum mean discrepancy} (MMD) associated with the kernel $\ker$ for a pair $\Pb,\Qb\in \Pc(\Xf,\Xc)$ is defined as, cf. \cite{gretton2012kernel}:
\begin{align*}
\MMD(\Pb,\Qb;\ker)  =& \sup_{h\in \Hc, \|h\|_{\Hc}=1} \left|\int_\Xc hd\Pb - \int_\Xc h d\Qb \right|.  
\end{align*}
Moreover, from the reproducing property it can be easily shown that
\begin{align}
\label{eq:mmd}
\MMD^2(\Pb,\Qb;\ker) = & \|\mu(\Pb)-\mu(\Qb)\|^2_\Hc \nonumber \\
= & \Eb \ker(Z,Z') - 2 \Eb \ker(Z,Y) + \Eb \ker(Y,Y'),
\end{align}
where $Z$ and $Z'$ are independent random variables with distribution $\Pb$, and $Y$ and $Y'$ independent random variables with distribution $\Qb$ \cite[Lemma 4, Lemma 6]{gretton2012kernel}. 

Next, a bounded measurable kernel is called a \emph{characteristic kernel} if the map $\mu: \Pc(\Xf,\Xc)\to\Hc$ is injective, i.e., $\MMD(\Pb,\Qb)=0$ if and only if $\Pb=\Qb\in \Pc(\Xf,\Xc)$. If a kernel is bounded, measurable and characteristic, then $\MMD$ is a valid a metric on $\Pc(\Xf,\Xc)$.  Thus, the map $\mu$ provides a natural embedding of the space of probability measures $\Pc(\Xf,\Xc)$ into the reproducing kernel Hilbert spaces $\Hc$ associated with the characteristic kernel $\ker(\cdot,\cdot)$.

Now we verify that the map $\mu: \Mc_b(\Xf,\Bc(\Xf))\to \Hc$ is injective. The first lemma generalizes \cite[Lemma 5]{gretton2012kernel}. 
Note that if $\ker(\cdot,\cdot)$ is bounded, then each member in $\Hc$ is a bounded function on $\Xf$ by \cite[Lemma 4.23]{steinwart2008support}.   Let $\bar{\Hc}$ denote the closure of $\Hc$ w.r.t. the uniform metric. Denote by $C_b(\Xf)$ the space of all bounded continuous functions on $\Xf$.

\begin{lem}
\label{lem:injectiveembedding}
   Suppose that $\Xf$ is metrizable. Consider a measurable bounded kernel $\ker(\cdot,\cdot)$. Suppose  $\bar{\Hc}$ contains $C_b(\Xf)$. Then the map $\mu: \Mc_b(\Xf,\Bc(\Xc))\to \Hc$ is injective. In particular, $\ker(\cdot,\cdot)$ is characteristic.
\end{lem}

The assumptions in Lemma \ref{lem:injectiveembedding} are weaker than the universal assumption in \cite[Lemma 4.23]{steinwart2008support}. Finer characterization can be found under further topological assumptions on the domain $\Xf$. For instance,  \cite{sriperumbudur2011universality} studies when the map $\mu$ is injective on the space of finite signed Radon measures on a locally compact Hausdorff space $\Xf$.  
\textcolor{blue}{In the appendix Lemma \ref{lem:injectiveembeddingRd}} are summarized from \cite[Theorem 6, Proposition 11 and Proposition 16]{sriperumbudur2011universality} which provides useful criteria for verifying whether $\mu$ is injective or not.
In particular, Lemma \ref{lem:injectiveembeddingRd} immediately implies that Gaussian and Laplace kernels have injective $\mu$ on $\Mc_b(\Rb^d,\Bc(\Rb^d))$. \cite{fukumizu2004dimensionality} also contains additional results on the injectivity of the map $\mu$; \myred{see also \cite{muandet2017kernel} for a survey paper on kernel embedding.}  

\paragraph{Minimum MMD estimators}
Consider a bounded measurable kernel $\ker(\cdot,\cdot)$ on the space $(\Xf,\Xc)$. Let $\Fc_1$ to be the unit ball in the associate RKHS $\Hc$, and take 
$$\myF = \myF_1= \biggr \{\theta \mapsto \int f_1(x) d\Pb_\theta(x) \, \biggr | \, f_1\in \Fc_1  \biggr \}.$$ 
Then 
\begin{align}
    \sup_{\myf\in \myF_1} |G \myf-H \myf| = \sup_{f_1 \in \Fc_1} \left|\int f_1 d\Pb_G-\int f_1 d\Pb_H\right|= \MMD(\Pb_G,\Pb_H;\ker),
\end{align}
the \emph{maximum mean discrepancy} between $\Pb_G$ and $\Pb_H$. When the kernel is clear, we write $\MMD(\Pb_G,\Pb_H)$ for $\MMD(\Pb_G,\Pb_H;\ker)$.

A minimum IPM estimator for this choice of $\myF_1$ is now called a minimum MMD estimator associated with the kernel function $\ker$, and takes the form
\begin{align*}
\hat{G}_n\in & \argmin_{G'\in \Gc_k(\Theta)} \MMD(\Pb_{G'},\hat{\Pb}_n) \\
= & \argmin_{G'\in \Gc_k(\Theta)} \MMD^2(\Pb_{G'},\hat{\Pb}_n) \\
= & \argmin_{G'\in \Gc_k(\Theta)} \int K(\theta,\theta') dG'(\theta) dG'(\theta') - 2  \int J_n(\theta)  dG'(\theta)
\end{align*}
where the last line follows from Eq.~\eqref{eq:mmd}, and
$
K(\theta,\theta') : = \int \ker(z,z') d\Pb_\theta(z) d\Pb_{\theta'}(z')
$, and $J_n (\theta) =  \frac{1}{n} \sum_{i\in [n]} \int \ker(x,X_i) d\Pb_\theta(x)$. A summary of the computation of a minimum MMD estimator is available in Algorithm 3, with $\Delta^{k-1}:=\{p'\in \Rb^k \mid p'_i\geq 0, \ \sum_{i\in [k]}p'_i=1  \}$  the probability simplex. 
\begin{algorithm}
\caption{Minimum MMD estimators} \label{alg:three}
\KwData{$X_1,\ldots, X_n\overset{\iidtext}{\sim} \Pb_{G^*}$}
\KwResult{$\hat{G}_n$}
$(\hat{p},\hat{\theta}_1,\ldots,\hat{\theta}_k)\in \argmin_{\substack{p'\in \Delta^{k-1}\\ \theta'_1,\ldots,\theta'_k\in  \Theta}}  \sum_{i,j\in [k]} p'_i p'_j K(\theta'_i,\theta'_j) - 2 \sum_{i\in [k]} p'_i J_n(\theta'_i)$; \\ 
$\hat{G}_n =\sum_{i\in [k]} \hat{p}_i\delta_{\hat{\theta}_i}$
\end{algorithm}

\myred{Minimum MMD distance estimators have been studied in \cite{briol2019statistical, cherief2022finite} where they focus on density estimation rate, not parameter convergence rates as in our paper. They did not specifically apply the estimators to mixture models; although it is worth to mention that \cite{cherief2022finite} did apply their results to dictionary, which can be viewed as a special case of mixture models with $\theta_i$ known. There are various computational algorithm like stochastic/projective gradient descent proposed in \cite{briol2019statistical, cherief2022finite} to compute the optimization problem above and since it is not the focus of our paper, we refer the interested readers to them for more details. }

\paragraph{Theoretical properties} We discuss next the properties of this estimator.

\begin{lem}
\label{lem:distinguishableMMD}
    Suppose that the map $\mu: \Mc_b(\Xf,\Xc)\to \Hc$ is injective, and that the mixture model is identifiable, then $\Gc_k(\Theta)$ is distinguishable by $\myF_1$, i.e., for any $G\neq H\in \Gc_k(\Theta)$, $\MMD(\Pb_G,\Pb_H)>0$. 
\end{lem}

The above lemma is straightforward and thus its proof is omitted. 
We next establish the inverse bounds \eqref{eqn:globalinversebound} and \eqref{eqn:localinversebound} by applying the general Theorem \ref{thm:inversebound}. To do this we need some regularity on the family of components $\{\Pb_{\theta}\}_{\theta\in \Theta}$, so that the sufficient condition for Theorem~\ref{thm:inversebound} on the function class $\Phi = \Phi_1$ can be verified.

\begin{assump}
Suppose that $\{\Pb_{\theta}\}_{\theta\in \Theta}$ has density $\{p(x \mid \theta)\}_{\theta\in \Theta}$ w.r.t. a dominating measure $\lambda$ on $(\Xf,\Xc)$.
The family $\{p(x \mid \theta)\}_{\theta\in \Theta}$ is said to satisfy Assumption $A2(m)$ if the following holds.
Suppose for any $\alpha\in \Ic_m$, $D^\alpha p(x \mid \theta)$ exists and as a function of $\theta$ is continuous on $\Theta$. Moreover, for any $\gamma\in \Ic_{m-1}$, any $i\in [q]$ and any $\theta\in \Theta$, there exists a constant $\Delta_\theta>0$ 
such that for any $0<|\Delta|<\Delta_\theta$:
$$
\left|\frac{D^\gamma p(x \mid \theta+\Delta e_i)- D^\gamma p(x \mid \theta)}{\Delta}\right| \leq \psi_\theta(x), \quad \lambda-a.e.\ x\in \Xf,
$$
with $\int_{\Xf} \psi_\theta(x) d\lambda<\infty$.  Moreover, for any $\gamma\in \Ic_m\setminus\Ic_{m-1}$, and any $\theta\in \Theta$, there exists a constant $\Delta'_\theta>0$ 
such that for any $0<\|\Delta'\|_2<\Delta'_\theta$:
$$
\left|D^\gamma p(x \mid \theta+\Delta' )- D^\gamma p(x \mid \theta)\right| \leq \psi'_\theta(x), \quad \lambda-a.e.\ x\in \Xf,
$$ 
with $\int_{\Xf} \psi'_\theta(x) d\lambda<\infty$.
\end{assump}
     
\begin{lem}
\label{lem:exchangealorder}
    If $\{p(x \mid \theta)\}_{\theta\in \Theta}$ satisfies Assumption $A2(m)$, then for any essentially bounded measurable function on $\Xf$, i.e., any $f_2\in L^\infty(\Xf,\Xc,\lambda)$,
    the function $\theta \mapsto \Psi(\theta) = \int f_2(x) p(x \mid \theta) d\lambda$ is $m$-th order continuously differentiable, and 
    $$
D^\alpha \Psi(\theta) = D^\alpha \int_{\Xf} f_2(x) p(x \mid \theta) d\lambda =  \int_{\Xf} f_2(x) D^\alpha p(x \mid \theta) d\lambda.
$$
\end{lem}

\begin{definition} \label{def:linearindependent}
The family $\{p(x \mid \theta)\}_{\theta\in \Theta}$ of functions on $\Xf$  is said to be a \textit{$(m,k_0,k)$ linear independent} if the following hold. 1) For $\lambda$-a.e. $x\in \Xf$, $p(x \mid \theta)$ is $m$-th order continuously differentiable on $\Theta$; and 2) Consider any integer $\ell\in [k_0,2k-k_0]$, and any vector $(m_1,m_2,\ldots,m_\ell)$ such that $1\leq m_i\leq m+1$ for $i\in [\ell]$ and $\sum_{i=1}^\ell m_i\in [2k_0, 2k]$.
For any distinct $\{\theta_i\}_{i\in [\ell] }\subset \Theta$, 
\begin{subequations}
\begin{align}
\sum_{i=1}^{\ell}\  \sum_{|\alpha|\leq m_i-1} a_{i\alpha} D^\alpha p(x \mid \theta_i )  = &0,  \quad \lambda-a.e.\ x\in \Xf \label{eqn:linearinda}\\
\sum_{i\in [\ell]}   a_{i\bm{0}}   = &0, \label{eqn:linearindb}
\end{align}
\end{subequations}
if and only if 
$$
a_{i\alpha}=0, \quad \forall \ 0\leq |\alpha|< m_i, \ i\in [\ell].
$$
\end{definition}

Specializing the $m$-strong identifiability from Definition \ref{def:mstrong} to $\{p(x\mid \theta)\}$ gives the following.

\begin{definition} 
\label{def:strideden}
The family $\{p(x \mid \theta)\}_{\theta\in \Theta, x\in \Xf}$ is said to be a \textit{$m$-strongly identifiable} if the following hold. 1) For $\lambda$-a.e. $x\in \Xf$, $p(x \mid \theta)$ is $m$-th order continuously differentiable on $\Theta$. 2) 
For any distinct $\{\theta_i\}_{i\in [\ell] }\subset \Theta$, 
\begin{subequations}
\begin{align*}
\sum_{i=1}^{\ell}\  \sum_{\alpha\in \Ic_m} a_{i\alpha} D^\alpha p(x \mid \theta_i )  = &0,  \quad \lambda-a.e.\ x\in \Xf 
\end{align*}
\end{subequations}
if and only if 
$$
a_{i\alpha}=0, \quad \forall \ 0\leq |\alpha|< m_i, \ i\in [\ell].
$$
\end{definition}
Clearly, $\{p(x \mid \theta)\}$ is $m$-strongly identifiable implies that it is $(m,k_0,k)$ linear independent. In previous work \cite{chen1995optimal, nguyen2013convergence,ho2016strong, holzmann2004identifiability} established the connection between mixture models and the  
 $m$-strong identifiability of $\{p(x \mid \theta)\}_{\theta\in\Theta}$ for the case $m=1,2$, and they also showed many density families are  $m$-strongly identifiable.   

\begin{lem}
\label{lem:invbouMMD}
Suppose that $\Theta$ is compact and the map $\mu: \Mc_b(\Xf,\Xc)\to \Hc$ given in Eq.~\eqref{eq:mu} is injective. 
\begin{enumerate}[label=(\alph*)]
     \item \label{lem:invbouMMDa}
     If $\{p(x \mid \theta)\}_{\theta\in \Theta}$ satisfies Assumption $A2(2d_1-1)$ and  $\{p(x\mid \theta)\}_{\theta\in \Theta}$ is $(2d_1-1,k_0,k)$ linear independent, then $\myF_1$ is a $(2d_1-1,k_0,k)$ linear independent domain, and hence the inverse bounds \eqref{eqn:localinversebound} and \eqref{eqn:localinverseboundmoment} hold for any $G_0\in \Ec_{k_0}(\Theta)$: \myred{for $D\in \{W_{2d_1-1}^{2d_1-1}, \mbf_{2d_1-1}\}$}
\begin{equation}
\liminf_{ \substack{G,H\overset{W_1}{\to} G_0\\ G\neq H\in \Gc_k(\Theta) }} \frac{\MMD(\Pb_G,\Pb_H)}{D(G,H)} >0. \label{eqn:localinverseboundMMD}
\end{equation}
     \item \label{lem:invbouMMDb}
      If $\{p(x \mid \theta)\}_{\theta\in \Theta}$ satisfies Assumption $A2(2k-1)$ and $\{p(x \mid \theta)\}_{\theta\in \Theta}$ is $(2k-1,1,k)$ linear independent, then $\myF_1$ is a $(2k-1,1,k)$ linear independent domain. Moreover if the mixture model is identifiable on $\Gc_k(\Theta)$, then the inverse bound \eqref{eqn:globalinversebound} holds: \myred{for $D\in \{W_{2k-1}^{2k-1}, \mbf_{2k-1}\}$}
    \begin{equation}
	\inf_{  G\neq H\in \Gc_k(\Theta) } \frac{\MMD(\Pb_G,\Pb_H)}{D(G,H)} >0. \label{eqn:globalinverseboundMMD}
\end{equation}
\end{enumerate}
\end{lem}
\begin{rem} \label{strongidentifiability}
Under similar assumptions of Lemma \ref{lem:invbouMMD}, by mimicking the proof of Lemma \ref{lem:invbouMMD}, we can also obtain that $\{p(x \mid \theta)\}_{\theta\in\Theta}$ is $m$-strongly identifiable implies that $\myF_1$ is $m$-strongly identifiable.  \myeoe
\end{rem}

It remains to control the governing empirical process: 
$$
\sup_{f_1 \in \Fc_1}\left|\int f_1 d\Pb_{G^*}-\frac{1}{n}\sum_{i\in [n]}f_1(X_i)\right| = \MMD(\Pb_{G^*}, \hat{\Pb}_n).
$$
The next lemma is a result that controls the empirical process for any probability measure $\Pb$, not necessarily the mixture probability measures $\Pb_G$. 
\begin{lem}
\label{lem:empiricalprocessMMD}
Consider a measurable bounded kernel $\ker(\cdot,\cdot)$.  Then
    $$ \sup_{\Pb\in \Pc(\Xf,\Xc)} \Eb_{\Pb} \MMD(\Pb, \hat{\Pb}_n) 
    \leq \frac{2\|\ker\|_\infty}{\sqrt{n}}, $$
    where $\Eb_\Pb$ denotes the expectation when the random variables $\{X_i\}_{i\in [n]} \overset{\iidtext}{\sim} \Pb$.
\end{lem}

\myred{A probabilistic version of Lemma \ref{lem:empiricalprocessMMD} is also available; see Lemma \ref{lem:empiricalprocessMMDconcentration}; see also \cite[Lemma 1]{briol2019statistical} or \cite[Theorem 3.2]{cherief2022finite}.} Combining Theorem \ref{thm:convergencerate}, Lemma \ref{lem:invbouMMD} and Lemma \ref{lem:empiricalprocessMMD} immediately gives the following theorem.

\begin{thm}
\label{thm:MMDuniformrate}
Suppose that $\Theta$ is compact and the mixture model is identifiable on $\Gc_k(\Theta)$. Let $\hat{G}_n$ be a minimum MMD estimator. Consider a bounded and measurable kernel $\ker(\cdot,\cdot)$ on $(\Xf,\Xc)$ and the map $\mu: \Mc_b(\Xf,\Xc)\to \Hc$ is injective. 
 \begin{enumerate}[label=(\alph*)]
\item 
If $\{p(x \mid \theta)\}_{\theta\in \Theta}$ satisfies Assumption $A2(2k-1)$ and $\{p(x \mid \theta)\}_{\theta\in \Theta}$ is $(2k-1,1,k)$ linear independent. Then there exists a constant $C$, where its dependence on $ \Theta, k$ and the probability kernel $\{\Pb_\theta\}$ is suppressed, such that \myred{for $D\in \{W_{2k-1}^{2k-1}, \mbf_{2k-1}\}$}
\begin{align*}
\sup_{G^*\in \Gc_k(\Theta)}	\Eb_{G^*} D(G^*,\hat{G}_n)
	\leq C \frac{\|\ker\|_\infty}{\sqrt{n}} . \label{eqn:globalconvergencerategeneralMMD}
\end{align*}

\item 
If $\{p(x\mid \theta)\}_{\theta\in \Theta}$ satisfies Assumption $A2(2d_1-1)$ and  $\{p(x \mid \theta)\}_{\theta\in \Theta}$ is $(2d_1-1,k_0,k)$ linear independent, 
then for any $G_0\in \Ec_{k_0}(\Theta)$, there exists $r(G_0)$, $C(G_0)$ and $c(G_0)$, where their dependence on $ \Theta, k_0,k$ and the probability kernel $\{\P_\theta\}$ are suppressed, such that \myred{for $D\in \{W_{2d_1-1}^{2d_1-1}, \mbf_{2d_1-1}\}$}

\begin{align*}
\sup_{G^*\in \Gc_k(\Theta): W_1(G_0,G^*)<r(G_0)} \Eb_{G^*} D(\hat{G}_n,G^*) 
\leq  C(G_0)  \frac{\|\ker\|_\infty}{\sqrt{n}}. 
\end{align*}

\end{enumerate}
\end{thm}

By a similar argument as given in Remark \ref{rem:minimaxrateKS}, we conclude that a minimum MMD estimator also achieves the minimax optimal rate, under a rather general setting of the data domain $\Xf$. 


\vspace{2mm}
\noindent
\textbf{Example: Multi-dimensional Gaussian mixture models} \\
\myred{Next, we apply the general theory to study multi-dimensional Gaussian mixture models. More specifically, the density for each component is $\Pb_\theta = \Nc(\theta,\Sigma) $ on $\Rb^d$ where $\Sigma$ is a known covariance matrix. The Gaussian mixture model is 
$$
\Pb_{G^*} =  \sum_{i}^k p^*_i \Nc(\theta_i^*,\Sigma), 
$$
and the goal is estimate $G^*$ based on i.i.d. samples $X_1,\ldots,X_n\sim \Pb_{G^*}$. Note in this case the dimension $q$ for parameter $\theta$ is the same as the dimension of the samples $d$. One can verify easily that $\{p(x\mid\theta)\}$ satisfies Assumption A2(m) for any $m\geq 0$. It follows from the classical result \cite[Proposition 1]{teicher1963identifiability} that this mixture model is identifiable. In fact, this family 
is $m$-strongly identifiable for any $m\geq 0$ due to Lemma \ref{eqn:locationmixture} when $d=1$. For general $d$, when restricted to Gaussian location mixtures, it is also straightforward.

Now consider the kernel $\ker(x,y)=\exp(-\gamma \|x-y|_2^2)$. It is clearly a bounded and measurable kernel on $(\Rb^d,\Bc(\Rb^d))$. Moreover the map $\mu: \Mc_b(\Rb^d,\Bc(\Rb^d))\to \Hc$ is injective, as discussed \myred{after Lemma \ref{lem:injectiveembeddingRd}.} Thus all assumptions in Theorem \ref{thm:MMDuniformrate} are verified for this example so then the uniform convergence rates are obtained.

Finally, on the computational aspect, notice that for this choice of kernel, we have analytical expressions for $K(\theta,\theta')$ and $J_n(\theta)$ from Algorithm \ref{alg:three}:
\begin{align*}
K(\theta,\theta') = & \frac{1}{\sqrt{\det(1+4\gamma\Sigma)}} e^{- \gamma (\theta- \theta')^\top (1+4\gamma\Sigma)^{-1}  (\theta- \theta')}, \\ 
J_n(\theta)= & \frac{1}{n}\sum_{i=1}^n \frac{1}{\sqrt{\det(1+2\gamma\Sigma)}} e^{- \gamma (X_i-\theta)^\top (1+2\gamma\Sigma)^{-1}  (X_i-\theta)}.
\end{align*}
So the gradient can also be computed in analytical form. The specific minimization in Algorithm \ref{alg:three} can be solved by various numerical optimization methods, say stochastic gradient descent, projective gradient descent, or coordinate descent (by viewing $\theta'_1,\ldots, \theta'_k$ as one coordinate and viewing $p'$ as the other coordinate). We leave the computational details to interested readers.}

\subsection{Moment based estimators}
\label{sec:momentdeviation}

In this section we consider the monomial family $\myF_2:=\{(\theta-\theta_0)^\alpha\}_{\alpha\in \Ic_{2k-1}}$, where $\theta_0$ is an arbitrarily chosen element in $\Rb^q$. The univariate case $q=1$ has been presented in Example \ref{exa:moment}. Unlike the $\myF$ for minimum IPM estimators in Section \ref{sec:integralprobabilitymetric},
we will see that $\myF_2$ already satisfies the inverse bounds, \myred{so it remains  to guarantee that $\myF_2$ is estimatable.} 

\paragraph{Inverse bounds} We first show that $\myF_2$ is a $(2k-1,1,k)$ linear independent domain. It is obvious that each monomial in $\myF_2$ is $2k-1$ differentiable. Consider any integer $\ell\in [1,2k-1]$, and any vector $(m_1,m_2,\ldots,m_\ell)$ such that $1\leq m_i\leq 
2k$ for $i\in [\ell]$ and $\sum_{i=1}^\ell m_i\in [2, 2k]$. Consider any distinct $\{\theta_i\}_{i\in [\ell] }\subset \Theta$. The equations \eqref{eqn:linearinddomaina} \eqref{eqn:linearinddomainb} become
\begin{align}
\sum_{i=1}^{\ell}\  \sum_{|\gamma|\leq m_i-1} a_{i\gamma} \frac{\alpha!}{(\alpha-\gamma)!}  (\theta_i-\theta_0)^{\alpha-\gamma}1_{\alpha\geq \gamma}  = &0,  \quad \alpha \in \Ic_{2k-1}. 
\end{align}
It then follows from Lemma \ref{lem:linearalgebra} 
that $a_{i\gamma}=0$ for any $\gamma\in \Ic_{m_i-1}, i\in [\ell]$. So $\myF_2$ is a $(2k-1,1,k)$ linear independent domain. (Note that it is straightforward to see that $\myF_2$ is not $(2k-1)$-strongly identifiable.) 
Provided that $\Theta$ is compact, then
we may apply Theorem \ref{thm:inversebound}, which yields that \eqref{eqn:localinversebound} holds for $\myF=\myF_2$, for any $G_0\in \Ec_{k_0}(\Theta)$ for any $k_0\in [k]$, and is as below: 
\begin{equation}
\liminf_{ \substack{G,H\overset{W_1}{\to} G_0\\ G\neq H\in \Gc_k(\Theta) }} \frac{\|\mbf_{2k-1}(G-\theta_0)-\mbf_{2k-1}(H-\theta_0)\|_{\infty}}{W_{2d_1-1}^{2d_1-1}(G,H)} >0. \label{eqn:localmomentbound}
\end{equation}
%
%
Since discrete distributions with $k$ support points are uniquely characterized by their first $2k-1$ moments, but not their first $2k-2$ moments by Lemma \ref{lem:discretedistributions}, 
we know that $\Gc_k(\Theta)$ is distinguishable by $\myF_2$.
By Lemma \ref{lem:localtoglobal}, \eqref{eqn:globalinversebound} holds for $\myF=\myF_2$ 
and is as below: 
\begin{equation}
	\inf_{  G\neq H\in \Gc_k(\Theta) } \frac{\|\mbf_{2k-1}(G-\theta_0)-\mbf_{2k-1}(H-\theta_0)\|_{\infty}}{W_{2k-1}^{2k-1}(G,H)} >0. \label{eqn:globalmomentbound}
\end{equation}
The next lemma summarizes the discussions up to this point in this subsection. 
\begin{lem} \label{lem:momentbound}
The family $\myF_2$ is a $(2k-1,1,k)$ linear independent domain, and $\Gc_k(\Theta)$ is distinguishable by $\myF_2$. If additionally $\Theta\subset \Rb^q$ is  compact, then the local inverse bound \eqref{eqn:localmomentbound} holds for any $G_0\in \Ec_{k_0}(\Theta)$ for any $k_0\in [k]$. Moreover,  \eqref{eqn:globalmomentbound} holds. 
\end{lem}
\begin{rem}
\label{rem:mominvcom}
 The univariate case $q=1$ of \eqref{eqn:globalmomentbound} was first established by \cite[Proposition 1]{wu2020optimal}. 
 \myred{ The \eqref{eqn:globalmomentbound} for general $q$ was implied by the Theorem 4.2 and Equation (4.49) in \cite{doss2020optimal}.  It is worth mentioning that both previous bounds specify the dependence on the parameters $k$ and $q$. 
 Lemma \ref{lem:momentbound} produces similar results by specializing Theorem \ref{thm:inversebound} 
 for $\myF=\myF_2$. 
 }
 Moreover, we also obtain the local version  \eqref{eqn:localmomentbound}, which is new to the best of our knowledge. $2k-1$ in the numerator is the smallest number for the local moment inverse bound \eqref{eqn:localmomentbound}; for details, see Lemma \ref{lem:optimallocalmomentinversebou} in the Appendix. \myred{One specific instance to demonstrate the usefulness of this result is to study mixture of multinomials is  in Example \ref{exa:mixmuldis} ahead, which cannot be deduced directly from existing results.} \myeoe  
\end{rem}

\paragraph{Estimation of $\Phi=\Phi_2$}
So far the discussions only concerns properties about discrete distributions in $\Gc_k(\Theta)$ and does not involve mixture models or the probability kernel $\{\Pb_{\theta}\}_{\theta\in\Theta}$. To ensure that $\myF_2$ is estimatable, it is required that for each $\myf = \left(\theta-\theta_0\right)^\alpha\in \myF_2$, where $\alpha\in \Ic_{2k-1}$, there exists a function $t_\alpha$ defined on $\Xf$ such that 
\begin{equation}
G \myf= m_\alpha(G-\theta_0) =\Eb_{G}t_\alpha(X), \quad \forall G\in \Gc_k(\Theta). \label{eqn:momenttj}
\end{equation}  
A minimum $\myF$-distance estimator in this case becomes 
\begin{equation}
\hat{G}_n\in \argmin_{G'\in \Gc_k(\Theta)} \sup_{\alpha\in \Ic_{2k-1}}\left|m_\alpha(G'-\theta_0) -\frac{1}{n}\sum_{i\in [n]}t_\alpha(X_i)\right|.  \label{eqn:momentmethod}
\end{equation}
This shall be called a \emph{generalized method of moments} (GMM).
A summary of the estimation procedure is Algorithm \ref{alg:moments}. In a standard moment-based estimation method, the statistic $t_\alpha$ may be  taken to be a power function or power product function (i.e., a monomial) of the variable $x \in \Xf$. 
For many standard families of probability kernels $\Pb_{\theta}$, this choice of statistic results in the expectation $\Eb_{G}t_\alpha(X)$ taking the form of monomials of the parameter vector $\theta$. In general, we may use any other choices of statistic function $t_\alpha$ as well, as long as they can be used to define the functions $\mbf_\alpha(G-\theta_0)$ in the sense of Eq.~\eqref{eqn:momenttj}. It is in this sense that we use the term "generalized". 

\begin{rem}
In the description of Algorithm \ref{alg:moments} note that the $\bar{\tbf}$ is an empirical estimate of $\mbf_{2k-1}(G-\theta_0)$, and might not lie in a valid moment space for a discrete distribution due to the randomness, but the parameter estimate may be obtained by finding the closest corresponding moment vector w.r.t. $\|\cdot\|_\infty$. Specializing Algorithm \ref{alg:moments} when the probability kernel $\Pb_\theta$ is a univariate Gaussian distribution, we obtain the "denoised method of moments" algorithm that was investigated by \cite{wu2020optimal}. 
\end{rem}

\begin{algorithm}
\caption{Generalized method of moments} \label{alg:moments}
\KwData{$X_1,\ldots, X_n\overset{\iidtext}{\sim} \Pb_{G^*}$}
\textbf{Parameter}: $\theta_0$ \\
\KwResult{$\hat{G}_n$}
$\bar{t}_\alpha(\theta_0) \gets \frac{1}{n}\sum_{i\in [n]}t_\alpha(X_i)$, for $\alpha\in \Ic_{2k-1}$\;
$\hat{G}_n\in \argmin_{G'\in \Gc_k(\Theta)} \| \mbf_{2k-1}(G'-\theta_0) - \bar{\tbf}  \|_\infty$, where $\bar{\tbf}=(\bar{t}_\alpha(\theta_0))_{\alpha\in \Ic_{2k-1}}$.
\end{algorithm}


\myred{To compute the minimizers of the above algorithm, one may consider projection based algorithm (especially for one dimensional mixtures); see \cite{wu2020optimal} for details. This is due to that the cardinality of $\myF_2$ is finite. }

We now state Theorem \ref{thm:convergencerate} specialized to the GMM estimators. 

\begin{thm} \label{thm:convergencerategeneralmoment}
	Suppose that $\Theta$ is compact. Suppose that for each $\alpha \in \Ic_{2k-1}$, there exists a real-valued function $t_\alpha$ defined on $\Xf$ such that \eqref{eqn:momenttj} holds. Let $\hat{G}_n$ be the output of Algorithm \ref{alg:moments}.
 \begin{enumerate}[label=(\alph*)] 
\item \label{itema:thm:convergencerategeneralmoment}
Then there exists  $C$, where its dependence on $ \Theta, k$ and the probability kernel $\{\Pb_\theta\}$ is suppressed, such that for any $G^*\in \Gc_k(\Theta)$, \myred{ and for $D\in \{W_{2k-1}^{2k-1}, \mbf_{2k-1}\}$}
\begin{align}
	\Pb_{G^*} \left(D(G^*,\hat{G}_n)\geq t \right) 
	\leq  \Pb_{G^*} \left( \sup_{\alpha\in \Ic_{2k-1}}\left| \frac{1}{n}\sum_{i\in [n]}t_\alpha(X_i)-m_\alpha(G^*-\theta_0)\right| \geq C  t  \right). \label{eqn:convergenceratemomentgeneral}
\end{align}
and
\begin{align*}
	\Eb_{G^*} D(G^*,\hat{G}_n)
	\leq C \Eb_{G^*}  \sup_{\alpha\in \Ic_{2k-1}}\left| \frac{1}{n}\sum_{i\in [n]}t_\alpha(X_i)-m_\alpha({G^*}-\theta_0)\right| . 
\end{align*}
 
\item \label{itemb:thm:convergencerategeneralmoment}
Fix $G_0\in \Ec_k(\Theta)$.  Then there exists $r(G_0)$, $C(G_0)$ and $c(G_0)$, where their dependence on $ \Theta, k_0,k$ and the probability kernel $\{\Pb_\theta\}$ is suppressed, such that for any ${G^*} \in \Gc_k(\Theta)$ satisfying $W_1(G_0,{G^*})<r(G_0)$,
\myred{ and for $D\in \{W_{2d_1-1}^{2d_1-1}, \mbf_{2d_1-1}\}$}
\begin{align} 
\Pb_{G^*} \left(D({G^*},\hat{G}_n)\geq t \right) 
\leq   \Pb_{G^*} \left( \sup_{\alpha\in \Ic_{2k-1}}\left| \frac{1}{n}\sum_{i\in [n]}t_\alpha(X_i)-m_\alpha({G^*}-\theta_0)\right| \geq C(G_0)  t  \right)\label{eqn:convergenceratemomentlocalgeneral}
\end{align}
and
\begin{align*}
\Eb_{G^*} D(\hat{G}_n,{G^*}) 
\leq  C(G_0)  \Eb_{G^*} \sup_{\alpha\in \Ic_{2k-1}}\left| \frac{1}{n}\sum_{i\in [n]}t_\alpha(X_i)-m_\alpha({G^*}-\theta_0)\right|. 
\end{align*}
\end{enumerate}
	\end{thm}

\paragraph{Examples: Location mixtures of exponential families with quadratic variance functions}

As an illustration of the applicability of the GMM and Theorem~\ref{thm:convergencerategeneralmoment}, we present a class of probability kernel $\{\Pb_{\theta}\}$ of which $t_\alpha$ with the property \eqref{eqn:momenttj} exists. By applying Theorem \ref{thm:convergencerategeneralmoment}, the right hand sides of  \eqref{eqn:convergenceratemomentgeneral} and \eqref{eqn:convergenceratemomentlocalgeneral} are also calculated to obtain convergence rates.
In particular, we consider \emph{the natural exponential families with quadratic variance functions} (NEF-QVF), where within each family the variance of the random variable is a quadratic function of the mean-value parameter. NEF-QVF is shown in  \cite{morris1982natural} to contain only six probability families and their linear transformations. They are
\begin{align*}
    \text{Gaussian:} &\quad  f(x \mid \xi,\sigma)  =  \frac{1}{\sqrt{2\pi}\sigma} e^{-\frac{(x-\xi)^2}{2\sigma^2}} \quad  \forall x\in \Rb, \xi\in \Rb, \sigma>0; \\
        \text{Poisson:} & \quad  f(x \mid \lambda)  = e^{-\lambda} \frac{\lambda^x}{x!}  \quad  \forall x\in \Nb, \lambda>0; \\
        \text{gamma:}& \quad  f(x \mid \alpha,\beta)  =  \frac{\beta^\alpha}{\Gamma(\alpha)} x^{\alpha-1}e^{-\beta x} \quad  \forall x>0, \alpha,\beta >0;\numberthis \label{eqn:NEFQVF}\\
        \text{binomial:} & \quad f(x \mid m,p) =  \binom{m}{x} p^x (1-p)^{m-x}, \quad \forall x\in \Nb, 0<p<1, n\in \Nb_+;\\
        \text{negative binomial:} & \quad f(x \mid r,p) = \frac{\Gamma(x+r)}{x! \Gamma(r)} (1-p)^x p^r, \quad \forall x\in \Nb, 0<p<1, r>0;\\
        \text{NEF-GHS:} & \quad f(x \mid r,\varphi) = e^{\varphi x +r \ln \cos(\varphi)}\  \frac{2^{r-2}}{\Gamma(r)} \prod_{j=0}^\infty \left( 1+ x^2/(r^2+2j)^2 \right)^{-1}, \\
        & \quad \forall x\in \Rb, r>0, \varphi\in \left(-\frac{\pi}{2},\frac{\pi}{2}\right).
\end{align*}
Each of the six univariate families has at most $2$ parameters, and therefore their linear transformations can have at most $4$ parameters. Further details on NEF-QVF can be found in \cite[Section 2]{morris1983natural} or \cite{morris1982natural}. Here in this paper we focus on the above $6$ families; the results on their linear transformations should readily be available following the same procedures.

Following the framework in \cite{lindsay1989moment}, let $p(x\mid \theta)$ be a generic family for the $6$ families in \eqref{eqn:NEFQVF} with the parameter $\theta$ being the mean of the distribution. 
Denote by $\tilde{\Theta}\subset \Rb$ the set of all possible values of $\theta$, which depends on the specific families of probability kernels.  
In particular, if $p(x \mid \theta)$ is the Gaussian family, then the parameter is $\theta=\xi$ and $\sigma$ is known; if $p(x \mid \theta)$ is the negative binomial family, then the parameter $\theta=r(1-p)/p$ is a reparametrization of the parameter $r$ or $p$ with the other known; if $p(x \mid \theta)$ is binomial family, then the parameter $\theta=m p$ is a reparametrization of $p$ while $m$ is fixed---$m$ is not considered as a parameter in this paper since it is discrete-valued. 

It follows from \cite{morris1982natural,lindsay1989moment} that there exists a function $b(\theta)$ and a constant $b_j$, where they both depend on the family of probability kernels $p(x \mid \theta)$ and $b_j$ additionally depends on $j$, such that 
$$t_j(x \mid \theta):= b_j \left(b(\theta)\right)^j \frac{\partial^j p(x \mid \theta)}{\partial \theta^j} \frac{1}{p(x\mid \theta)}$$ 
satisfies 
\begin{equation}
    \Eb_\theta t_j(Y \mid \theta_0) = (\theta-\theta_0)^j   \label{eqn:gammamean}
\end{equation} 
where $Y\sim p(x \mid \theta)$ and $\theta_0\in \tilde{\Theta}^\circ$. 
We may write $t_j(x \mid \theta)=\sum_{i=0}^j a_{j i}(\theta)x^i$ as a polynomial of $x$ where $a_{j i}(\theta)$ is a polynomial of $\theta$ depending on $j$ and the specific family of probability kernels. It follows from \eqref{eqn:gammamean} that for $X\sim \int p(x\mid \theta) dG$ with $G=\sum_{i\in [k]}p_i\delta_{\theta_i}$,
\begin{equation}
\Eb_G t_j(X|\theta_0) = \sum_{i\in [k]} p_i (\theta_i-\theta_0)^j=m_j(G-\theta_0), \label{eqn:polynomialmeanmixture}
\end{equation}
which is the target property \eqref{eqn:momenttj} in the univariate case.

Given $\iidtext$ data $\{X_i\}_{i\in [n]}$, the sample version of the left hand side of \eqref{eqn:polynomialmeanmixture} is  
$$\bar{t}_j(\theta_0):= \frac{1}{n} \sum_{i\in [n]} t_j(X_i \mid \theta_0).$$ The suitable summary statistic for data $\{X_i\}_{i\in [n]}$ is the vector $\bar{\tbf} = (\bar{t}_1(\theta_0),\ldots,\bar{t}_{2k-1}(\theta_0))$. 
\begin{rem} 
\label{rem:theta0}
It can be shown that $\bar{\tbf}$ contains the same information as the first $2k-1$ sample moments $\{\frac{1}{n}\sum_{i\in[n]} X_i^j \}_{j\in [2k-1]}$ since $\{t_j(x \mid \theta_0)\}_{j\in [2k-1]}$ form a family of orthogonal polynomials w.r.t. $p(x \mid \theta_0)$ \cite[Theorem 4]{morris1982natural}.
The choice of $\theta_0\in \tilde{\Theta}^\circ$ is arbitrary and has no theoretical impact on the solution. One convenient choice is $\theta_0=0$, provided that $0\in \tilde{\Theta}^\circ$.
\myeoe
\end{rem}

The parameter space is $\Theta=[M_1,M_2] \subset\tilde{\Theta}$, i.e., the mean parameters $\theta_i$ are assumed to lie in a known compact interval $[M_1,M_2]$. The next lemma analyzes the deviation of $\bar{t}_j(\theta_0)$ from its mean. 
\begin{lem} \label{lem:concentrationboualljNEFQVF} 
	Consider any of the $6$ NEF-QVF families \eqref{eqn:NEFQVF} and let $t_j(\cdot|\theta_0)$ and $\bar{t}_j(\theta_0)$ be defined as above for each specific family of probability kernels $p(x \mid \theta)$. Then there exist  $C$ and $c$, where their dependences on $ \Theta, k,\theta_0$ and the specific NEF-QVF family $\{p(x \mid \theta)\}$ are suppressed, such that for any $\epsilon>0$, 
	\begin{align*}
	\sup_{G\in \Pc(\Theta)}	\Pb_G\left(\max_{j\in [2k-1]} \mid \bar{t}_j(\theta_0)-\Eb_G t_j(X \mid \theta_0)|\geq \epsilon\right) 
		\leq  e^2 (2k-1) \exp\left( -C \min \left\{ n\epsilon^2,  (n\epsilon)^{\frac{1}{2k-1}} \right\}  \right),
	\end{align*}
	and consequently
	$$
	\sup_{G\in \Pc(\Theta)}	\Eb_G \max_{j\in [2k-1]}|\bar{t}_j(\theta_0)-\Eb_G t_j(X \mid \theta_0)|\leq cn^{-\frac{1}{2}}.
	$$
\end{lem}

By combining Lemma \ref{lem:concentrationboualljNEFQVF} and Theorem \ref{thm:convergencerategeneralmoment}, we immediately obtain the following proposition. 
 
\begin{prop}  \label{prop:NEFQVF}
	Consider any of the $6$ NEF-QVF families \eqref{eqn:NEFQVF} and let $t_j(\cdot \mid \theta_0)$ and $\bar{t}_j(\theta_0)$ be defined as above for each specific family of probability kernels $p(x \mid \theta)$. Suppose that $\Theta$ is a compact interval.
		\begin{enumerate}[label=(\alph*)]
 \item \label{itema:prop:NEFQVF}
 Then there exist positive constants $C$ and $c$, where their dependence on $ \Theta, k,\theta_0$ and the probability kernel $\{p(x \mid \theta)\}$ are suppressed, such that \myred{for $D\in \{W_{2k-1}^{2k-1}, \mbf_{2k-1}\}$}
 \begin{align*}
  \sup_{G^*\in \Gc_k(\Theta)}	\Pb_{G^*} \left(D(G^*,\hat{G}_n)\geq t \right) 
 	\leq   e^2 (2k-1) \exp\left( -C \min \left\{ nt^{2},  \left(n t\right)^{\frac{1}{2k-1}} \right\}  \right), 
 \end{align*}
 and consequently,
 \begin{align*}
 \sup_{G^*\in \Gc_k(\Theta)}	\Eb_{G^*} D(G^*,\hat{G}_n )
 	\leq   c n^{-\frac{1}{2}}. 
 \end{align*}

		\item \label{itemb:prop:NEFQVF}
		Fix any $G_0\in \Gc_k(\Theta)$. Then there exists $r(G_0)$,  $C(G_0)$ and $c(G_0)$,  where their dependence on $ \Theta, k,k_0,\theta_0$ and the probability kernel $\{p(x \mid \theta)\}$ are suppressed, such that 
        \myred{for $D\in \{W_{2d_1-1}^{2d_1-1}, \mbf_{2d_1-1}\}$}
	\begin{align*}
		&\sup_{G^*\in \Gc_k(\Theta): W_1(G_0,G^*)<r(G_0) }	\Pb_{G^*} \left(D(G^*,\hat{G}_n)\geq t \right) \\
		\leq  &  e^2 (2k-1) \exp\left( -C(G_0) \min \left\{ nt^{2},  \left(nt\right)^{\frac{1}{2k-1}} \right\}  \right), 
	\end{align*}
     and consequently,
     \begin{align*}
     	\sup_{G^*\in \Gc_k(\Theta): W_1(G_0,G^*)<r(G_0) }	\Eb_{G^*} D(G^*,\hat{G}_n)
     	\leq   c(G_0) n^{-\frac{1}{2}}. 
     \end{align*}
 \end{enumerate}
\end{prop}
By a similar argument as Remark \ref{rem:minimaxrateKS}, we conclude that GMM estimators achieve the minimax optimal rate.

\begin{rem} 
\label{rem:1dgmm}
In Algorithm \ref{alg:moments}, the $\|\cdot\|_\infty$ can be replaced with any other norm while the same conclusion as Proposition \ref{prop:NEFQVF} holds since all norms on $\Rb^{2k-1}$ are equivalent up to a factor of constant. The paper \cite[Theorem 1]{wu2020optimal} obtained the same conclusion as Proposition \ref{prop:NEFQVF} \ref{itema:prop:NEFQVF} for the special case of the univariate location Gaussian mixture. Proposition \ref{prop:NEFQVF} \ref{itema:prop:NEFQVF} extends the previous result to  location mixtures of any NEF-QVF families. \myred{Moreover, the local uniform convergence result   Proposition \ref{prop:NEFQVF} \ref{itemb:prop:NEFQVF}  states that the uniform convergence rate decreases to $n^{-\frac{1}{4d_1 -2}}$ 
once the true mixing measure is constrained to be in a neighborhood of a known $G_0$; a similar local convergence rate when $G^*$ is constrained to have exactly $k_0$ atoms, with mixing weights that are bounded below and atoms that are also well separated, was developed in \cite[Theorem 2]{wu2020optimal} for univariate location Gaussian mixtures.} 
\myeoe
\end{rem}

\myred{
\paragraph{Examples: Multi-dimensional Gaussian mixture models} 
Next, we apply the general theory to study multi-dimensional Gaussian mixture models. More specifically, the density for each component is $\Pb_\theta = \Nc(\theta,\Sigma) $ on $\Rb^d$ where $\Sigma$ is a known covariance matrix. The Gaussian mixture model is 
$$
\Pb_{G^*} =  \sum_{i}^k p^*_i \Nc(\theta_i^*,\Sigma), 
$$
and the goal is estimate $G^*$ based on an $n$-i.i.d. sample $X_1,\ldots,X_n\sim \Pb_{G^*}$. Note in this case the dimension $q$ for parameter $\theta$ is the same as the dimension of the samples $d$.  For the sake of clean presentation in high dimensions, we consider $\theta_0=\bm{0}$, i.e. $\myF_2=\{\theta^\alpha\}_{\alpha\in \Ic_{2k-1}}$, but it is easy to generalize the result to the case of non-centered monomials. 

In the previous example, we have presented for the case $d=q=1$ the existence of polynomials $t_\myf$ such that the family $\myF_2$ of monomials is estimatable. For general $d$ it turns out the multinomials also exist and they are best described in terms of tensor notation. For a discrete distribution $G= \sum_{i\in [k]} p_i \delta_{\theta_i}$, the $\ell$-th moment tensor is defined as
$$
M_\ell (G)  =:  \sum_{i\in [k]} p_i \theta_i^{\otimes \ell},
$$
where $\otimes$ denotes the tensor product and $\otimes \ell$ in the exponent denotes the tensor power. We also use the notation $\operatorname{sym}(\cdot)$ to denote the symmetrization operation of a tensor. For location Gaussian mixture models, we have the following lemma adapted from \cite[Theorem 5.1]{pereira2022tensor}. 

\begin{lem} \label{lem:denoisedmomenttensor}
For $X\sim \Pb_G$, location Gaussian mixture models with any mixing measure $G$ (that may be continuous or discrete), we have: for any positive integer $\ell$,
$$
M_\ell (G) = \Eb_{G} \sum_{j=0}^{\lfloor \ell/2 \rfloor } A_{\ell,j} (-1)^j \operatorname{sym}\left( X^{\otimes \ell-2j} \otimes \Sigma^{\otimes j} \right) := \Eb_{G} F_\ell(X)  , 
$$
where $A_{\ell,j} =\binom{\ell}{2j} \frac{(2j)!}{j!2^j}$. 
\end{lem}


We now use the above lemma to establish that $\myF_2$ is estimatable. For any $\beta\in [d]^\ell$, denote $\pi_{\ell,i}(\beta)= \# \{j\in [\ell]: \beta^{(j)} = i  \} $. Then $\pi_\ell(\beta)=(\pi_{\ell,1}(\beta),\ldots, \pi_{\ell,d}(\beta)) \in \Omega_{d,\ell}:=\{ \alpha\in \Nb^d: |\alpha| = \ell  \} $. Now consider any $\myf(\theta) = \theta^\alpha $  for some $\alpha \in \Omega_{d,\ell}$ for some $\ell$.   Choose any $\beta\in \pi_\ell^{-1}(\alpha)$ and define
$$
t_\alpha(X) :=  (F_\ell(X))_{\beta}, 
$$
where $(F_\ell(X))_{\beta}$ is the $\beta$-coordinate of the tensor $F_\ell(X)$. Since $F_\ell(X)$ is a symmetric tensor, $(F_\ell(X))_{\beta}$ remains the same for any $\beta\in \pi_\ell^{-1}(\alpha)$ and thus $t_\alpha(X)$ is well-defined. Then for any mixing measure $G$
$$
\Eb_G t_\alpha(X) = (\Eb_G F_\ell(X))_\beta = (M_\ell(G))_{\beta} =  G \theta^\alpha, 
$$
which shows that the space of all monomials is estimatable on the space of all mixing measures. 

Now our estimators with $\myF_2$ from \eqref{eqn:momentmethod} is equivalent to
\begin{equation}
\hat{G}_n \in \argmin_{G'\in G_k(\Theta)} \max_{\alpha\in \Ic_{2k-1}} \left|m_\alpha(G) - \frac{1}{n} \sum_{i\in [n]}t_\alpha(X_i)   \right| =  \argmin_{G'\in G_k(\Theta)} \max_{\ell \in [2k-1] } \left\|M_\ell(G) - \frac{1}{n} \sum_{i\in [n]}F_\ell(X_i)   \right\|_\infty, \label{eqn:tensorestimator}
\end{equation}
where $\|\cdot\|_\infty$ of a tensor is defined to be the largest magnitude of its entries. The above estimator is already studied in \cite[Section 5]{pereira2022tensor} with a small difference being that they use Frobenius norm of tensors instead $\|\cdot\|_\infty$. Interested readers may refer to their paper for computational methods to calculate the optimization problems, but they do not provide a statistical theoretical guarantee for the estimator, which we will discuss as a special case of our general framework. Note the results below can be easily modified to Frobenius norm of tensors.

\begin{lem} \label{lem:tensorcon}
For $X_1,\ldots,X_n \overset{iid}{\sim} \Pb_G$, location Gaussian mixture models on $\Rb^d$ with any mixing measure $G$  (that may be continuous or discrete) on compact $\Theta$, we have: for any $\epsilon>0$, 
\begin{align*}
& \sup_{G} \Pb_{G}\left( \max_{\ell\in [2k-1]}  \left\|M_\ell(G) - \frac{1}{n} \sum_{i\in [n]}F_\ell(X_i)   \right\|_\infty > \epsilon \right) \\
\leq & C(d,k) \exp \left( -C(\Theta,\|\Sigma^{\frac{1}{2}}\|_2, k ) \min\{ n\epsilon^2, (n\epsilon^2)^{\frac{1}{2k-1}} \} \right). 
\end{align*}
Consequently, 
$$
\sup_{G} \Eb_{G} \max_{\ell\in [2k-1]}  \left\|M_\ell(G) - \frac{1}{n} \sum_{i\in [n]}F_\ell(X_i)   \right\|_\infty \leq    C(d,\Theta,\|\Sigma^{\frac{1}{2}}\|_2, k ) n^{-\frac{1}{2}}. 
$$
\end{lem}

\begin{thm}  \label{prop:tensor} 

For $X_1,\ldots,X_n \overset{iid}{\sim} \Pb_{G^*}$, location Gaussian mixture models on $\Rb^d$ with mixing measure $G^*$  on compact $\Theta$.  Let $\hat{G}_n$ be a GMM estimator as in \eqref{eqn:tensorestimator}. 
		\begin{enumerate}[label=(\alph*)]
 \item \label{itema:thm:tensor}
 Then there exist positive constants $C, C'$ and $c$, where their dependence on $ \Theta, k,d, \Sigma$  are suppressed, such that \myred{for $D\in \{W_{2k-1}^{2k-1}, \mbf_{2k-1}\}$}
 \begin{align*}
  \sup_{G^*\in \Gc_k(\Theta)}	\Pb_{G^*} \left(D(G^*,\hat{G}_n)\geq t \right) 
 	\leq   C' \exp\left( -C \min \left\{ nt^{2},  \left(n t^2\right)^{\frac{1}{2k-1}} \right\}  \right), 
 \end{align*}
 and consequently,
 \begin{align*}
 \sup_{G^*\in \Gc_k(\Theta)}	\Eb_{G^*} D(G^*,\hat{G}_n )
 	\leq   c n^{-\frac{1}{2}}. 
 \end{align*}

		\item \label{itemb:prop:tensor}
		Fix any $G_0\in \Gc_k(\Theta)$. Then there exists $r(G_0)$,  $C(G_0)$, $C'(G_0)$ and $c(G_0)$,  where their dependence on $ \Theta, k,k_0,d, \Sigma$  are suppressed, such that 
        \myred{for $D\in \{W_{2d_1-1}^{2d_1-1}, \mbf_{2d_1-1}\}$}
	\begin{align*}
		&\sup_{G^*\in \Gc_k(\Theta): W_1(G_0,G^*)<r(G_0) }	\Pb_{G^*} \left(D(G^*,\hat{G}_n)\geq t \right) \\
		\leq  &  C'(G_0) \exp\left( -C(G_0) \min \left\{ nt^{2},  \left(nt\right)^{\frac{1}{2k-1}} \right\}  \right), 
	\end{align*}
     and consequently,
     \begin{align*}
     	\sup_{G^*\in \Gc_k(\Theta): W_1(G_0,G^*)<r(G_0) }	\Eb_{G^*} D(G^*,\hat{G}_n)
     	\leq   c(G_0) n^{-\frac{1}{2}}. 
     \end{align*}
 \end{enumerate}
\end{thm}
Applying the same argument given in Remark \ref{rem:minimaxrateKS}, we conclude that GMM estimators for high-dimensional location Gaussian mixtures achieve the minimax optimal rate. It is worth mentioning that \cite{doss2020optimal} studies multi-dimensional location Gaussian mixtures using different estimators by projecting it to the univariate case and measures the error by sliced Wasserstein distance. 
}

\section{Pointwise convergence analysis}
\label{sec:pointwiserate}

We have seen in the previous sections how the optimal minimax estimation rate for the mixing measure deteriorates with the overfit index $d_1=k-k_0+1$. In many statistical applications where the data sample can be reasonably assumed to be draw from a \emph{single} unknown distribution, the pointwise convergence rate of the unknown parameters may be more meaningful. We shall show that the family of minimum $\phi$-distance estimator achieves the pointwise optimal rate of convergence under relatively milder conditions.
We consider the setting where the number of support points $k^*$ for the true mixing measure $G^*$ is unknown. 
The estimator consists of the two steps: first, a consistent estimate of $k^*$ will be obtained, and second, a plug-in estimate for $G^*$. Both steps make essential use of the $\Phi$-distance.

\subsection{Estimating the number of mixture components}

For a positive sequence $a_n$, define the following estimator
\begin{equation}
\hat{k}_n:= \inf \left\{\ell \geq 1: \sup_{\myf\in\myF}\left|\hat{G}_n(\ell)  \myf-\frac{1}{n}\sum_{i\in [n]}t_\myf(X_i)\right|\leq a_n  \right\},  \label{eqn:orderestimate}
\end{equation}
with the convention that $\inf \emptyset = \infty$, and recall that $\hat{G}_n(\ell)$ is defined in Section \ref{sec:mindisest}. 
For any discrete distribution $G\in \Gc_k$, let $k(G)$ denote its number of support points and define $b_G$ to be the distance between $G$ and $\Gc_{k(G)-1}$, the set of all discrete measures with fewer supporting atoms than that of $G$, i.e., 
    $$
    b_G:= \inf_{\substack{G'\in \Gc_{k(G)-1} }   } \sup_{\myf\in\myF} |G'\myf-G \myf| .  
    $$
 Since $\Gc_{k(G)-1}$ is compact due to the compactness of $\Theta$, we have $b_G>0$ provided that $\Gc_{k(G)}$ is distinguishable by $\myF$. 
 The following lemma provides a basic template for the design and analysis of the 
 estimate $\hat{k}_n$.

\myred{
\begin{lem} \label{lem:componentnumberestimation}
Consider any discrete measure $G$ on $\Theta$.  Suppose that $\Theta$ is compact, $\Gc_{k(G)}(\Theta)$ is distinguishable by $\myF$, and $\myF$ is estimatable on $\Gc_{k(G)}(\Theta)$.  There holds  
    $$
    \{\hat{k}_n\neq k(G)\}   \subset \left\{ \sup_{\myf\in\myF}\left|G  \myf-\frac{1}{n}\sum_{i\in [n]}t_\myf(X_i)\right| \geq \min\{a_n,b_G-a_n\} \right\}.
    $$ 
\end{lem}
\begin{rem}
Note that the right hand side in the above statement depends on $G$, so the deduced convergence rate result for $\hat{k}_n$ will be pointwise. Moreover, $G$ can be any discrete measure, not necessarily the true mixing measure $G^*$. Since we are "free" to choose both $a_n$ (for the method design), in order to derive a meaningful bound,  one should have $a_n \leq b_G$ asymptotically. To make the set on the right hand side as small as possible, one would ideally choose $a_n=\frac{b_G}{2}$. However, since $b_G$ is generally unknown, such a choice is not possible. 
    One can choose $a_n=o(1)$ to guarantee $a_n \leq  b_G$ asymptotically. It then follows that $\min\{a_n,b_G-a_n\} \asymp a_n$ so we want to choose the rate $a_n$ converging to zero as slow as possible so that the event on the right hand side is small for large $n$.  On the other hand, if $a_n$ converges to $0$ too slowly, then $a_n\leq b_G$ might not hold for small $n$ and thus the result is non-trivial only for large $n$. In summary, the choice for $a_n$ should be $a_n=o(1)$ and the convergence rate to $0$ represents a trade-off between an asymptotic result and a non-asymptotic one. Making the result non-trivial for small $n$ favors choosing a fast decaying sequence $a_n$, while making the result tighter asymptotically favors a slower decaying $a_n$. 
    We will show in the sequel that such a sequence can be chosen to derive optimal pointwise rates of convergence for the mixing measure.

\myred{
The above inclusion conclusion with $G=G^*$ naturally yields a bound on the estimation error probability of $\hat{k}_n$; see Examples \ref{exa:KSpointwise}, \ref{exa:MMDpointwise} and \ref{exa:GMMfixedG_0} ahead for rates for specific examples. One benefit of this result is the absence of an upper bound on $k^*$, which is typically required in the literature (e.g. \cite{manole2021estimating} and the references therein). 
}
    \myeoe
\end{rem}
}

\subsection{Inverse bounds with one argument fixed}
\label{sec:invboufix}

The key and the most technical part for deriving the uniform convergence rate under our general framework is to establish the local inverse bounds \eqref{eqn:localinversebound} and \eqref{eqn:localinverseboundmoment} as shown in Theorem \ref{thm:inversebound}. However if one only intends to establish a pointwise convergence rate for a particular true mixing measure, it   suffices to have an inverse bound with one argument fixed: given $G_0\in \Ec_{k_0}(\Theta)$,
\begin{equation}
\liminf_{ \substack{G\overset{W_1}{\to} G_0\\ G\in \Gc_k(\Theta) }} \frac{\sup_{\myf\in \myF} |G \myf-G_0 \myf|}{W_{2}^{2}(G,G_0)} >0. \label{eqn:localinverseboundfixargument}
\end{equation}
 Such an inverse bound can be established under a suitable strong identifiability condition which is considerably weaker than those required for establishing the uniform inverse bounds presented in the general Theorem \ref{thm:inversebound}.

\begin{definition} \label{def:2ndfixoneargument}
The family $\myF$ is said to be a \textit{$(G_0,k)$ second-order linear independent domain} for $G_0=\sum_{i=1}^{k_0}p_i^0\delta_{\theta_i^0}\in \Ec_{k_0}(\Theta)$ if the following hold: 1) Each $\myf\in \myF$ is second-order continuously differentiable at $\theta_i^0$ for each $i\in [k_0]$; and  2) Consider any integer $\ell_1\in [k_0]$, and $\ell\in [k_0,k]$. Set $m_i=2$ for $i\in [\ell_1]$, $m_i=1$ for $\ell_1<i\leq k_0$ and $m_i=0$ for $k_0<i\leq \ell$. For any distinct $\{\theta_i^0\}_{i=k_0+1}^{\ell}\subset \Theta \setminus \{\theta_i^0\}_{i\in [k_0]}$, the  operators $\{D^\alpha|_{\theta=\theta_i^0}\}_{ 0\leq |\alpha|\leq m_i, i\in [\ell]}$  on $\myF$ are linearly independent, i.e.,  
\begin{subequations}
\begin{align}
\sum_{i=1}^{\ell}\  \sum_{|\alpha|\leq m_i} a_{i\alpha} D^\alpha \myf(\theta_i^0)  = &0,  \quad \forall \myf\in \myF \label{eqn:linearinddomainafixed}\\
\sum_{i\in [\ell]}   a_{i\bm{0}}   = &0, \label{eqn:linearinddomainbfixed}
\end{align}
\end{subequations}
if and only if 
$$
a_{i\alpha}=0, \quad \forall \ 0\leq |\alpha|\leq m_i, \ i\in [\ell].
$$
\end{definition}

It is clear that $\myF$ is $m$-strongly identifiable for $m=2$ implies that $\myF$ is a $(G_0,k)$ second-order linear independent domain for any $k\geq 1$ and $G_0\in \Ec_{k_0}(\Theta)$ with $k_0\in [k]$. 

\begin{rem}
\myred{In principle, one can also have a stronger inverse bound (and upper bound) in terms of moment difference when $G_0$ is fixed, in the spirit of Theorem \ref{thm:inverseboundmoment}. The proof should be similar and we leave the details to interested readers. \myeoe}
\end{rem}

\begin{lem} \label{lem:inverseboundfixoneargument}
Consider a $G_0=\sum_{i=1}^{k_0}p_i^0\delta_{\theta_i^0}\in \Ec_{k_0}(\Theta)$. Suppose that $\myF$ is a $(G_0,k)$ second-order linear independent domain and that $\Theta\subset \Rb^q$ is compact. Then \eqref{eqn:localinverseboundfixargument} holds.
\end{lem}

    The case that $k=k_0$ is known as the inverse bound for the exact-fitted case \cite{ho2016strong,wei2022convergence}.  In this case the local inverse bound \eqref{eqn:localinverseboundfixargument} can be improved: given $G_0\in \Ec_{k_0}(\Theta)$,
\begin{equation}
\liminf_{ \substack{G\overset{W_1}{\to} G_0\\ G\in \Gc_{k_0}(\Theta) }} \frac{\sup_{\myf\in \myF} |G \myf-G_0 \myf|}{W_{1}(G,G_0)} >0. \label{eqn:localinverseboundfixargumentexactfitted}
\end{equation}

\begin{definition} \label{def:1stfixoneargument}
The family $\myF$ is said to be a \textit{$(G_0,k_0)$ first-order linear independent domain} for $G_0=\sum_{i=1}^{k_0}p_i^0\delta_{\theta_i^0}\in \Ec_{k_0}(\Theta)$ if the following hold. 1) Each $\myf\in \myF$ is first-order continuously differentiable at $\theta_i^0$ for each $i\in [k_0]$.  2) The  operators $\{D^\alpha|_{\theta=\theta_i^0}\}_{ 0\leq |\alpha|\leq 1, i\in [k_0]}$  on $\myF$ are linearly independent, i.e.,  
\begin{subequations}
\begin{align}
\sum_{i=1}^{k_0}\  \sum_{|\alpha|\leq 1} a_{i\alpha} D^\alpha \myf(\theta_i^0)  = &0,  \quad \forall \myf\in \myF \label{eqn:linearinddomainafixed1st}\\
\sum_{i\in [k_0]}   a_{i\bm{0}}   = &0, \label{eqn:linearinddomainbfixed1st}
\end{align}
\end{subequations}
if and only if 
$$
a_{i\alpha}=0, \quad \forall \ 0\leq |\alpha|\leq 1, \ i\in [k_0].
$$
\end{definition}

It is clear that if $\myF$ is a $(2d_1-1,k_0,k)$ linear independent domain then $\myF$ is a $(G_0,k_0)$ first-order linear independent domain for any $G_0\in \Ec_{k_0}(\Theta)$. It also follows that $\myF$ is $m$-strongly identifiable for $m=1$ implies that $\myF$ is a $(G_0,k_0)$ first-order linear independent domain for any $G_0\in \Ec_{k_0}(\Theta)$ for any $k_0\geq 1$. 

\begin{lem} \label{lem:inverseboundfixoneargumentexactfitted}
Consider a $G_0=\sum_{i=1}^{k_0}p_i^0\delta_{\theta_i^0}\in \Ec_{k_0}(\Theta)$. Suppose that $\myF$ is a $(G_0,k_0)$ first-order linear independent domain. Then \eqref{eqn:localinverseboundfixargumentexactfitted} holds.
\end{lem}

\begin{rem} 
\label{rem:compactness}
    Lemma \ref{lem:inverseboundfixoneargument} extends the existing results \cite{nguyen2013convergence,ho2016strong}, while Lemma \ref{lem:inverseboundfixoneargumentexactfitted} extends the existing results \cite{ho2016strong,wei2022convergence} to the general $\Phi$-distance.  \myred{Unlike the results of \cite{ho2016strong}, the
pointwise inverse bounds in this section hold when $\Pb_\theta$ is not necessarily absolutely continuous
with respect to the Lebesgue measure.} 
   Note the compactness of $\Theta$ is not required when $k=k_0$, while for the case $k>k_0$ in general compactness assumption is in fact necessary for inverse bounds to hold (see Lemma \ref{lem:compactness}). There are also some relevant inequalities (c.f. \cite[Theorem 2]{nguyen2013convergence}) that hold for a subset of mixing measure satisfying some moment constraints, and unlike the inverse bounds in this paper they hold for all mixing measures on $\Gc_k(\Theta)$; for such inequalities the compactness is not needed. 
    \myeoe
\end{rem}

\begin{lem} 
\label{lem:compactness}
Suppose that $\Theta=\Rb^q$ and the function class $\myF$ is uniformly bounded, i.e.  $\sup_{\myf\in \myF} \sup_{\theta\in \Theta}|\myf(\theta)|<\infty$.  Consider $G_0\in \Ec_{k_0}(\Theta)$ and $k>k_0$. Then for any $r>0$, 
\begin{equation}
\liminf_{ \substack{G\overset{W_1}{\to} G_0\\ G\in \Gc_k(\Theta) }} \frac{\sup_{\myf\in \myF} |G \myf-G_0 \myf|}{W_{r}^{r}(G,G_0)} =0. 
\end{equation}
\end{lem}


\subsection{Optimal pointwise convergence for mixing measures}

Let $\hat{k}_n$ be any estimator for the number of mixture components. In this subsection we study the plug-in estimate $\hat{G}_n(\hat{k}_n)$, a minimum $\myF$-distance estimator combining with the estimated number of mixture component $\hat{k}_n$. We first
state a general theorem and then specialize it to the examples considered in Section \ref{sec:examples}. \myred{The main message is that, to improve the convergence rates, one should perform model selection first, and then do the parameter estimation. }

\begin{thm}
\label{thm:fixedG_0estimator}
    Consider a $G_0=\sum_{i=1}^{k_0}p_i^0\delta_{\theta_i^0}\in \Ec_{k_0}(\Theta)$. Suppose that $\Theta$ is compact, that $\Gc_{k_0}(\Theta)$ is distinguishable by $\myF$ and that $\myF$ is estimatable on $\Gc_{k_0}(\Theta)$. Suppose further that inverse bound \eqref{eqn:localinverseboundfixargumentexactfitted} holds for $G_0\in \Ec_{k_0}(\Theta)$.  
    \begin{enumerate}[label=(\alph*)]
        \item 
        \label{thm:fixedG_0estimatora}
        Consider any estimator $\hat{k}_n$ for the number of mixture components.  Then there exist positive constants $\epsilon_1, \epsilon'_1, C(G_0)>0$ that depend on $G_0$, $ \Theta$ and $\myF$, such that for any $t>0$,
\begin{align*}
    \{W_{1}(G_0,\hat{G}_n(\hat{k}_n)) \geq t \} 
    \subset &  \left\{ \sup_{\myf\in\myF}\left|G_0  \myf-\frac{1}{n}\sum_{i\in [n]}t_\myf(X_i)\right| \geq \min\left\{\epsilon_1 t,\epsilon_1\right\} \right\} \bigcup \{\hat{k}_n\neq k_0\},
\end{align*}
and 
\begin{align*}
& \Eb_{G^*} W_{1}(G_0,\hat{G}_n(\hat{k}_n)) \\
\leq & C(G_0)\Eb_{G^*} \sup_{\myf\in\myF}\left|G_0  \myf-\frac{1}{n}\sum_{i\in [n]}t_\myf(X_i)\right| +\diam(\Theta) \Pb_{G^*} \left( \sup_{\myf\in\myF}\left|G_0  \myf-\frac{1}{n}\sum_{i\in [n]}t_\myf(X_i)\right| \geq \epsilon'_1  \right) \\
& \quad +\diam(\Theta) \Pb_{G^*} \left( \hat{k}_n\neq k_0 \right).
\end{align*}

        \item
         \label{thm:fixedG_0estimatorb}
        Let $\hat{k}_n$ be the estimator defined in \eqref{eqn:orderestimate}. 
        Then there exist positive constants $\epsilon_0, \epsilon'_0>0$ that depend on $G_0$, $ \Theta$ and $\myF$, such that for any $t>0$,
\begin{align*}
    \{W_{1}(G_0,\hat{G}_n(\hat{k}_n)) \geq t \} 
    \subset & \left\{ \sup_{\myf\in\myF}\left|G_0  \myf-\frac{1}{n}\sum_{i\in [n]}t_\myf(X_i)\right| \geq \min\left\{\epsilon_0 t, a_n,\epsilon_0-a_n\right\} \right\},
\end{align*}
and 
\begin{align*}
& \Eb_{G^*} W_{1}(G_0,\hat{G}_n(\hat{k}_n)) \\
\leq & C(G_0) \Eb_{G^*} \sup_{\myf\in\myF}\left|G_0  \myf-\frac{1}{n}\sum_{i\in [n]}t_\myf(X_i)\right| + \\
& \diam(\Theta) \Pb_{G^*} \left( \sup_{\myf\in\myF}\left|G_0  \myf-\frac{1}{n}\sum_{i\in [n]}t_\myf(X_i)\right| \geq \min\{a_n,\epsilon'_0-a_n\}  \right).
\end{align*}
    \end{enumerate}

\end{thm}
\begin{rem}
It is emphasized that the above theorem is stated for any $G_0$ for which the inverse bound \eqref{eqn:localinverseboundfixargumentexactfitted} holds, and $G_0$ is not necessarily the true mixing measure $G^*$. Thus, the theorem is applicable to deriving the rates of convergence for mixing measures in the setting of model mis-specification.
Moreover, it applies to any estimator $\hat{k}_n$, not just the one studied in Part \ref{thm:fixedG_0estimatorb}.  Estimating the number of mixture component (or the order of the mixture) is an important question that attracts continued attention (cf. e.g., recent papers \cite{guha2021posterior,manole2021estimating,cai2020power,cai2021finite} and references therein). 
\myeoe
\end{rem}

It is worth to point out that, unlike the minimax rate setting, the pointwise convergence rate result Theorem \ref{thm:fixedG_0estimator} does not require the knowledge of an upper bound $k$ for the order $k^*$ of the true mixing measure $G^*$. 

\begin{exa}[Minimum KS-distance estimator combined with $\hat{k}_n$]
\label{exa:KSpointwise}
Consider the example studied in Section \ref{sec:minimumdistance estimator}. As in Theorem \ref{thm:convergencerateminimumdistanceestimator}, suppose that $\Theta$ is compact, and suppose  that the mixture model is identifiable on $\Gc_k(\Theta)$. 
Let $a_n=c_1 \sqrt{\frac{\ln n}{n}} $ for some constant $c_1$ and let $\hat{k}_n$ be the estimator defined in \eqref{eqn:orderestimate}.   Applying Lemma \ref{lem:componentnumberestimation} with $G=G^*$, we then have 
$$
\Pb_{G^*}(\hat{k}_n\neq k(G^*)) \leq  \Pb_{G^*}(\KS(\hat{\Pb}_n,\Pb_{G^*})\geq \min\{a_n,b_{G^*}-a_n\}   ) \leq c_2(G^*,c_1) \Pb_{G^*}(\KS(\hat{\Pb}_n,\Pb_{G^*})\geq a_n  ),
$$
where $c_2(G^*,c_1)$ is a constant that depends on $G^*,d,c_1, b_{G^*}$ and the model ($\Theta$, $\{\Pb_{\theta}\}$ etc). By \cite[Lemma 4.1]{naaman2021tight}, we then have 
$$
\Pb_{G^*}(\hat{k}_n\neq k(G^*))  \leq C(G^*,c_1) d(n+1)e^{-2n a_n^2} = C(G^*,c_1) d(n+1) n^{-2c_1^2},
$$
which converges to $0$ with $c_1> \frac{1}{\sqrt{2}}$, and $C(G^*,c_1)$ is a constant that depends on $G^*,d,c_1$, $b_{G^*}$ and the model ($\Theta$, $\{\Pb_{\theta}\}$ etc).

Suppose additionally that $\myF_0$ is $1$-strongly identifiable  as in Section \ref{sec:minimumdistance estimator}. Then \eqref{eqn:localinverseboundfixargumentexactfitted} holds for any $G_0\in \Ec_{k_0}(\Theta)$ for any $k_0$. Applying Theorem \ref{thm:fixedG_0estimator} with $G_0=G^*$, we obtain:
\begin{align*}
& \Eb_{G^*} W_{1}(G^*,\hat{G}_n(\hat{k}_n)) \\
\leq & \Eb_{G^*} \KS(\hat{\Pb}_n,\Pb_{G^*}) +\diam(\Theta) \Pb_{G^*} \left( \KS(\hat{\Pb}_n,\Pb_{G^*}) \geq \min\{a_n,\epsilon'_0-a_n\}  \right) \\
\leq & C(G^*,c_1) (n^{-\frac{1}{2}}+ (n+1)n^{-2c_1^2}) \\
\leq & C(G^*,c_1) n^{-\frac{1}{2}},
\end{align*}
where the second inequality follows from Lemma \ref{lem:KSempiricalprocess} and \cite[Lemma 4.1]{naaman2021tight} with $C(G^*,c_1)$ is a constant that depends on $G^*,d,c_1$, $b_{G^*}$ and the model ($\Theta$, $\{\Pb_{\theta}\}$ etc), and the last step follows by choosing $c_1\geq \frac{\sqrt{3}}{2}$. The convergence rate of the estimator $\hat{G}_n(\hat{k}_n)$ in this example under the setting of univariate case $q=d=1$  was firstly studied in \cite[Theorem 4.1]{heinrich2018strong} (with $a_n=n^{\frac{1}{2}+\kappa}$ for some $\kappa>0$). Note that to establish pointwise convergence rate above we do not require the knowledge of an upper bound $k$ for $k^*$. 
Despite the slow uniform rate  $n^{-\frac{1}{2k-1}}$ or $n^{-\frac{1}{2(2d_1-1)}}$ with $d_1=k-k_0+1$ discussed in Remark \ref{rem:minimaxrateKS}, the pointwise convergence rate can be much better ---- in this example $n^{-\frac{1}{2}}$ in particular. \myeoe
\end{exa}

\begin{exa}[Minimum MMD estimator combined with $\hat{k}_n$]
\label{exa:MMDpointwise}
Consider the example studied in Section \ref{sec:MMD}. As in Theorem \ref{thm:MMDuniformrate}, suppose that $\Theta$ is compact and that the map $\mu: \Mc_b(\Xf,\Xc)\to \Hc$ is injective. 
Let $a_n=c_1 \sqrt{\frac{\ln n}{n}} $ for some constant $c_1$ and let $\hat{k}_n$ be the estimator defined in \eqref{eqn:orderestimate}. Again applying Lemma \ref{lem:componentnumberestimation} with $G=G^*$, we then have 
$$
\Pb_{G^*}(\hat{k}_n\neq k(G^*)) \leq  \Pb_{G^*}(\MMD(\hat{\Pb}_n,\Pb_{G^*})\geq \min\{a_n,b_{G^*}-a_n\}   ) \leq c_2(G^*,c_1) \Pb_{G^*}(\MMD(\hat{\Pb}_n,\Pb_{G^*})\geq a_n  ),
$$
where $c_2(G^*,c_1)$ is a constant that depends on $G^*,d,c_1, b_{G^*}$ and the model ($\Theta$, $\{\Pb_{\theta}\}$, $\ker(\cdot,\cdot)$ etc). By Lemma \ref{lem:empiricalprocessMMDconcentration} below, we then have 
\begin{align*}
\Pb_{G^*}(\hat{k}_n\neq k(G^*)) \leq & 
 c_2(G^*,c_1) 2\exp\left(-\frac{n(a_n-\frac{2\|\ker\|_\infty}{\sqrt{n}})^2}{2 \|\ker\|^2_\infty}\right) \\
= & c_2(G^*,c_1) 2\exp\left(-\frac{(c_1\sqrt{\ln n}-2\|\ker\|_\infty)^2}{2 \|\ker\|^2_\infty}\right)\\
\leq & C(G^*,c_1) n^{-\frac{c_1^2}{8\|\ker\|_\infty^2}}, \numberthis  \label{eqn:exammdcontemp}
\end{align*}
where $C(G^*,c_1)$ is a constant that depends on $G^*,d,c_1$, $b_{G^*}$ and the model ($\Theta$, $\{\Pb_{\theta}\}$, $\ker(\cdot,\cdot)$ etc). 

Suppose additionally that $\{p(x \mid \theta)\}_{\theta\in \Theta}$ satisfies Assumption $A2(2k-1)$ and $\{p(x \mid \theta)\}_{\theta\in \Theta}$ is $1$-strongly identifiable as in Section \ref{sec:MMD}. By Remark \ref{strongidentifiability} and Lemma \ref{lem:inverseboundfixoneargumentexactfitted}, \eqref{eqn:localinverseboundfixargumentexactfitted} holds for any $G_0\in \Ec_{k_0}(\Theta)$ for any $k_0\leq k$. Applying Theorem \ref{thm:fixedG_0estimator} with $G_0=G^*$, we obtain:
\begin{align*}
& \Eb_{G^*} W_{1}(G^*,\hat{G}_n(\hat{k}_n)) \\
\leq & \Eb_{G^*} \MMD(\hat{\Pb}_n,\Pb_{G^*}) +\diam(\Theta) \Pb_{G^*} \left( \MMD(\hat{\Pb}_n,\Pb_{G^*}) \geq \min\{a_n,\epsilon'_0-a_n\}  \right) \\
\leq & C(G^*,c_1) (n^{-\frac{1}{2}}+ n^{-\frac{c_1^2}{8\|\ker\|_\infty^2}}) \\
\leq & C(G^*,c_1) n^{-\frac{1}{2}},
\end{align*}
where the second inequality follows from Lemma \ref{lem:empiricalprocessMMDconcentration} and $\eqref{eqn:exammdcontemp}$ with $C(G^*,c_1)$ is a constant that depends on $G^*,d,c_1$, $b_{G^*}$ and the model ($\Theta$, $\{\Pb_{\theta}\}$, $\ker(\cdot,\cdot)$ etc), and the last step follows by choosing $c_1\geq 2\|\ker\|_\infty$. Note that to establish pointwise convergence rate above we do not require the knowledge of an upper bound $k$ for $k^*$.
Despite the slow uniform rate  $n^{-\frac{1}{2k-1}}$ or $n^{-\frac{1}{2(2d_1-1)}}$ with $d_1=k-k_0+1$ established in Theorem \ref{thm:MMDuniformrate}, the pointwise convergence rate can be much better---in this example $n^{-\frac{1}{2}}$ in particular. \myeoe
\end{exa}

\begin{lem}
\label{lem:empiricalprocessMMDconcentration}
Consider a measurable bounded kernel $\ker(\cdot,\cdot)$.  Then for $\epsilon>0$,
    $$ \Pb  \left( \MMD(\Pb, \hat{\Pb}_n) \geq \frac{2\|\ker\|_\infty}{\sqrt{n}} + \epsilon \right)
    \leq 2\exp\left(-\frac{n\epsilon^2}{2 \|\ker\|^2_\infty}\right) , $$
    where the random variables $\{X_i\}_{i\in [n]} \overset{\iidtext}{\sim} \Pb$.
\end{lem}

For pointwise convergence rate for minimum GMM estimators we do need to assume the upper bound $k$ of $k^*$. In fact, the definition of the function class $\myF_2$ already involve $k$. 

\begin{exa}[Minimum GMM estimator combined with $\hat{k}_n$]
\label{exa:GMMfixedG_0}
Consider the example studied in Section \ref{sec:momentdeviation}. Suppose that $\Theta$ is compact and that the mixture model is univariate with $p(x \mid \theta)$ belonging to NEF-QVF. Assume $k$ is an upper bound for $k^*$. Let $a_n=c_1 \sqrt{\frac{\ln n}{n}} $ for some constant $c_1$ and let $\hat{k}_n$ be the estimator defined in \eqref{eqn:orderestimate}. Applying Lemma \ref{lem:componentnumberestimation} with $G=G^*$, we then have 
\begin{align*}
\Pb_{G^*}(\hat{k}_n\neq k(G^*)) \leq & \Pb_{G^*}\left( \max_{j\in [2k-1]}|\bar{t}_j(\theta_0)-\Eb_{G^*} t_j(X \mid \theta_0)| \geq \min\{a_n,b_{G^*}-a_n\}   \right) \\
\leq & c_2(G^*,c_1) \Pb_{G^*}\left(\max_{j\in [2k-1]}|\bar{t}_j(\theta_0)-\Eb_{G^*} t_j(X \mid \theta_0)|\geq a_n  \right),
\end{align*}
where $c_2(G^*,c_1)$ is a constant that depends on $G^*,d,c_1, b_{G^*}$ and the model ($\Theta$, $\{\Pb_{\theta}\}$, $\ker(\cdot,\cdot)$ etc). By Lemma \ref{lem:concentrationboualljNEFQVF}, we then have 
\begin{align*}
\Pb_{G^*}(\hat{k}_n\neq k(G^*)) \leq & \Pb_{G^*}\left(\max_{j\in [2k-1]}|\bar{t}_j(\theta_0)-\Eb_{G^*} t_j(X \mid \theta_0)| \geq \min\{a_n,b_{G^*}-a_n\}   \right) \\
\leq & c_2(G^*,c_1) e^2 (2k-1) \exp\left( -C \min \left\{ n a_n^2,  (n a_n)^{\frac{1}{2k-1}} \right\}  \right) \\
\leq & c_3(G^*,c_1) \exp\left(- C n a_n^2  \right)\\
= & c_3(G^*,c_1) n^{-C c_1^2}, \numberthis  \label{eqn:exammdcontempmoment}
\end{align*}
where $C$,$c_3(G^*,c_1)$ are positive constants that depends on $d$, $b_{G^*}$ and the model ($\Theta$, $\{\Pb_{\theta}\}$, $\ker(\cdot,\cdot)$ etc), and $c_3(G^*,c_1)$ additionally depends on $G^*$ and $c_1$.

 By Lemma \ref{lem:momentbound} and Lemma \ref{lem:inverseboundfixoneargumentexactfitted}, \eqref{eqn:localinverseboundfixargumentexactfitted} holds for any $G_0\in \Ec_{k_0}(\Theta)$ for any $k_0\leq k$. Applying Theorem \ref{thm:fixedG_0estimator} with $G_0=G^*$, we obtain:
\begin{align*}
& \Eb_{G^*} W_{1}(G^*,\hat{G}_n(\hat{k}_n)) \\
\leq & \Eb_{G^*} \max_{j\in [2k-1]}|\bar{t}_j(\theta_0)-\Eb_{G^*} t_j(X|\theta_0)| +\\
& \diam(\Theta) \Pb_{G^*} \left( \max_{j\in [2k-1]}|\bar{t}_j(\theta_0)-\Eb_{G^*} t_j(X|\theta_0)| \geq \min\{a_n,\epsilon'_0-a_n\}  \right) \\
\leq & c_4(G^*,c_1) \left(n^{-\frac{1}{2}}+ n^{-C c_1^2} \right) \\
\leq & c_4(G^*,c_1) n^{-\frac{1}{2}},
\end{align*}
where the second inequality follows from Lemma \ref{lem:concentrationboualljNEFQVF} and $\eqref{eqn:exammdcontemp}$ with $c_4(G^*,c_1)$ a constant that depends on $G^*,d,c_1$, $b_{G^*}$ and the model ($\Theta$, $\{\Pb_{\theta}\}$, $\ker(\cdot,\cdot)$ etc), and the last step follows by choosing $c_1\geq \frac{1}{\sqrt{2C}}$. 
Despite the slow uniform rate  $n^{-\frac{1}{2k-1}}$ or $n^{-\frac{1}{2(2d_1-1)}}$ with $d_1=k-k_0+1$ established in Theorem \ref{thm:MMDuniformrate}, again the pointwise convergence rate can be much better --- in this example $n^{-\frac{1}{2}}$ in particular. \myeoe
\end{exa}

\section{Discussion}
\label{sec:discussion}

In this paper we proposed a general estimation framework for finite mixing measures and analyzed the convergence rates. While the minimum $\Phi$-distance estimation framework is very general, as demonstrated in this paper,  we note that there are certain minimum distance or divergence-type estimators which do not belong to our framework \cite{ho2020robust, edelman1988estimation, jana2022optimal}. 

There are a number of interesting open questions that are worth exploring. A direction is to generalize our distance to more a general form, e.g., one which accommodates the $f$-divergence. Another direction is to remove the assumption of a known upper bound for the true number of mixture components. 
One may also further investigate different choices of test function classes $\myF$ and possibly find an optimal one in a certain sense (the one with smallest cardinality for instance). \myred{One may also investigate the dependence of the constant in the inverse bounds on different parameters, say $k$, $d$ and $q$ (it is worth to mention \cite{wei2022convergence} managed to derive the dependence of the constant in the inverse bounds on $m$ for general mixtures of $m$-product distribution); see also Remark \ref{rem:infgeocom} for some related literature. Finally, one attractive property of minimum MMD estimators is that they can potentially be applied to mixture distributions that are non-Euclidean, and thus one can explore this direction to study mixtures on non-Euclidean space, say mixtures of von Mises-Fisher distributions \cite{banerjee2005clustering} or mixture of general product distributions \cite[Section 7.4]{wei2022convergence}. }

One of the key components of the theory in such an effort is the development of inverse bounds that go beyond the $\sup$ norm associated with the $\Phi$ function class. In the following we describe some relevant results that may be of independent interest.

\subsection{Inverse bounds: beyond $\sup$ norm }
\label{sec:beyondsup}

In the previous sections in the paper we have considered minimum distance estimators where the distance between two mixing measures is given by $\sup_{\myf\in \myF} |G \myf-H \myf|$. The particular form $\sup_{\myf\in \myF} |G \myf-H \myf|$ taken is due to its generality but there are other alternatives. Suppose there is a measure $\mathscr{T}$ on $\myF$. Then one alternative is $\int_{\myF} |G \myf-H \myf| d\mathscr{T} $, the average of the absolute difference between the two mixing measure applying to each member $\myf$. Similar to Definition \ref{def:distin}, we have the following definition of distinguishability. 

\begin{definition}
	$\Gc_k(\Theta)$ is said to be \textit{distinguishable} by $(\myF,\mathscr{T})$  if for any $G \neq H\in \Gc_k(\Theta)$, $\int_{\myF} |G \myf-H \myf| d\mathscr{T}>0$. 
\end{definition}
If $\Gc_k(\Theta)$ is distinguishable by $(\myF,\mathscr{T})$, then it is easy to see that $\int_{\myF} |G \myf-H \myf| d\mathscr{T}$ is a distance on $\Gc_k(\Theta)$.

\begin{exa}[Total variational distance between mixtures]
\label{exa:intIPM}
Assume that $\{\Pb_{\theta}\}_{\theta\in \Theta}$ has density $\{p(x \mid \theta)\}_{\theta\in \Theta}$ w.r.t. a dominating measure $\lambda$ on $(\Xf,\Xc)$.
Consider $\myF_3=\left\{ \theta \mapsto p(x \mid \theta) | x\in \Xf \right\}$. For each $x\in \Xf$, $p(x \mid \theta)$ is a function of $\theta$. Note that $\lambda$ on $(\Xf,\Xc)$ induces a measure $\mathscr{T}$ on $\myF$. Then
\begin{equation}
\int_{\myF} |G \myf-H \myf| d\mathscr{T} = \int_{\Xf} \left|p_G(x)-p_H(x)\right| d\lambda=2 V(\Pb_G,\Pb_H), \label{eqn:integraldistanceint}
\end{equation}
twice of the total variation distance between the mixtures $\Pb_G$ and $\Pb_H$. 
   \myeoe
\end{exa}

Next we discuss the corresponding inverse bounds. The local inverse bound becomes:
\begin{equation}
\liminf_{ \substack{G,H\overset{W_1}{\to} G_0\\ G\neq H\in \Gc_k(\Theta) }} \frac{\int_{ \myF} |G \myf-H \myf| d\mathscr{T} }{W_{2d_1-1}^{2d_1-1}(G,H)} >0. \label{eqn:localinverseboundint}
\end{equation}

\begin{definition} \label{def:linearindependentdomainint}
The family $(\myF,\mathscr{T})$ is said to be a \textit{$(m,k_0,k)$ linear independent domain} if the following hold. 1) For $\Ts$-a.e. $\myf\in \myF$, $\myf$ is $m$-th order continuously differentiable on $\Theta$.  2) Consider any integer $\ell\in [k_0,2k-k_0]$, and any vector $(m_1,m_2,\ldots,m_\ell)$ such that $1\leq m_i\leq m+1$ for $i\in [\ell]$ and $\sum_{i=1}^\ell m_i\in [2k_0, 2k]$.
For any distinct $\{\theta_i\}_{i\in [\ell] }\subset \Theta$, the  operators $\{D^\alpha|_{\theta=\theta_i}\}_{ 0\leq |\alpha|< m_i, i\in [\ell]}$  on $\myF$ are linear independent, i.e.,  
\begin{subequations}
\begin{align}
\sum_{i=1}^{\ell}\  \sum_{|\alpha|\leq m_i-1} a_{i\alpha} D^\alpha \myf(\theta_i)  = &0,  \quad \mathscr{T}-a.e. \ \  \myf\in \myF \label{eqn:linearinddomainaint}\\
\sum_{i\in [\ell]}   a_{i\bm{0}}   = &0, \label{eqn:linearinddomainbint}
\end{align}
\end{subequations}
if and only if 
$$
a_{i\alpha}=0, \quad \forall \ 0\leq |\alpha|< m_i, \ i\in [\ell].
$$
\end{definition}

\begin{thm} 
\label{thm:inverseboundint}
Consider $\Theta\subset \Rb^q$ is compact.
\begin{enumerate}[label=(\alph*)]

\item \label{thm:inverseboundaint}
If that $(\myF,\mathscr{T})$ is a $(2d_1-1,k_0,k)$ linear independent domain, then \eqref{eqn:localinverseboundint} holds for any $G_0\in \Ec_{k_0}(\Theta)$.
\item \label{thm:inverseboundbint}
If that $(\myF,\mathscr{T})$ is a $(2k-1,1,k)$ linear independent domain, then \eqref{eqn:localinverseboundint} holds for any $G_0\in \Gc_{k_0}(\Theta)$ for any $k_0\in [k]$. 
\end{enumerate}
\end{thm}
The proof of Theorem \ref{thm:inverseboundint} is a simple and straightforward modification of the proof of Theorem \ref{thm:inversebound} and is thus omitted. Note also an entirely analogous change from $\sup_{\myf\in \myF} |G \myf-H \myf|$ to $\int_{ \myF} |G \myf-H \myf| d\mathscr{T}$  for inverse bounds with one argument fixed presented in Section \ref{sec:invboufix} can be carried out and is omitted in this paper. Next we apply the above theorem to the total variational distance presented in Example \ref{exa:intIPM}, for which Definition \ref{def:linearindependentdomainint} specializes to Definition  \ref{def:linearindependent}. 

\begin{thm} 
\label{thm:inverseboundintTV}
Consider $\Theta\subset \Rb^q$ is compact.
\begin{enumerate}[label=(\alph*)]

\item 
If $\{p(x \mid \theta)\}_{\theta\in \Theta}$ is a $(2d_1-1,k_0,k)$ linear independent, then it holds for any $G_0\in \Ec_{k_0}(\Theta)$:
\begin{equation}
\liminf_{ \substack{G,H\overset{W_1}{\to} G_0\\ G\neq H\in \Gc_k(\Theta) }} \frac{V(\Pb_G,\Pb_H) }{W_{2d_1-1}^{2d_1-1}(G,H)} >0.  \label{eqn:localinvbouTV}
\end{equation}
\item 
If $\{p(x \mid \theta)\}_{\theta\in \Theta}$ is $(2k-1,1,k)$ linear independent, then \eqref{eqn:localinvbouTV} holds for any $G_0\in \Ec_{k_0}(\Theta)$ for any $k_0\in [k]$. Moreover, if the mixture model is identifiable, then it holds: 
\begin{equation}
	\inf_{  G\neq H\in \Gc_k(\Theta) } \frac{V(\Pb_G,\Pb_H)}{W_{2k-1}^{2k-1}(G,H)} >0. \label{eqn:globalinverseboundtv}
\end{equation}
\end{enumerate}
\end{thm}

The inverse bound \eqref{eqn:localinvbouTV} has been studied and used to establish convergence rate for parameters in the literature for Bayesian and likelihood-based methods \cite{nguyen2013convergence, wei2022convergence}. Here we obtain the results as a special example of the general result Theorem  \ref{thm:inverseboundint}. Note one could also apply Theorem \ref{thm:inversebound} \ref{thm:inversebounda} with $\myF = \myF_4=\{\Pb_\theta(B)| B\in \Xc \}$ (as in Example \ref{exa:divationbetweenmixturedistributions}) to establish \eqref{eqn:localinvbouTV}, but now the assumption $\myF_4$ is a linear independent domain is relatively more difficult to work with since $\myF_4$ is indexed by all measurable sets from $\Xc$. This specific example demonstrates one instance in which using $\int_{\myF} |G \myf-H \myf| d\mathscr{T} $ is preferable to $\sup_{\myf\in \myF} |G \myf-H \myf|$. \myred{Note when $\Xf=\Rb^d$, since KS distance is a lower bound for total variation distance, \eqref{eqn:localinvbouTV} and \eqref{eqn:globalinverseboundtv} may also be deduced from results in Section \ref{sec:minimumdistance estimator}.}

\myred{
\subsection{Mixture of multinomials}

To demonstrate the novelty of Lemma \ref{lem:momentbound}, we consider the model of mixture of multinomial distributions. 

\begin{exa}[Inverse bound for mixture of multinomials] 
\label{exa:mixmuldis}
A $q$-dimensional multinomial distribution with parameter $N\in \Zb_{\geq 1}$, the set of positive integers, and parameter $\theta \in \Theta:= \{\theta\in \Rb^q| \sum_{i=1}^q \theta^{(i)}\leq 1,  \theta^{(i)}\geq 0, \forall i\} $ has the probability mass function (p.m.f.): $\forall  x\in \Ic_N$,
\begin{align*}
p(x|\theta,N) 
= &\binom{N}{x^{(1)},\ldots,x^{(q)},x^{(q+1)}} \prod_{j=1}^{q+1} (\theta^{(j)})^{x^{(j)}}, \numberthis \label{eqn:mulnomden} 
\end{align*}
where $\theta^{(q+1)}:=1-\sum_{i=1}^q \theta^{(i)} $ and $y^{(q+1)}:=N-\sum_{i=1}^q y^{(i)} $. We denote the multinomial distribution with probability mass function \eqref{eqn:mulnomden} by $\Mul(N,\theta)$. Note that when $q=1$, it reduces to the binomial distribution.

Consider $k_0=1$ and $m=2k-1$. Consider any integer $\ell\in [k_0,2k-k_0]$, and any vector $(m_1,m_2,\ldots,m_\ell)$ such that $1\leq m_i\leq m+1$ for $i\in [\ell]$ and $\sum_{i=1}^\ell m_i\in [2k_0, 2k]$.
For any distinct $\{\theta_i\}_{i\in [\ell] }\subset \Theta$, the functions $\{\frac{\partial^\alpha p}{\partial \theta^{\alpha}}(x \mid \theta_i)\}_{ 0\leq |\alpha|< m_i, i\in [\ell]}$  are linear independent, i.e.,  
\begin{subequations}
\begin{align*}
\sum_{i=1}^{\ell}\  \sum_{|\alpha|\leq m_i-1} a_{i\alpha} \frac{\partial^\alpha p}{\partial \theta^{\alpha}}(x \mid \theta_i, N) = &0,  \quad \forall x\in \Ic_N,  \\ 
\sum_{i\in [\ell]}   a_{i\bm{0}}   = &0. 
\end{align*}


\end{subequations}
\label{sec:mixofmul}

Since $\sp\left(\{p(x|\theta,N,s)\}_{x\in \Ic_{N}}\right)$, viewing as functions of $\theta$, is all multinomials of degree at most $N$. The above linear system is equivalent to: for any multinomial $P(\theta)$ of degree at most $N$,
\begin{align*}
\sum_{i=1}^{\ell}\  \sum_{|\alpha|\leq m_i-1} a_{i\alpha} \frac{\partial^\alpha P}{\partial \theta^{\alpha}}( \theta_i) = &0  . 
\end{align*}
By Lemma \ref{lem:linearalgebra} \ref{lem:linearalgebraa}, when $N\geq 2k-1$, we have
$$
a_{i\alpha}=0, \quad \forall \ 0\leq |\alpha|< m_i, \ i\in [\ell].
$$
That is, $\{p(x \mid \theta)\}_{\theta\in \Theta}$ is a $(2k-1,1,k)$ linear independent and thus by Theorem \ref{thm:inverseboundintTV}, \eqref{eqn:localinvbouTV} holds for any $G_0\in \Ec_{k_0}(\Theta)$ for any $k_0\in [k]$. Moreover, when $N\geq 2k-1$, the mixture of multinomial distributions is identifiable, which yields by Lemma \ref{lem:localtoglobal} the following: 
\begin{equation}
	\inf_{  G\neq H\in \Gc_k(\Theta) } \frac{V(\Pb_G,\Pb_H)}{W_{2k-1}^{2k-1}(G,H)} >0. \label{eqn:globalinverseboundtvmul}
\end{equation}
Since the mixture of multinomial distributions is not identifiable when $N<2k-1$, it follows that \eqref{eqn:globalinverseboundtvmul} does not hold when $N<2k-1$ for mixture of multinomial distributions. As a result, the inverse bound \eqref{eqn:globalinverseboundtvmul} holds if and only if $N\geq 2k-1$. 

As a comparison, \cite[Proposition 1 and Corollary 1]{Manole2021-nw} established that mixture of binomial distribution (special case of our case with $q=1$) satisfies the $m$-strongly identifiability as in Definition \ref{def:strideden} if and only if $N\geq (m+1)k-1$, and then use the $m$-strongly identifiability to establish the inverse bounds. Note that the inverse bounds are what matter in the analysis of the convergence rates, not the sufficient condition $m$-strong identifiable. As one of the key contributions of our paper, a better sufficient condition to guarantee inverse bounds is Definition \ref{def:linearindependent} instead of the $m$-strong identifiability. Indeed, as shown above, as long as $N\geq 2k-1$, the weaker sufficient condition holds and thus inverse bounds hold, which significantly improves the previous results when $m>1$. Our results also hold for mixture of multinomial distributions (any $q$) beyond mixture of binomial distributions ($q=1$). In fact, \cite[Corollary 1]{Manole2021-nw} claims mixture of multinomial distributions is $m$-strongly identifiable when $N\geq 3k-1$ but there is an error in their proof. Our result presented herein outperforms their claimed conclusion by establishing the inverse bound if and only if $N\geq 2k-1$. 
\myeoe

\end{exa}

}

\section*{Acknowledgements}
\myred{We thank Pierre Alquier for bringing the paper \cite{cherief2022finite} to our attention. We thank Dat Do for for bring \cite{pereira2022tensor} to our attention. We also want to thank the anonymous referees and the associate editors for various suggestions and comments that significantly improve our manuscript.}
Yun Wei would like to acknowledge partial
funding from SAMSI and
NSF DMS
17-13012.
 Long Nguyen was partially supported by the NSF Grant DMS-2015361 and a research gift from Wells Fargo. 

\bibliography{reference}{}
	\bibliographystyle{plain}

\appendix


\section{Proofs for Section \ref{sec:unifiedframework}}

\subsection{Proof of Theorem \ref{thm:minimax} \ref{item:thm:minimaxa}}
\label{sec:proofofminimaxtheorem}

\begin{lem}\label{lem:discretedistributions}
    \begin{enumerate}[label=(\alph*)]
    \item
    Consider any $G,H\in \Gc_k(\Rb^q)$. If $\mbf_{2k-1}(G)=\mbf_{2k-1}(H)$, then $G=H$.
    \item \label{itemb:lem:discretedistributions}
    For any $G\in \Ec_k(\Rb)$, there exist infinitely many $H\in \Ec_k(\Rb)$ such that $\mbf_{2k-2}(G)=\mbf_{2k-2}(H)$. Consider any $\psi\in \Rb^q$. For any $G=\sum_{i\in [k]}p_i\delta_{\theta_i}\in \Ec_k(\Rb^q)$ with $\theta_i\in \sp(\psi)$, there exist infinitely many $H\in \Ec_k(\Rb^q)$ with supporting points in  $\sp(\psi)$, such that $\mbf_{2k-2}(G)=\mbf_{2k-2}(H)$.
    \end{enumerate}
\end{lem}
\begin{proof}
(a) Consider $X\sim G$ and $Y\sim H$. Then for any $b\in \Rb^q$,
    $$
    \sum_{|\alpha|=i} \binom{i}{\alpha} m_\alpha(G) b^\alpha = \sum_{|\alpha|=i} \binom{i}{\alpha} \Eb X^\alpha b^\alpha = \Eb \langle X,b\rangle^i.
    $$
    By the above observation and $\mbf_{2k-1}(G)=\mbf_{2k-1}(H)$, we have $\mbf_{2k-1}(\langle X,b\rangle)=\mbf_{2k-1}(\langle Y,b\rangle)$. 
    It then follows from \cite[Lemma 4]{wu2020optimal} that the univariate discrete random variables $\langle X,b\rangle$ and $\langle Y,b\rangle$ have the same distributions. Since $b$ is arbitrary, $X$ and $Y$ have the same distributions by the Cram\'er-Wold device \cite[Section 8.6]{pollard_2001}. \\
(b) Firstly consider the case $q=1$. Write $G=\sum_{i\in [k]}p_i\delta_{\theta_i}$ and $H=\sum_{i\in [k]}\pi_i\delta_{\eta_i}$. Then
$\mbf_{2k-2}(G)=\mbf_{2k-2}(H)$ means 
$$
\sum_{j\in [k]} p_i\theta_i^j = \sum_{j\in [k]} \pi_i\eta_i^j \quad \forall j=0,1,\ldots, 2k-2.
$$
By \cite[Lemma C.4, c)]{wei2022convergence} there are infinite many solutions  $(\pi_1,\ldots,\pi_k,\eta_1,\ldots,\eta_k)$ with $\pi_i>0$ for the above system of equations. That is, there exist infinitely many $H\in \Ec_k(\Rb)$ such that $\mbf_{2k-2}(G)=\mbf_{2k-2}(H)$. 

For any $G=\sum_{i\in [k]}p_i\delta_{\theta_i}\in \Ec_k(\Rb^q)$ with $\theta_i\in \sp(\psi)$, we can write $\theta_i=a_i\psi$ with $a_i\in \Rb$. 
Define $G'=\sum_{i\in [k]}p_i\delta_{a_i}\in \Ec_k(\Rb)$. Then by the last paragraph there exist infinitely many $H'=\sum_{i\in [k]}\pi_i\delta_{b_i}$ such that $\mbf_{2k-2}(G')=\mbf_{2k-2}(H')$. Now consider $H=\sum_{i\in [k]}\pi_i\delta_{\eta_i}\in \Ec_k(\Rb^q)$ with $\eta_i=b_i\psi$. Then for any $\alpha\in \Ic_{2k-2}$, $m_\alpha(G)=\sum_{i\in [k]}p_i a_i^{|\alpha|}\gamma^\alpha =\sum_{i\in [k]}\pi_i b_i^{|\alpha|}\gamma^\alpha=m_\alpha(H)$.
\end{proof}

\begin{lem} \label{lem:technicallemma}
Consider any $k_0\leq k$ and any $G_0 = \sum_{i\in [k_0-1]}p_i^0\delta_{\theta_i^0} + p_{k_0}^0\delta_{\theta_0} \in \Ec_{k_0}(\Theta)$. For any $a>0$, any $b>0$, any sequence $\epsilon_n=o(1)$, and any unit vector $\psi\in \Rb^q$, there exist $G=\sum_{i=1}^{d_1}p_i\delta_{\theta_i}\in \Ec_{d_1}(\Rb^q)$ and $H=\sum_{i=1}^{d_1}\pi_i\delta_{\eta_i}\in \Ec_{d_1}(\Rb^q)$ with
	 $\theta_i,\eta_i\in \sp(\psi)\bigcap \{\theta\in \Theta: \|\theta-\theta_0\|_2< b\}$ for any $i\in [d_1]$ such that:
1) $G_n= \sum_{i=1}^{k_0-1}p_i^0\delta_{\theta_i^0}  + p_{k_0}^0 \sum_{j=1}^{d_1} p_j\delta_{\theta_{0}+\epsilon_n\theta_j}\in \Ec_k(\Theta)$, $H_n= \sum_{i=1}^{k_0-1}p_i^0\delta_{\theta_i^0}  + p_{k_0}^0 \sum_{j=1}^{d_1} \pi_j\delta_{\theta_{0}+\epsilon_n\eta_j}\in \Ec_k(\Theta)$; 
2) $W_1(G_n,G_0)<a\epsilon_n$, $W_1(H_n,G_0)<a\epsilon_n$ and $W_1(G_n,H_n)=\epsilon_n p_{k_0}^0 W_1(G,H)$; 
3) for any function $\myf(\theta)$ that is $(2d_1-1)$-th order continuously differentiable  on  $\{\theta\in \Theta: \|\theta-\theta_0\|_2< b\}$,
	\begin{align*}
		&\left(\frac{\int \myf(\theta)dG_n  - \int \myf(\theta)dH_n}{\epsilon_n^{2d_1-1}}\right)^2 \\
		\leq &  C(d_1,q) \sum_{|\alpha|=2d_1-1} \int_0^1  \sum_{j\in [d_1]} p_j \left( D^{\alpha} \myf (\theta_0+t\epsilon_n\theta_j) \right)^2 + \sum_{j\in [d_1]}\pi_j \left( D^{\alpha} \myf (\theta_0+t\epsilon_n\eta_j) \right)^2 dt. \numberthis \label{eqn:diffsquareuppbou}
	\end{align*}

\end{lem}

\begin{proof}
By Lemma \ref{lem:discretedistributions} \ref{itemb:lem:discretedistributions}, there exist $G=\sum_{i=1}^{d_1}p_i\delta_{\theta_i}\in \Ec_{d_1}(\Rb^q)$ and $H=\sum_{i=1}^{d_1}\pi_i\delta_{\eta_i}\in \Ec_{d_1}(\Rb^q)$ 
	such that $\mbf_{2d_1-2}(G)=\mbf_{2d_1-2}(H)$ and $\theta_i,\eta_i\in \sp(\psi)$ for any $i\in [d_1]$. 
	Denote $\delta_{\bm{0}}$ the Dirac measure at the origin $\bm{0}\in \Rb^q$. We may assume that $W_1(G,\delta_{\bm{0}})<a$, $W_1(H,\delta_{\bm{0}})<a$ and $\|\theta_i\|_2\vee \|\eta_i\|_2\leq b\wedge 1$ for any $i\in [d_1]$; otherwise, simply replace $G$ and $H$ respectively with $S_w G$ and $S_w H$ for small enough $w>0$. 
	Without loss of generality, write $G_0=\sum_{i=1}^{k_0}p_i^0\delta_{\theta_i^0}$ with $\theta_{k_0}^0=\theta_0$. Set $\rho=\frac{1}{2}\min_{1\leq i<j\leq k_0}\|\theta_i^0-\theta_j^0\|_2$. Following the same reasoning as above, we may further require that $\max_{i\in [d_1]}\left\|\theta_i\right\|_2<\rho$ and $\max_{i\in [d_1]}\left\|\eta_i\right\|_2<\rho$.
	
	Consider $G_n= \sum_{i=1}^{k_0-1}p_i^0\delta_{\theta_i^0}  + p_{k_0}^0 \sum_{j=1}^{d_1} p_j\delta_{\theta_{0}+\epsilon_n\theta_j}  = \sum_{i=1}^{k_0-1}p_i^0\delta_{\theta_i^0} + p_{k_0}^0 \left( S_{\epsilon_n}G +\theta_{0}\right) $. Similarly, define $H_n=\sum_{i=1}^{k_0-1}p_i^0\delta_{\theta_i^0} + p_{k_0}^0 \left( S_{\epsilon_n}H +\theta_{0}\right)$. It is clear that $G_n,H_n\in \Ec_{k}(\Theta)$ for any $n\geq 1$ from our construction of $G$ and $H$. Moreover, $G_n, H_n\overset{W_1}{\to} G_0$. Thus we may view that $G_n$ as a sequence on the curve $\{\sum_{i=1}^{k_0-1}p_i^0\delta_{\theta_i^0} + p_{k_0}^0 \left( S_{\epsilon}G +\theta_{0}\right)|\epsilon\in [0,1] \}$ specified 
	by a fixed direction $G$. A similar viewpoint applies to $H_n$. By the definition of Wasserstein distance, $W_1(G_n,G_0)=\epsilon_n p_{k_0}^0 W_1(G,\delta_{\bm{0}})<a\epsilon_n$, $W_1(H_n,G_0)=\epsilon_n p_{k_0}^0 W_1(H,\delta_{\bm{0}})<a\epsilon_n$, and $W_1(G_n,H_n)=\epsilon_n p_{k_0}^0 W_1(G,H)$.
	
	It follows by Taylor's theorem with integral remainder that,
	\begin{align*}
		&\sum_{j\in [d_1]} p_j \myf\left({\theta_{0}+\epsilon_n\theta_j}\right) \\ 
		= & \sum_{0\leq |\alpha|\leq 2d_1-2 } \frac{1}{\alpha!} D^{\alpha} \myf(\theta_0) \epsilon_n^{|\alpha|} m_\alpha(G) + \epsilon_n^{2d_1-1} (2d_1-1) \sum_{|\alpha|=2d_1-1}   \int_0^1 (1-t)^{2d_1-2} \sum_{j\in [d_1]} p_j\psi_{n,\alpha}(t|\theta_j)  dt,
	\end{align*}
	where $\psi_{n,\alpha}(t|\theta)=\frac{\theta^{\alpha}}{\alpha!} D^{\alpha} \myf (\theta_0+t\epsilon_n\theta)$. A similar formula holds for $\sum_{j\in [d_1]} \pi_j \myf\left({\theta_{0}+\epsilon_n\eta_j}\right) $. Thus 
	\begin{align*}
	&	\int \myf(\theta)dG_n  - \int \myf(\theta)dH_n\\
		= &p_{k_0}^0\left(\sum_{j\in [d_1]} p_j \myf\left({\theta_{0}+\epsilon_n\theta_j}\right) -\sum_{j\in [d_1]} \pi_j \myf\left({\theta_{0}+\epsilon_n\eta_j}\right)\right)\\
		= & p_{k_0}^0 \epsilon_n^{2d_1-1} (2d_1-1) \sum_{|\alpha|=2d_1-1}   \int_0^1 (1-t)^{2d_1-2} \left(\int \psi_{n,\alpha}(t|\theta)d\left(G-H\right)\right)  dt  .
	\end{align*}
	Then 
	\begin{align*}
		&\left(\frac{\int \myf(\theta)dG_n  - \int \myf(\theta)dH_n}{\epsilon_n^{2d_1-1}}\right)^2 \\
		\leq & C(d_1,q) \sum_{|\alpha|=2d_1-1} \left(\int_0^1 (1-t)^{2d_1-2} \left(\int \psi_{n,\alpha}(x,t|\theta)d\left(G-H\right)\right)  dt\right)^2 \\
		\overset{(*)}{\leq} &  C(d_1,q) \sum_{|\alpha|=2d_1-1}  \int_0^1  \left(\int \psi_{n,\alpha}(x,t|\theta)d\left(G-H\right)\right)^2 dt \\
		\overset{(**)}{\leq} &  C(d_1,q) \left(b\wedge 1 \wedge \rho\right)^{2d_1-1} \sum_{|\alpha|=2d_1-1} \int_0^1  \sum_{j\in [d_1]} p_j \left( D^{\alpha} \myf (\theta_0+t\epsilon_n\theta_j) \right)^2 + \sum_{j\in [d_1]}\pi_j \left( D^{\alpha} \myf (\theta_0+t\epsilon_n\eta_j) \right)^2 dt 
		 \\
		\leq &  C(d_1,q)  \sum_{|\alpha|=2d_1-1} \int_0^1  \sum_{j\in [d_1]} p_j \left( D^{\alpha} \myf (\theta_0+t\epsilon_n\theta_j) \right)^2 + \sum_{j\in [d_1]}\pi_j \left( D^{\alpha} \myf (\theta_0+t\epsilon_n\eta_j) \right)^2 dt,
	\end{align*}
	where step $(*)$ follows by Cauchy-Schwartz formula for the integral, and step $(**)$ follows Cauchy-Schwartz formula for the integrand, and the fact that $\|\theta_j\|_2 \vee \|\eta_j\| \leq b\wedge 1 \wedge \rho$ for any $j\in [d_1]$.
\end{proof}

\begin{proof}[Proof of Theorem \ref{thm:minimax}]

	By the two-point Le Cam bound (see (15.14) in \cite{wainwright2019high}\footnote{Strictly speaking, their setting with parameter as a functional of the probability measure does not directly applies to our setting. If we assume the map $G'\to P_{G'}$ is injective on $\Gc_k(\Theta)$, or equivalently the mixture model is identifiable, then it is safe to view $G'$ as a functional of $P_{G'}$ and hence the cited result directly applies. But a proof following the proof of the cited result line by line produces the same conclusion in our setting (that probability measure is a functional of the parameter $G'$), without requiring the identifiability assumption.}),  
	for any $G_n,H_n\in \Gc_k(\Theta)$ satisfying $W_1(G_n,G_0)<a \epsilon_n$ and $W_1(H_n,G_0)<a \epsilon_n$,
	\begin{align}
		\inf_{\hat{G}_n\in \Ef_n}\ \sup_{ \substack{G^*: W_1(G^*,G_0)<a \epsilon_n }  } \Eb_{G^*} W_1(\hat{G}_n,G^*) \geq & \frac{W_1(G_n,H_n)}{4} \left(1- V\left(\bigotimes{^n}\Pb_{G_n},\bigotimes{^n} \Pb_{H_n} \right) \right), \label{eqn:lecamlowerbou}
	\end{align}
where $\bigotimes{^n}\Pb_{G_n}$ denotes the product measure on the product space $(\Xf^n,\Xc^n)$.

It then suffices to choose $G_n$ and $H_n$ such that the right hand side of \eqref{eqn:lecamlowerbou} is large. Let $b>0$ be the constant and $\psi\in \Rb^q$ be the unit vector in the definition of the assumption $A(\theta_0,d_1)$. For $\epsilon_n=n^{-\frac{1}{2d_1-1}}$, let $G,H, G_n,H_n$ be specified in Lemma \ref{lem:technicallemma}. Then by the property 2) in Lemma \ref{lem:technicallemma},

\begin{align*}
	\inf_{\hat{G}_n\in \Ef_n}\ \sup_{ \substack{G^*: W_1(G^*,G_0)<a \epsilon_n }  } \Eb_{G^*} W_1(\hat{G}_n,G^*)	\geq & \frac{\epsilon_n p_{k_0}^0 W_1(G,H)}{4} \left(1- h\left(\bigotimes{^n}\Pb_{G_n},\bigotimes{^n} \Pb_{H_n} \right)  \right). \numberthis \label{eqn:minimaxlowboutemp1}
\end{align*}

	It suffices now to bound $h\left(\bigotimes{^n}\Pb_{G_n},\bigotimes{^n} \Pb_{H_n} \right) $. Note that 
	\begin{align*}
		& n h^2\left(\Pb_{G_n}, \Pb_{H_n} \right)\\
		= &\frac{ h^2\left(\Pb_{G_n}, \Pb_{H_n} \right)}{\left(\epsilon_n^{2d_1-1}\right)^2}\\
		=  & 
		\frac{ 1}{2\left(\epsilon_n^{2d_1-1}\right)^2} 
		\int \frac{\left(\int p(x \mid \theta)dG_n  -\int p(x \mid \theta)dH_n\right)^2}{\left(\sqrt{\int p(x \mid \theta)dG_n}  + \sqrt{\int p(x \mid \theta)dH_n}\right)^2}  d\lambda  \\
		\overset{(*)}{\leq}  & 
	C(d_1,q) \sum_{|\alpha|=2d_1-1} \int \frac{\int_0^1  \sum_{i\in [d_1]} p_i \left( D^{\alpha} p (x|\theta_0+t\epsilon_n\theta_i) \right)^2 + \sum_{i\in [d_1]}\pi_i \left( D^{\alpha} p (x|\theta_0+t\epsilon_n\eta_i) \right)^2 dt}{\int p(x \mid \theta)dG_n  + \int p(x \mid \theta)dH_n} d\lambda\\
		{\leq}  &  
		C(d_1,q) \sum_{|\alpha|=2d_1-1} \int \int_0^1 \sum_{i\in [d_1]} \frac{ \left( D^{\alpha} p (x|\theta_0+t\epsilon_n\theta_i) \right)^2  }{ p(x|\theta_0+\epsilon_n\theta_i) } +
		\sum_{i\in [d_1]} \frac{   \left( D^{\alpha} p (x|\theta_0+t\epsilon_n\eta_i) \right)^2  }{ p(x|\theta_0+\epsilon_n\eta_i) } dt d\lambda \\
		\overset{(**)}{=}  &  
		C(d_1,q) \sum_{|\alpha|=2d_1-1} \int_0^1 \int  \sum_{i\in [d_1]} \frac{ \left( D^{\alpha} p (x|\theta_0+t\epsilon_n\theta_i) \right)^2  }{ p(x|\theta_0+\epsilon_n\theta_i) } +
		\sum_{i\in [d_1]} \frac{   \left( D^{\alpha} p (x|\theta_0+t\epsilon_n\eta_i) \right)^2  }{ p(x|\theta_0+\epsilon_n\eta_i) } d\lambda dt  \\
		\overset{(***)}{\leq} &C_0(d_1,q,A), \numberthis \label{eqn:hellingeruppbound}
	\end{align*}
	where step $(*)$ follows from \eqref{eqn:diffsquareuppbou} with $\myf(\theta)=p(x \mid \theta)$, step $(**)$ follows from Tonelli Theorem and the joint Lebesgue measurability of the integrand (due to \cite[Lemma 4.51]{guide2006infinite}),
	and step $(***)$ follows from \eqref{eqn:minimaxhypothesis2} since $\theta_i,\eta_i\in \sp(\psi)$, with $C_0(d_1,q,b,G_0,A)$ a positive constant. By \eqref{eqn:hellingeruppbound}, when $n> C_0(d_1,q,A)$,
	\begin{align*}
		1-h^2\left(\bigotimes{^n}\Pb_{G_n},\bigotimes{^n} \Pb_{H_n} \right) = \left(1- h^2\left(\Pb_{G_n},\Pb_{H_n}\right)\right)^n \geq \left(1- \frac{C_0(d_1,q,A)}{n}\right)^n \geq c_0(d_1,q,A),
	\end{align*}
	where the last step follows from $\left(1- \frac{C_0(d_1,q,A)}{n}\right)^n\to e^{-C_0(d_1,q,A)}>0$ and $c_0(d_1,q,A)$ is a positive constant. 
	The above inequality immediately implies that when $n> C_0(d_1,q,A)$
	$$
	h\left(\bigotimes{^n}\Pb_{G_n},\bigotimes{^n} \Pb_{H_n} \right)\leq \sqrt{1-c_0(d_1,q,A)}<1.
	$$

	Plugging the preceding inequality into \eqref{eqn:minimaxlowboutemp1} yields \eqref{eqn:uniformintegraluppbou} for $n> C_0(d_1,q,A)$. \eqref{eqn:uniformintegraluppbou} for $n\leq C_0(d_1,q,A)$ can be obtained directly by tuning the constant coefficient in its lower bound.
\end{proof}

\subsection{Proofs of Lemma \ref{lem:localtoglobal} and Lemma \ref{lem:convergencerate}}

\begin{proof}[Proof of Lemma \ref{lem:localtoglobal}]
	Suppose that \eqref{eqn:globalinversebound} does not hold. Then there exists $G_n\neq H_n\in \Gc_k(\Theta)$, such that 
	\begin{equation}
	\frac{\sup_{\myf\in \myF} |G_n \myf-H_n \myf|}{W_{2k-1}^{2k-1}(G_n,H_n)} \to 0. \label{eqn:localtoglobal1}
	\end{equation}
	Since $\Theta$ is compact, $\Gc_k(\Theta)$ is compact. Then by considering subsequence if necessary, we may require $G_n\overset{W_1}{\to} G_\infty\in \Gc_k(\Theta)$ and $H_n \overset{W_1}{\to} H_\infty\in \Gc_k(\Theta)$. If $G_\infty = H_\infty$, then \eqref{eqn:localtoglobal1} contradicts with \eqref{eqn:localinversebound} for $G_0=G_\infty$ since $W_{2k-1}^{2k-1}(G_n,H_n)\leq \left(\diam(\Theta) \right)^{2(k-d_1)} W_{2d_1-1}^{2d_1-1}(G_n,H_n)$. Thus we have $G_\infty\neq H_\infty$, but then \eqref{eqn:localtoglobal1} implies that $\sup_{\myf\in \myF} |G_\infty \myf-H_\infty \myf|=0$, which contradicts with the assumption that $\Gc_k(\Theta)$ is distinguishable by $\myF$.
\end{proof}

\myred{
\begin{proof}[Proof of Lemma \ref{lem:equinvbou}]
\noindent
Parts \ref{lem:equinvboua} and \ref{lem:equinvboub} are trivial. \\

\noindent
	\ref{lem:equinvbouc}. 
By \eqref{eqn:localinversebound}, there exists $r>0$ such that for any $G_1,H\in B_{W_1}(G_0,r)$, the $W_1$-ball centering at $G_0$ of radius $r$ in $\Gc_k(\Theta)$, we have
\begin{align}
\sup_{\myf\in \myF} |G_1 \myf-H \myf|\geq C(G_0,\myF,\Theta,k_0,k) W_{2d_1-1}^{2d_1-1}(G_1,H).  \label{eqn:invboudircon}
\end{align}

Define 
$$
z := \inf_{ \substack{G_1\in \bar{B}_{W_1}(G_0,r/2)\\ H \in \Gc_k(\Theta) \setminus {B_{W_1}(G_0,r)} } }\ \  \sup_{\myf\in \myF} |G_1 \myf-H \myf|,
$$
where $\bar{B}_{W_1}(G_0,r/2)$ is the closed ball. 
Since $\sup_{\myf\in \myF} |G_1 \myf-H \myf|$ is lower semicontinuous 
on the compact set $\bar{B}_{W_1}(G_0,r/2) \times \left(\Gc_k(\Theta) \setminus B_{W_1}(G_0,r)
\right)$, the infimum is attained. Since $G_k(\Theta)$ is distinguishable by $\myF$, $z>0$.  Since $W_{2d_1-1}^{2d_1-1}(G_1,H)\leq \diam^{2d_1-1}(\Theta)$, we have $G_1\in \bar{B}_{W_1}(G_0,r/2)$ and $H \in \Gc_k(\Theta) \setminus {B_{W_1}(G_0,r)}$:
\begin{align}
\sup_{\myf\in \myF} |G_1 \myf-H \myf|\geq \frac{z}{\diam^{2d_1-1}(\Theta)} W_{2d_1-1}^{2d_1-1}(G_1,H).  \label{eqn:comdisdircon}
\end{align}

Combining \eqref{eqn:invboudircon} and \eqref{eqn:comdisdircon} completes the proof. The other direction follows since \ref{lem:equinvbouc} implies \ref{lem:equinvboub}. 

\end{proof}
}

\subsection{Proof of Theorem \ref{thm:inversebound} \ref{thm:inversebounda}}
\label{sec:proofoflocalinversebound}

\textbf{Notation for this subsection.} When comparing sequences, we will write $a_n\preccurlyeq b_n$ or $a_n=O(b_n)$ for $a_n\leq Cb_n$ where $C>0$ does not depend on $n$ but may depend on other parameters. We also write $a_n\succcurlyeq b_n$ if $ b_n \preccurlyeq a_n$. We will furthermore use $a_n \asymp b_n$ if $b_n\preccurlyeq a_n \preccurlyeq b_n$.

\begin{proof}[Proof of Theorem \ref{thm:inversebound} \ref{thm:inversebounda}] 

The proof is divided into the following steps.

\noindent
\textbf{Step 1}: (Proof by contradiction and subsequences) 
Suppose that \eqref{eqn:localinversebound} does not hold. Then there exists  $G_n\neq H_n\in \Gc_k(\Theta)$ and $G_n,H_n \overset{W_1}{\to} G_0$ 
such that 
\begin{equation}
\lim_{ n\to\infty} \frac{\sup_{\myf\in \myF} |G_n \myf-H_n \myf|}{W_{2d_1-1}^{2d_1-1}(G_n,H_n)} = 0. \label{eqn:limitzero}
\end{equation}

Since $\Theta$ is compact, by taking subsequence if necessary, we have that for each $n$: 1) $G_n\in \Ec_{m_1}(\Theta)$ and $H_n\in \Ec_{m'}(\Theta)$ with $m_1,m'\in [k_0,k]$ independent of $n$; 2) $G_n=\sum_{j\in [m_1]}p_{jn}\delta_{\theta_{jn}}$ and $H_n=\sum_{j \in [m']}\pi_{jn}\delta_{\eta_{jn}}$ with
\begin{align*}
    \sum_{j\in [m_1]} p_{jn} = 1, &  \sum_{j\in [m']} \pi_{jn} = 1, \\
    p_{jn}>0,\ \theta_{jn} \text{ all distinct}, & \quad  \pi_{jn}>0,\ \eta_{jn} \text{ all distinct},   \\
    \theta_{jn} \to \theta_j, & \quad \eta_{jn} \to \eta_j.
    \numberthis \label{eqn:supporpointslimit}
\end{align*}
For each $n$, set
$$
(\omega_{jn}, \nu_{jn}) = 
\begin{cases}  
(p_{jn},\theta_{jn}), & \text{ if } j\leq m_1, \\
(-\pi_{(j-m_1)n},\eta_{(j-m_1)n}), & \text{ if } m_1<j\leq m_1+m'.
\end{cases}
$$

\noindent
\textbf{Step 2}: (Decreasing rate of Wasserstein distance) 
Each member in the sequence of sets $(\{\nu_{jn} | j \in [m_1+m']\})_{n=1}^{\infty}$ defined in the previous step contains the supporting atoms from the pair of measures $G_n$ and $H_n$, which tend to $G_0$ under the Wasserstein distance $W_{2d_1-1}$ (and also $W_1$). 
These sets of atoms can be partitioned into groups using a useful tree structure introduced by \cite{heinrich2018strong}. 
This step of the proof proceeds by adapting from  
\cite[Lemma 7.1, Definition 7.2, Lemma 7.3]{heinrich2018strong} and hence the proofs are omitted here. 
\myred{We also $|\nu_{in}-\nu_{jn}|$ to represent some norm between $\nu_{in}$ and $\nu_{jn}$ on $\Rb^q$ (to be concrete, one can think that it is the $\|\cdot\|_\infty$). }
First, it is simple to note the following:
\begin{lem}[Discrepancy orders of $\nu_{jn}$] \label{lem:rateofatoms}
By taking a subsequence of $\{\nu_{jn}, j \in [m_1+m']\}_{n=1}^{\infty}$ if necessary, there exists a finite number $S \leq m_1 m'$ of ``scaling'' sequences
$$
0:=\epsilon_0(n)<\epsilon_1(n)<\cdots<\epsilon_S(n):=1 \quad \text{ with } \epsilon_s(n)=o(\epsilon_{s+1}(n)),
$$
such that, for any $i, j\in [m_1+m']$ there is a unique $s(i,j)\in [S]\cup\{0\}$ satisfying $|\nu_{in}-\nu_{jn}|\asymp \epsilon_{s(i,j)}(n)$. 
\end{lem}

Note that $S$ and $s(\cdot,\cdot)$ are independent of $n$. Moreover, $s(\cdot,\cdot)$ is a ultrametric on $[m_1+m']$, i.e., $s(\cdot,\cdot)$ satisfies all the requirements of a distance except the triangle inequality, which is replaced by $s(i,j)\leq \max\{s(i,\ell),s(\ell,j)\}$. It follows immediately that on the space $([m_1+m'],s(\cdot,\cdot))$, the closed balls $\bar{B}_s(i,r)$ w.r.t. the ultrametric $s(\cdot,\cdot)$, with center $i\in [m_1+m']$ and radius $r\in [S]\cup \{0\}$, are either disjoint or in the case that one is a subset of the other. This leads to the following definition.

\begin{definition}[coarse-grained tree] The vertices of the coarse-grained tree $\Tc$ are the balls $\{\bar{B}_s(i,r)|i\in [m_1+m'], r\in [S]\cup \{0\} \}$. The root of $\Tc$ is $J_r=[m_1+m']$. For each vertex $J\neq J_r$, its parent $J^\uparrow$ is the vertex that (as a set) contains $J$ as a subset and has the smallest cardinality.
\end{definition}

For a given vertex $J$, the set of its children, descendents are respectively denoted by $\Child(J)$, $\Desc(J)$. The diameter of a vertex $J$ is $s(J):=\max_{i,j\in J} s(i,j)$, which is also the radius since $s(\cdot,\cdot)$ is a ultrametric; in fact, $J=\bar{B}_s(i,s(J))$ for any $i\in J$. Note that $\Tc$ is constructed based on the sequence $\{\nu_{jn}\}_{j\in [m_1+m'], n\geq 1}$ but does not depend on $n$. 

One essential property of $\Tc$ is that for any $i\in K, j \in K'$ where $K\neq K'\in \Child(J)$, $s(i,j)=s(J)$ since $s(K)<s(J)$ and $s(K')<s(J)$. Translating the previous sentence in terms of the $\nu_{jn}$, it means $|\nu_{in}-\nu_{jn}|\asymp \epsilon_{s(J)}(n)$. Thus the coarse-grained tree $\Tc$ is a device to keep track of the partitioning of the supporting atoms $\{\nu_{jn}\}$ into groups in terms of the decreasing rate of their pairwise distances. The following simple facts about this device are useful. 

\begin{lem}
\label{lem:treecount}
\begin{enumerate}[label=(\alph*)]
\item
 $|\Child(J_r)| \in [k_0, m_1+m'-k_0]\subset [k_0,2k-k_0]$. Moreover, $\sum_{J\in \Child(J_r)}|J| = |J_r| = m_1+m'\in [2k_0,2k]$.
 \item  \label{lem:treecountb}
 If $\epsilon_J(n)=o(1)$, then $|J|\leq m_1+m'-2(k_0-1) \leq 2d_1$.
 \end{enumerate}
\end{lem}
\begin{proof} (a) Trivial.
(b) If $\epsilon_{J_r}(n)=o(1)$, then $k_0=1$ and thus the statement holds. 
If $\epsilon_{J_r}(n)=1$, then it suffices to prove that $|J|\leq m_1+m'-2(k_0-1)$ for any $J\in \Child(J_r)$. Since $G_n,H_n\overset{W_1}{\to}G_0$, it then holds that there are at least $k_0$ children of $J_r$ having cardinality at least $2$. So
$$m_1+m'=|J_r| = \sum_{J\in \Child(J_r)}|J| \geq  \max_{J\in \Child(J_r)}|J|  +  2(k_0-1). $$
\end{proof}

Set for short
$$
\bar{\omega}_{Jn} := \sum_{j\in J}\omega_{jn}, \text{ and }  \quad \epsilon_{J}(n) := \epsilon_{s(J)}(n).
$$
\begin{lem}[Characterization of the decreasing rate of $W_\ell^\ell(G_n,H_n)$] 
\label{lem:wassersteincharacterization}
For any $\ell\geq 1$, we have 
$$
W_\ell^\ell(G_n,H_n) \asymp \max_{J\in \Desc(J_r)} |\bar{\omega}_{Jn}|\left(\epsilon_{J^\uparrow}(n)\right)^\ell.
$$
\end{lem}

\noindent
\textbf{Step 3}: (Expansion of the integral of $\myf$ w.r.t. a signed measure)


Consider the signed measure $GH_{Jn}:=\sum_{j\in J}\omega_{jn}\nu_{jn}$ and let $GH_{Jn} \myf := \int \myf d GH_{Jn}= \sum_{j\in J} \omega_{jn} \myf(\nu_{jn})$.  

\begin{lem} 
\label{lem:numeratortaylormul}
 For each vertex $J$ of $\Tc$, choose an index $i(J)\in J$ (that is independent of $n$) and denote $\nu_{Jn}=\nu_{i(J) n}$.
Consider any vertex $J$ of $\Tc$ with $\epsilon_{J}(n)=o(1)$, any $m\in [|J|-1,2d_1-1]$. 
Then for any $m$-th order continuously differentiable function $\myf$ defined on $\Theta$,
\begin{equation}
GH_{Jn} \myf =\sum_{\alpha\in\Ic_m} a(\alpha|J,\nu_{Jn}) \left(\epsilon_J(n)\right)^{|\alpha|} D^\alpha \myf(\nu_{Jn}) + R(\myf,J,\nu_{Jn}), \label{eqn:numeratortylortheoremmul}
\end{equation}
  with $a(\alpha|J,\nu_{Jn})=\frac{m_\alpha(GH_{Jn}-\nu_{Jn})}{\alpha!(\epsilon_J(n))^{|\alpha|}}$ (in particular, $a(\bm{0}|J,\nu_{Jn})=\bar{\omega}_{Jn}$) satisfying the following: 
\begin{enumerate}[label=(\alph*)]
\item \label{itema:lem:numeratortaylormul}
If $J$ is a leaf vertex, then  $a(\alpha|J,\nu_{Jn})=0$ for $1\leq |\alpha| \leq m$. If $J$ is not a leaf vertex, then
$$ 
\max_{|J|\leq |\alpha|\leq m } |a(\alpha|J,\nu_{Jn})|   \preccurlyeq \max_{K\in \Child(J)}M_{m,K} (\nu_{Kn}) \asymp \max_{K\in \Child(J)}M_{|K|-1,K} (\nu_{Kn}) \asymp \max_{0\leq |\alpha|\leq |J|-1 } |a(\alpha|J,\nu_{Jn})|,
$$ 
where $M_{p,K}(\nu_{Kn}):=\max_{0\leq |\gamma|\leq p}\left|a(\gamma|K,\nu_{Kn}) \left(\frac{\epsilon_K(n)}{\epsilon_{K^\uparrow}(n)}\right)^{|\gamma|} \right|$. \\
\item \label{itemb:lem:numeratortaylormul}
Denote $a(J,\nu_{Jn}):=(a(\alpha|J,\nu_{Jn}))_{\alpha\in \Ic_m}$. If $J$ is not a leaf vertex, then 
$$
\| a(J,\nu_{Jn}) \|_\infty \succcurlyeq \max_{K\in \Desc(J)}\left| \bar{\omega}_{Kn} \left(\frac{\epsilon_{K^\uparrow}(n)}{\epsilon_J(n)}\right)^{|J|-1} \right|.
$$
\item \label{itemc:lem:numeratortaylormul}
If $J$ is a leaf vertex, then  $ R(\myf,J,\nu_{Jn})=0$. If $J$ is a not a leaf vertex, then $ R(\myf,J,\nu_{Jn})=o\left( \left\| a\left(J,\nu_{Jn}\right) \right\|_\infty \left(\epsilon_J(n)\right)^{m} \right)$ and thus $R(\myf,J,\nu_{Jn})= o\left( \max_{\alpha\in \Ic_m} \left|a(\alpha|J,\nu_{Jn})\left(\epsilon_{J}(n)\right)^{|\alpha|}\right| \right)$.
\myred{\item \label{itemd:lem:numeratortaylormul}
Suppose in addition, that there is a uniform continuity modulus $w(\cdot)$ such that: for any $\alpha$ with $|\alpha|=m$,  
$$\sup_{\myf\in \myF}|D^\alpha \myf(\theta) - D^\alpha \myf(\theta')|\leq w(\theta-\theta')$$ 
with $\lim_{h\to 0}  w(h) =0$.
Then  $ \sup_{\myf\in \myF} |R(\myf,J,\nu_{Jn})|=o\left( \left\| a\left(J,\nu_{Jn}\right) \right\|_\infty \left(\epsilon_J(n)\right)^{m} \right)$.

Thus $\sup_{\myf\in \myF} |R(\myf,J,\nu_{Jn})|= o\left( \max_{\alpha\in \Ic_m} \left|a(\alpha|J,\nu_{Jn})\left(\epsilon_{J}(n)\right)^{|\alpha|}\right| \right)$.
}
\end{enumerate}
\end{lem}

Lemma \ref{lem:numeratortaylormul} is a multivariate version of the univariate result \cite[Lemma 7.4]{heinrich2018strong}. Moreover, Lemma \ref{lem:numeratortaylormul} is an improvement of \cite[Lemma 7.4]{heinrich2018strong}, as the former requires less differentiability assumption and no assumption on uniform continuity on the derivative in comparison to the latter (see  \cite[third bullet point in Assumption B(k) on Page 2850]{heinrich2018strong})\myred{; essentially that is equivalent to the additional assumption in part \ref{itemd:lem:numeratortaylormul}, but part \ref{itemd:lem:numeratortaylormul} is not needed in the proof of Theorem \ref{thm:inversebound}}.   An important observation 
is that $a(\alpha|J,\nu_{Jn})$ does not depends on $\myf$. The proof of Lemma \ref{lem:numeratortaylormul} is deferred to Section \ref{sec:proofoflem:numeratortaylor}.\\

\noindent 
\textbf{Step 4}: (Deriving contradiction with that $\myF$ is a $(2d_1-1,k_0,k)$ linear independence domain). There are two cases: either
$\epsilon_{J_r}(n)=1$ or $\epsilon_{J_r}(n) = o(1)$. \\

\noindent
\textit{Case 1}: Suppose $\epsilon_{J_r}(n)=1$. Notice that by Lemma \ref{lem:wassersteincharacterization},
\begin{align*}
    W_{2d_1-1}^{2d_1-1}(G_n,H_n) 
    \asymp & \max\left\{ \max_{J\in \Child(J_r)} |\bar{\omega}_{Jn}|, \max_{\substack{J\in \Child(J_r)\\ J \text{ non-leaf} }} \max_{K\in \Desc(J)} |\bar{\omega}_{Kn}|\left(\epsilon_{K^\uparrow}(n)\right)^{2d_1-1}\right\}\\
    \preccurlyeq & 
    \max\left\{ \max_{J\in \Child(J_r)} |\bar{\omega}_{Jn}|, \max_{\substack{J\in \Child(J_r)\\ J \text{ non-leaf} }} \left(\epsilon_{J}(n)\right)^{|J|-1} \|a(J,\nu_{Jn})\|_{\infty} \right\} \\
    \leq &   
    \underbrace{\max_{J\in \Child(J_r)}\ \max_{|\alpha|\leq |J|-1} \left|a(\alpha|J,\nu_{Jn}) \left(\epsilon_J(n)\right)^{|\alpha|}\right|}_{:=d_n},
\end{align*}
where the ``$\preccurlyeq$'' step follows from Lemma \ref{lem:numeratortaylormul} \ref{itemb:lem:numeratortaylormul} and  $|J|\leq 2d_1$ due to Lemma \ref{lem:treecount} \ref{lem:treecountb}, and in the last step $a(\bm{0}|J,\nu_{Jn})=\bar{\omega}_{Jn}$ is used.  

Since $G_n\neq H_n$ and $\epsilon_{J_r}(n)=1$, $\Child(J_r)$ is not empty. 
Since $\epsilon_J(n)=o(1)$ for any $J\in \Child(J_r)$, by Lemma \ref{lem:numeratortaylormul} with $m=|J|-1$ for each $J$, 
\begin{align*}
    \frac{|G_n \myf-H_n \myf|}{W_{2d_1-1}^{2d_1-1}(G_n,H_n)} \succcurlyeq & \frac{|G_n \myf-H_n \myf|}{d_n}\\
    = &
     \left|\sum_{J\in \Child(J_r)} \left(\sum_{|\alpha|\leq {|J|-1}}  D^\alpha \myf(\nu_{Jn}) \frac{a(\alpha|J,\nu_{Jn}) \left(\epsilon_J(n)\right)^{|\alpha|}}{d_n} + \frac{R(\myf,J,\nu_{Jn})}{d_n}\right)\right|. \numberthis \label{eqn:childexpansion}
\end{align*}
 It follows from Lemma \ref{lem:numeratortaylormul} \ref{itemc:lem:numeratortaylormul} and the condition $m= |J|-1$ that
\begin{align*}
\frac{R(\myf,J,\nu_{Jn})}{d_n}    = o(1). \numberthis \label{eqn:remainderlimit0}
\end{align*}
By taking subsequence if necessary, we have that 
\begin{equation}
\frac{a(\alpha|J,\nu_{Jn}) \left(\epsilon_J(n)\right)^{|\alpha|}}{d_n} \to b_{J\alpha} \label{eqn:coefficientlimit}
\end{equation} 
for some $b_{J\alpha}\in [-1,1]$. Moreover, at least one of $\{b_{J\alpha}\}$ has magnitude $1$. We also have 
\begin{equation}
    \sum_{J\in \Child(J_r)} b_{J\bm{0}} = 0 \label{eqn:sum0}
\end{equation}
since $\sum_{J\in \Child(J_r)} a(\bm{0}|J,\nu_{Jn})=\sum_{J\in \Child(J_r)}\bar{\omega}_{Jn}=\sum_{j\in [m_1]}p_{jn}-\sum_{j\in [m']}\pi_{jn}=0$. 

Then following \eqref{eqn:limitzero},
\begin{align*}
0= &\lim_{ n\to\infty} \frac{\sup_{\myf\in \myF} |G_n \myf-H_n \myf|}{W_{2d_1-1}^{2d_1-1}(G_n,H_n)}\\
\geq & \sup_{\myf\in \myF} \liminf_{ n\to\infty} \frac{ |G_n \myf-H_n \myf|}{W_{2d_1-1}^{2d_1-1}(G_n,H_n)}\\
\succcurlyeq & \sup_{\myf\in \myF} \left| \sum_{J\in \Child(J_r)} \ \sum_{|\alpha|\leq |J|-1} b_{J\alpha} D^\alpha \myf(\nu_J)  \right|, \numberthis \label{eqn:differentiatedomain}
\end{align*}
where the last step follows from \eqref{eqn:childexpansion}, \eqref{eqn:remainderlimit0}, \eqref{eqn:coefficientlimit} and that $\nu_{Jn}\to \nu_J$, due to our choice of $\nu_{Jn}$ in Lemma \ref{lem:numeratortaylormul}, and the limit $\nu_J$ exists due to \eqref{eqn:supporpointslimit}.

Since $\epsilon_{J_r}(n)\asymp 1$, $\nu_J$ for different $J\in \Child(J_r)$ are all distinct. Moreover, by Lemma \ref{lem:treecount},
 $|\Child(J_r)| \in  [k_0,2k-k_0]$ and $\sum_{J\in \Child(J_r)}|J| \in [2k_0,2k]$.
That the equations \eqref{eqn:sum0} and \eqref{eqn:differentiatedomain} hold with at least one $b_{J\alpha}$ nonzero contradicts with the hypothesis that $\myF$ is a $(2d_1-1,k_0,k)$ linear independence domain. \\

\noindent
\textit{Case 2}: $\epsilon_{J_r}(n)=o(1)$. 
This implies that $G_0=\delta_{\theta}$ for some $\theta\in \Theta$ and $\nu_{jn}\to \theta$ for any $j\in J_r$. Notice that by Lemma \ref{lem:wassersteincharacterization},
\begin{align*}
    W_{2d_1-1}^{2d_1-1}(G_n,H_n) 
    \asymp &   \max_{K\in \Desc(J_r)} |\bar{\omega}_{Kn}|\left(\epsilon_{K^\uparrow}(n)\right)^{2d_1-1}\\
    \preccurlyeq &   \left(\epsilon_{J_r}(n)\right)^{2d_1-1} \|a(J_r,\nu_{J_rn})\|_{\infty}  \\
    \leq &  \max_{|\alpha|\leq 2d_1-1} \left|a(\alpha|J_r,\nu_{J_r n}) \left(\epsilon_{J_r}(n)\right)^{|\alpha|}\right|
\end{align*}
where the ``$\preccurlyeq$'' step follows from Lemma \ref{lem:numeratortaylormul} \ref{itemb:lem:numeratortaylormul}.

By Lemma \ref{lem:numeratortaylormul}
\begin{align*}
    \frac{G_n \myf-H_n \myf}{W_{2d_1-1}^{2d_1-1}(G_n,H_n)} = \sum_{|\alpha|\leq 2d_1-1}  D^\alpha \myf(\nu_{J_r n}) \frac{a(\alpha|J_r,\nu_{J_rn}) \left(\epsilon_{J_r}(n)\right)^{|\alpha|}}{W_{2d_1-1}^{2d_1-1}(G_n,H_n)} + \frac{R(\myf,J_r,\nu_{J_rn})}{W_{2d_1-1}^{2d_1-1}(G_n,H_n)}.
\end{align*}
 The remainder of the proof for this case involves deriving a contradiction, which is done in the same manner as that of Case 1 above.

\end{proof}

\subsection{Proof of auxiliary lemmas in Section \ref{sec:proofoflocalinversebound}}
\label{sec:proofoflem:numeratortaylor}

\begin{proof}[Proof of Lemma \ref{lem:numeratortaylormul}]
If $J$ is a leaf vertex of $\Tc$, then $\epsilon_{J}(n)=0$ for all $n$, that is, $\nu_{jn}$ for any $j\in J$ are all the same, which we denote $\nu_{Jn}$. (In fact, $J$ has cardinality either $1$ or $2$, where the second case corresponds to $\theta_{in}=\eta_{jn}$ for all $n$, for some $i,j\in J$.) Thus $GH_{Jn} \myf=\bar{\omega}_{Jn} \myf(\nu_{Jn})$, i.e., $a(\alpha|J,\nu_{Jn})=0$ for $0<|\alpha| \leq m$, and $R(\myf,J,\nu_{Jn})=0$. 

Now suppose that the statements \ref{itema:lem:numeratortaylormul}, \ref{itemb:lem:numeratortaylormul}, \ref{itemc:lem:numeratortaylormul}   hold for any $K\in \Desc(J)$ where $J$ is a not a leaf vertex and $\epsilon_J(n)=o(1)$. It suffices to prove \ref{itema:lem:numeratortaylormul}, \ref{itemb:lem:numeratortaylormul}, \ref{itemc:lem:numeratortaylormul} hold also for $J$. (If it is proved, then by mathematical induction, the proof is completed.) 

By assumption  $\epsilon_J(n)=o(1)$, we have that $\nu_{jn}\to\nu_0$ for any $j\in J$. Consider any $K\in \Child(J)$. For any multi-index $\gamma$ such that $|\gamma| \leq m$, applying Taylor's theorem to the function $D^\gamma \myf(\nu_{Kn})$ we have
\begin{align*}
 D^\gamma \myf(\nu_{Kn}) = & \sum_{\alpha\geq \gamma, |\alpha|\leq m-1}  \frac{1}{(\alpha-\gamma)!}(\nu_{Kn}-\nu_{Jn})^{\alpha-\gamma} D^\alpha \myf(\nu_{Jn}) + \sum_{\alpha\geq \gamma, |\alpha|= m}  \frac{r_\alpha(\nu_{Jn},\nu_{Kn})}{(\alpha-\gamma)!}(\nu_{Kn}-\nu_{Jn})^{\alpha-\gamma} \\
=& \sum_{\alpha\geq \gamma, |\alpha|\leq m}  \frac{1}{(\alpha-\gamma)!}(\nu_{Kn}-\nu_{Jn})^{\alpha-\gamma} D^\alpha \myf(\nu_{Jn}) + \bar{R}\left(D^\gamma \myf(\nu_{Kn}),\nu_{Kn},\nu_{Jn} \right),
\end{align*}
where 
\begin{align*}
r_\alpha(\nu_{Jn},\nu_{Kn},\gamma) = & (m-|\gamma|) \int_0^1 (1-t)^{m-|\gamma|-1} D^\alpha \myf(\nu_{Jn}+t(\nu_{Kn}-\nu_{Jn})) dt  \\
\bar{R}\left(D^\gamma \myf(\nu_{Kn}),\nu_{Kn},\nu_{Jn} \right) = & \sum_{\alpha\geq \gamma,|\alpha|=m} \frac{h_\alpha(\nu_{Jn},\nu_{Kn},\gamma)}{(\alpha-\gamma)!}(\nu_{Kn}-\nu_{Jn})^{\alpha-\gamma}
\end{align*} 
with 
\begin{equation}
\lim_{n\to \infty}h_{\alpha}(\nu_{Jn},\nu_{Kn},\gamma)= \lim_{n\to\infty}(m-|\gamma|) \int_0^1 (1-t)^{m-|\gamma|-1} D^\alpha \myf(\nu_{Jn}+t(\nu_{Kn}-\nu_{Jn})) dt -D^\alpha \myf(\nu_{Jn}) =0 \label{eqn:hlimitmul}
\end{equation} 
by the dominated convergence theorem due to the continuity of $D^\alpha\myf$.  
Now, by the induction hypothesis,
\begin{align*}
GH_{Kn} \myf
= & \sum_{|\gamma|\leq m} a(\gamma|K,\nu_{Kn}) \left(\epsilon_K(n)\right)^{|\gamma|} D^\gamma \myf(\nu_{Kn}) + R(\myf,K,\nu_{Kn}) \\
= & \sum_{|\alpha|\leq m} D^\alpha \myf(\nu_{Jn}) \sum_{\gamma \leq \alpha}  a(\gamma|K,\nu_{Kn}) \left(\epsilon_K(n)\right)^{|\gamma|}   \frac{1}{(\alpha-\gamma)!}(\nu_{Kn}-\nu_{Jn})^{\alpha-\gamma} + \tilde{R}(\myf,K).
\end{align*}
where $\tilde{R}(\myf,K) = R(\myf,K,\nu_{Kn}) + \sum_{|\gamma|\leq m} a(\gamma|K,\nu_{Kn}) \left(\epsilon_K(n)\right)^{|\gamma|} \bar{R}\left(D^\gamma \myf(\nu_{Kn}),\nu_{Kn},\nu_{Jn} \right)  $.  
 Consequently, 
\begin{align*}
& GH_{Jn} \myf \\
= & \sum_{K\in \Child(J)} GH_{Kn} \myf \\
= & \sum_{|\alpha|\leq m} D^\alpha \myf(\nu_{Jn}) (\epsilon_J(n))^{|\alpha|} \sum_{K\in \Child(J)} \sum_{\gamma \leq \alpha}  a(\gamma|K,\nu_{Kn}) \left(\frac{\epsilon_K(n)}{\epsilon_J(n)}\right)^{|\gamma|}   \frac{1}{(\alpha-\gamma)!}\left(\frac{\nu_{Kn}-\nu_{Jn}}{\epsilon_J(n)}\right)^{\alpha-\gamma} + \sum_{K\in \Child(J)} \tilde{R}(\myf,K) \\
\end{align*}
By the inductive hypothesis about $a(\gamma|K,\nu_{Kn})$ and simple calculations using the binomial formula, 
\begin{align*}
    & \sum_{K\in \Child(J)} \sum_{\gamma \leq \alpha}  a(\gamma|K,\nu_{Kn}) \left(\frac{\epsilon_K(n)}{\epsilon_J(n)}\right)^{|\gamma|}   \frac{1}{(\alpha-\gamma)!}\left(\frac{\nu_{Kn}-\nu_{Jn}}{\epsilon_J(n)}\right)^{\alpha-\gamma} \\
    = & \frac{1}{\alpha! (\epsilon_J(n))^{|\alpha|}} \sum_{K\in \Child(J)} \   \sum_{\gamma\leq \alpha} \frac{\alpha!}{\gamma!(\alpha-\gamma)!}  m_\gamma(GH_{Kn}-\nu_{Kn})   (\nu_{Kn}-\nu_{Jn})^{\alpha-\gamma} \\
    = & \frac{1}{\alpha! (\epsilon_J(n))^{|\alpha|}}  \sum_{K\in \Child(J)} \   \sum_{j\in K} \bar{\omega}_{jn} \left(\nu_{jn}-\nu_{Jn}\right)^\alpha\\
    = & a(\alpha|J,\nu_{Jn}). \numberthis \label{eqn:recursionamul}
\end{align*}
In addition, $R(\myf,J,\nu_{Jn})=\sum_{K\in \Child(J)} \tilde{R}(\myf,K)$. We have now represented all the quantities for $J$ in $\eqref{eqn:numeratortylortheoremmul}$ in terms of the corresponding ones of its children vertices. It remains to verify their estimates. \\

\noindent
\textit{Proof of \ref{itema:lem:numeratortaylormul}}: By \eqref{eqn:recursionamul},
$$
a(\alpha|J,\nu_{Jn}) \preccurlyeq \max_{K\in \Child(J)} M_{|\alpha|,K}(\nu_{Kn}). 
$$
Moreover, $M_{p,K}(\nu_{Kn})$ is increasing in $p$, and, for any $p\geq |K|$, $M_{p,K}(\nu_{Kn})\preccurlyeq M_{|K|-1,K}(\nu_{Kn})$ since 
\begin{align*}
& \max_{|K|\leq |\gamma| \leq p}\left|a(\gamma|K,\nu_{Kn}) \left(\frac{\epsilon_K(n)}{\epsilon_{K^\uparrow}(n)}\right)^{|\gamma|} \right| \\
\leq  &  \max_{|K|\leq |\gamma| \leq p}\left|a(\gamma|K,\nu_{Kn})  \right| \left(\frac{\epsilon_K(n)}{\epsilon_{K^\uparrow}(n)}\right)^{|K|} \\
\preccurlyeq  & \max_{ 0\leq 
 |\gamma|< |K|}\left|a(\gamma|K,\nu_{Kn})  \right| \left(\frac{\epsilon_K(n)}{\epsilon_{K^\uparrow}(n)}\right)^{|K|} \\
 \leq  & \left(\frac{\epsilon_K(n)}{\epsilon_{K^\uparrow}(n)}\right) M_{|K|-1,K}(\nu_{Kn}), \numberthis \label{eqn:secondpartislessmul}
\end{align*}
where the ``$\preccurlyeq$'' step follows from the induction hypothesis \ref{itema:lem:numeratortaylormul} for $K$.
 It remains to establish that 
\begin{equation}
\max_{0\leq |\alpha|\leq |J|-1 } |a(\alpha|J,\nu_{Jn})|\succcurlyeq
\max_{K\in \Child(J)}M_{|K|-1,K} (\nu_{Kn}).  \label{eqn:coefficientslowerboumul}
\end{equation}

Write $a(\alpha|J,\nu_{Jn})= \hat{a}(\alpha|J,\nu_{Jn}) + \check{a}(\alpha|J,\nu_{Jn})$ with 
\begin{align*}
    \hat{a}(\alpha|J,\nu_{Jn}) = & \sum_{K\in \Child(J)} \sum_{|\gamma|\leq |K|-1} a(\gamma|K,\nu_{Kn})  \left(\frac{\epsilon_K(n)}{\epsilon_J(n)}\right)^{|\gamma|}    \frac{1}{(\alpha-\gamma)!}\left(\frac{\nu_{Kn}-\nu_{Jn}}{\epsilon_J(n)}\right)^{\alpha-\gamma} 1_{\alpha \geq \gamma}, \numberthis  \label{eqn:majorpartlinearsystemmul} \\
    \check{a}(\alpha|J,\nu_{Jn}) = & \sum_{K\in \Child(J)} \sum_{|K|\leq |\gamma|\leq m } a(\gamma|K,\nu_{Kn})  \left(\frac{\epsilon_K(n)}{\epsilon_J(n)}\right)^{|\gamma|}   \frac{1}{(\alpha-\gamma)!}\left(\frac{\nu_{Kn}-\nu_{Jn}}{\epsilon_J(n)}\right)^{\alpha-\gamma} 1_{\alpha \geq \gamma} \numberthis \label{eqn:minorpartlinearsystemmul}
\end{align*}
Recall the index set $\Ic_m:= \{\alpha\in \Nb^q| |\alpha|\leq m\}$.  
Denote $\hat{a}(J,\nu_{Jn})=\left(\hat{a}(\alpha|J,\nu_{Jn})\right)_{\alpha\in \Ic_{|J|-1}}\in \Rb^{|\Ic_{|J|-1}|}$. Set $\lambda_{K,\gamma}(n)=a(\gamma|K,\nu_{Kn})  \left(\frac{\epsilon_K(n)}{\epsilon_J(n)}\right)^{|\gamma|}$ and $\lambda(n)=(\lambda_{K,\gamma}(n))_{K\in \Child(J),\gamma\in \Ic_{|K|-1}}$. Thus we may view \eqref{eqn:majorpartlinearsystemmul} for $\alpha\in \Ic_{|J|-1}$ in matrix form as $\hat{a}(J,\nu_{Jn}) = A(n) \lambda(n)$, for a suitable matrix $A(n)$ defined as below. 

Set $\psi_K(n):=\frac{\nu_{Kn}-\nu_{Jn}}{\epsilon_J(n)}$ for $K \in \Child(J)$. Then $A(n)=A(\psi_{K_1}(n), \psi_{K_2}(n), \ldots, \psi_{K_{|\Child(J)|}}(n) )$ where $K_i\in \Child(J)$ and the function $A(\cdot,\cdots,\cdot)$ is defined in Lemma \ref{lem:linearalgebra} below (with $j=\Child(J)$, $i$ replacing by $K_i$ and $d_i=|K_i|$). Since for any $K\neq K'\in \Child(J)$
$$\|\psi_K(n)-\psi_{K'}(n)\|_2 = \left\|\frac{\nu_{Kn}-\nu_{K'n}}{\epsilon_J(n)}\right\|_2 \asymp 1, $$ 
we have for any $n$ and for any $K,K'\in \Child (J)$
$$
\|\psi_K(n)-\psi_{K'}(n)\|_2 \geq c
$$
for some positive constant $c$. Notice that for any $K\in \Child(J)$, $\|\psi_K(n)\|_2 \preccurlyeq 1$ by the definition of $\epsilon_J(n)$, which then yields $\psi_K(n)\in B(C)$, the closed ball of radius $C$. So for any $n$,
\begin{align*}
(\psi_{K_1}(n),\ldots,\psi_{K_{|\Child(J)|}}(n)) \in & B\\
:= & \{(\theta_1,\ldots,\theta_{|\Child(J)|}) \in (B(C))^{|\Child(J)|} | \theta_i \in 
\Rb^q, \|\theta_i-\theta_j \|_2\geq c,
\; \forall i\neq j \in [|\Child(J)|]\}, 
\end{align*}
a compact set. By Lemma \ref{lem:linearalgebra}, 
$$ \inf_n \inf_{\|w\|_\infty=1} \|A(n)w\|_{\infty} > 0. $$ 
It then follows that 
$$
\max_{0\leq |\alpha|<|J|}|\hat{a}(\alpha|J,\nu_{Jn})| =\|\hat{a}(J,\nu_{Jn})\|_\infty = \| A(n) \lambda(n)\|_\infty \succcurlyeq  \|\lambda(n)\|_{\infty} = \max_{K\in \Child(J)} M_{|K|-1,K}.
$$

By \eqref{eqn:minorpartlinearsystemmul}, for any $0\leq p<|J|$
$$
|\check{a}(\alpha|J,\nu_{Jn})| \preccurlyeq  \max_{K\in \Child(J)} \  \max_{|K|\leq |\gamma|\leq p}\left|a(\gamma|K,\nu_{Kn}) \left(\frac{\epsilon_K(n)}{\epsilon_{K^\uparrow}(n)}\right)^{|\gamma|} \right|=o\left( \max_{K\in \Child(J)} M_{|K|-1,K} \right)
$$
where the ``$=$'' step follows from \eqref{eqn:secondpartislessmul}. Combining the previous two equations proves \eqref{eqn:coefficientslowerboumul}.

\noindent
\textit{Proof of \ref{itemb:lem:numeratortaylormul}}: 
By \ref{itema:lem:numeratortaylormul} for $J$,
\begin{align*}
    \| a(J,\nu_{Jn}) \|_\infty \asymp & \max\left\{ \max_{K\in \Child(J)} \left| \bar{\omega}_{Kn} \right|, \max_{\substack{K\in \Child(J)\\ K \text{ non-leaf} }}\ \max_{|\gamma|< |K|}\left|a(\gamma|K,\nu_{Kn}) \left(\frac{\epsilon_K(n)}{\epsilon_{J}(n)}\right)^{|\gamma|} \right| \right\}  \\
    \geq & \max\left\{ \max_{K\in \Child(J)} \left| \bar{\omega}_{Kn} \right|, \max_{\substack{K\in \Child(J)\\ K \text{ non-leaf} }} \|a(K,\nu_{Kn})\|_\infty \left(\frac{\epsilon_K(n)}{\epsilon_{J}(n)}\right)^{|J|-1}  \right\} \\
    \succcurlyeq & \max\left\{ \max_{K\in \Child(J)} \left| \bar{\omega}_{Kn} \right|, \max_{\substack{K\in \Child(J)\\ K \text{ non-leaf} }}  \max_{F\in \Desc(K)}\left| \bar{\omega}_{Fn} \left(\frac{\epsilon_{F^\uparrow}(n)}{\epsilon_J(n)}\right)^{|J|-1} \right| \right\}\\
    = 
    & \max_{K\in \Desc(J)}\left| \bar{\omega}_{Kn} \left(\frac{\epsilon_{K^\uparrow}(n)}{\epsilon_J(n)}\right)^{|J|-1} \right|,
\end{align*}
where the ``$\succcurlyeq$'' step follows from  the induction hypothesis \ref{itemb:lem:numeratortaylormul} for $K$.

\noindent
\textit{Proof of \ref{itemc:lem:numeratortaylormul}}: By the formula of $R(\myf,J,\nu_{Jn})$ after \eqref{eqn:recursionamul},
\begin{align*}
  &  \frac{R(\myf,J,\nu_{Jn})}{\left(\epsilon_J(n)\right)^{m}} \\
  = & \sum_{K\in \Child(J)} \sum_{|\gamma|\leq m} a(\gamma|K,\nu_{Kn}) \left( \frac{\epsilon_K(n)}{\epsilon_J(n)} \right)^{|\gamma|}   \frac{ \bar{R}\left(D^\gamma \myf(\nu_{Kn}),\nu_{Kn},\nu_{Jn} \right)}{\left(\epsilon_J(n) \right)^{m-|\gamma|} } \\
 & \quad + \sum_{\substack{K\in \Child(J)\\ K \text{ non-leaf}}} \|a(K,\nu_{Kn})\|_\infty \left( \frac{\epsilon_K(n)}{\epsilon_J(n)} \right)^{m} \frac{R(\myf,K,\nu_{Kn})}{\|a(K,\nu_{Kn})\|_\infty \left(\epsilon_K(n)\right)^{m}} \numberthis \label{eqn:remexp}
 \\
  \leq &  
  \left(\sum_{K\in \Child(J)} M_{m,K}(\nu_{Kn})\right)\   o(1)  +  \|a(J,\nu_{Jn})\|_\infty\ o(1)\\
  \leq & \|a(J,\nu_{Jn})\|_\infty\ o(1) ,
\end{align*}
where the first inequality follows from \eqref{eqn:hlimitmul}, parts \ref{itemb:lem:numeratortaylormul} for $J$ and \ref{itemc:lem:numeratortaylormul} for $K$, and the last inequality follows from \ref{itema:lem:numeratortaylormul} for $J$. 

\myred{
\noindent
\textit{Proof of \ref{itemd:lem:numeratortaylormul}}: Notice that 
\begin{align*}
\sup_{\myf\in \myF} \left|\bar{R}\left(D^\gamma \myf(\nu_{Kn}),\nu_{Kn},\nu_{Jn} \right)\right| 
\leq & \sup_{t\in [0,1]} \left|w(t(\nu_{Kn}-\nu_{Jn})) \right|  \sum_{\alpha\geq \gamma, |\alpha|=m} \frac{1}{(\alpha-\gamma)!} \left|(\nu_{Kn}-\nu_{Jn})^{\alpha-\gamma}\right| \\
\preccurlyeq & \sup_{t\in [0,1]} \left|w(t(\nu_{Kn}-\nu_{Jn})) \right|
(\epsilon_J(n))^{m-|\gamma|}.   \numberthis \label{eqn:barremupp}
\end{align*}

Then following \eqref{eqn:remexp}, 
\begin{align*}
  &  \frac{\sup_{\myf\in \myF} \left| R(\myf,J,\nu_{Jn}) \right|}{\left(\epsilon_J(n)\right)^{m}} \\
   \leq & \sum_{K\in \Child(J)} \sum_{|\gamma|\leq m} |a(\gamma|K,\nu_{Kn})| \left( \frac{\epsilon_K(n)}{\epsilon_J(n)} \right)^{|\gamma|}   \frac{ \sup_{\myf\in \myF} \left|\bar{R}\left(D^\gamma \myf(\nu_{Kn}),\nu_{Kn},\nu_{Jn} \right)\right| }{\left(\epsilon_J(n) \right)^{m-|\gamma|} } \\
 & \quad + \sum_{\substack{K\in \Child(J)\\ K \text{ non-leaf}}} \|a(K,\nu_{Kn})\|_\infty \left( \frac{\epsilon_K(n)}{\epsilon_J(n)} \right)^{m} \frac{\sup_{\myf\in \myF} \left| R(\myf,K,\nu_{Kn})\right|}{\|a(K,\nu_{Kn})\|_\infty \left(\epsilon_K(n)\right)^{m}}  \\
  \preccurlyeq &  
  \left(\sum_{K\in \Child(J)} M_{m,K}(\nu_{Kn})\right)\   \max_{K\in \Child(J)} \sup_{t\in [0,1]} \left|w(t(\nu_{Kn}-\nu_{Jn})) \right|  +  \|a(J,\nu_{Jn})\|_\infty\ o(1)\\
  \leq & \|a(J,\nu_{Jn})\|_\infty\ o(1) ,
\end{align*}
where the "$\preccurlyeq$" follows from \eqref{eqn:barremupp},  parts \ref{itemb:lem:numeratortaylormul} for $J$ and \ref{itemd:lem:numeratortaylormul} for $K$, and the last inequality follows from \ref{itema:lem:numeratortaylormul} for $J$ and the property of $w(\cdot)$.  
}
\end{proof}

\myred{
\begin{lem} \label{lem:linearalgebra}
    Let $j,d_i$ be positive integers. Consider  $\theta_1,\ldots,\theta_j\in \Rb^q$ all distinct. Write $\Ic=\{(i,\gamma)| i\in [j], \gamma\in \Ic_{d_i-1} \}$. Denote $d=\sum_{i\in[j]}d_i$.
   \begin{enumerate}[label=(\alph*)]
   \item \label{lem:linearalgebraa}
   If for any multinomial $P(x)$ of degree $d-1$
   $$
 \sum_{(i,\gamma)\in \Ic} \lambda_{i,\gamma} D^\gamma P(\theta_i) =0, 
   $$
   then
   $$
   \lambda_{i,\gamma} =0, \quad (i,\gamma)\in \Ic. 
   $$
    \item  \label{lem:linearalgebrab}
    Define for each $(i,\gamma)\in \Ic$, a $|\Ic_{d-1}|$-dimensional column vector $a_{i,\gamma}=(a_{i,\gamma}(\alpha))_{\alpha\in \Ic_{d-1}}$ with
    $$
    a_{i,\gamma}(\alpha)= \frac{\theta_i^{\alpha-\gamma}}{(\alpha-\gamma)!}1_{\alpha\geq \gamma},
    $$
    and stack these vectors in a $|\Ic_{d-1}|\times |\Ic|$ matrix $A(\theta_1,\ldots,\theta_j)=(a_{i,\gamma})_{(i,\gamma)\in \Ic}$. Then $A(\theta_1,\ldots,\theta_j)$ is of full column rank. Moreover, for any compact subset $B$ of $\{(\theta_1,\ldots,\theta_j)| \theta_i\in \Rb^q, \theta_i\neq \theta_\ell,\ \forall i\neq \ell \in [j]\}$, 
    $$
    \inf_{(\theta_1,\ldots,\theta_j)\in B}\ \inf_{\|w\|_\infty=1}\|A(\theta_1,\ldots,\theta_j)w\|_\infty>0.$$
    \end{enumerate}
\end{lem}
\begin{proof}
(a)
Fix arbitrary $i\in [j]$. Since $\theta_i\neq \theta_\ell$ for $\ell\neq i$, there exists $\beta_\ell$ with $|\beta_\ell|=d_\ell$ such that $(\theta_i-\theta_\ell)^{\beta_\ell}\neq 0$. (Indeed, say for $m$-th coordinate,  $\theta_i^{(m)}\neq \theta_\ell^{(m)}$, and then one can set $\beta_\ell $ to be $d_i$ on the $m$-th coordinate and zero on other coordinates.) Consider arbitrary $\beta_i$ with $|\beta_i|=d_i-1$. 
    Now, we apply to \eqref{eqn:fullranktemp1} the polynomial 
    $
    P(x)= \prod_{\ell\in [j]} (x-\theta_\ell)^{\beta_\ell}
    $. Note that the highest multi-index power of such $P(x)$ has magnitude $\sum_{i\in[j]}|\beta_i|\leq d-1$.  With this particular choice of $P(x)$, \eqref{eqn:fullranktemp1} implies that $\lambda_{i,\beta_i}=0$. Since $\beta_i$ is arbitrary with $|\beta_i|=d_i-1$, we then have  $\lambda_{i,\gamma}=0$ for any $|\gamma|=d_i-1$. Next if we choose arbitrary $\beta_i$ with $|\beta_i|=d_i-2$ (while keeping the choice of $\beta_\ell$ for $\ell\neq i$), we can obtain  $\lambda_{i,\gamma}=0$ for any $|\gamma|=d_i-2$. Repeating this process yields $\lambda_{i,\gamma}=0$ for any $|\gamma|\leq d_i-1$. Repeating this for $i\in [j]$ completes the proof.\\

    (b)
    Set for short $A=A(\theta_1,\ldots,\theta_j)$. Let $\Lambda=(\lambda_{i,\gamma})_{(i,\gamma)\in \Ic}$ be a column vector such that $A\Lambda=0$. To show that $A$ is of full column rank is equivalent to prove that $\Lambda=0$. Note that for each $\alpha\in \Ic_{d-1}$, 
    $$
    0= (A\Lambda)_{\alpha} = \sum_{(i,\gamma)\in \Ic} \lambda_{i,\gamma} \frac{\theta_i^{\alpha-\gamma}}{(\alpha-\gamma)!}1_{\alpha\geq\gamma}.
    $$
    Then for any multinomial $P(x)=\sum_{\alpha\in \Ic_{d-1}} b_\alpha \frac{x^\alpha}{\alpha!} $, we have 
    \begin{equation}
    0=bA\Lambda = \sum_{\alpha\in \Ic_{d-1}} b_\alpha (A\Lambda)_\alpha = \sum_{(i,\gamma)\in \Ic} \lambda_{i,\gamma} D^\gamma P(\theta_i), \label{eqn:fullranktemp1}
    \end{equation}
    where $b=(b_\alpha)_{\alpha\in \Ic_{d-1}}$. Then by part \ref{lem:linearalgebraa}, $\lambda=0$.

    Consider $f(A)=\inf_{\|w\|_\infty=1}\|Aw\|_\infty$. 
    It is easy to verify that $|f(A)-f(A')|\leq f(A-A')\leq \|A-A'\|_\infty$, and thus $f$ is continuous. 
    Since $A(\theta_1,\ldots,\theta_j)$ is continuous on $(\Rb^q)^j$, $g(\theta_1,\ldots,\theta_j)=f(A(\theta_1,\ldots,\theta_j))$ is continuous. Moreover, $g$ is  positive on $B$ since $A(\theta_1,\ldots,\theta_j)$ is of full column rank. Then $g$ has a positive minimum by compactness of $B$.
\end{proof}
}

\subsection{Optimality of Theorem \ref{thm:inversebound}}
\label{sec:optimality}

In this subsection we show that the exponent $2d_1-1$ of the denominator in \eqref{eqn:localinversebound} is optimal.

\begin{lem} Consider any $k_0\leq k$ and any $G_0 = \sum_{i\in [k_0-1]}p_i^0\delta_{\theta_i^0} + p_{k_0}^0\delta_{\theta_0} \in \Ec_{k_0}(\Theta)$. Suppose each $\myf\in \myF$ is $(2d_1-1)$-th order continuously differentiable on $\{\theta\in \Theta: \|\theta-\theta_0\|_2< b\}$. Suppose furthermore that
$$
  A':=\max_{|\alpha|=2d_1-1}\  \sup_{\substack{\theta'\in\sp(\psi) \\ \|\theta'\|_2\leq b}}\ \sup_{t\in [0,1]}\ \sup_{\myf\in \myF}\  \left|D^{\alpha} \myf (\theta_0+t\theta')\right|<\infty.
$$
Then there exists $G_n\neq H_n\in \Ec_k(\Theta)$ such that $G_n,H_n\overset{W_1}{\to} G_0$ and for any $s<2d_1-1$,
\begin{align*}
    \frac{\sup_{\myf\in \myF}\left|\int \myf(\theta)dG_n  - \int \myf(\theta)dH_n\right|}{W^s_1(G_n,H_n)}\to 0.
\end{align*}
\end{lem}
\begin{proof}
    Let $G,H,G_n,H_n$ be the same as in the proof Lemma \ref{lem:technicallemma}. Then it follows from \eqref{eqn:diffsquareuppbou}, 
\begin{align*}
		&\left|\frac{\int \myf(\theta)dG_n  - \int \myf(\theta)dH_n}{\epsilon_n^{2d_1-1}}\right| 
		\leq  C(d_1,q,A'). 
	\end{align*}
Thus for any $s<2d_1-1$,
\begin{align*}
    \frac{\sup_{\myf\in \myF}\left|\int \myf(\theta)dG_n  - \int \myf(\theta)dH_n\right|}{W^s_1(G_n,H_n)}= \frac{\sup_{\myf\in \myF}\left|\int \myf(\theta)dG_n  - \int \myf(\theta)dH_n\right|}{\left(\epsilon_n p_{k_0}^0W_1(G,H)\right)^s}\to 0.
\end{align*}
\end{proof}

\myred{

\subsection{Proof of Theorem \ref{thm:inverseboundmoment}}
\label{sec:momentinvbou}

\begin{proof}[Proof of Theorem \ref{thm:inverseboundmoment} \ref{itema:thm:inverseboundmoment}] 

The proof is divided into the following steps.

\noindent
\textbf{Step 1}: (Proof by contradiction and subsequences) 
Suppose that \eqref{eqn:localinverseboundmoment} does not hold. Then there exists  $G_n\neq H_n\in \Gc_k(\Theta)$ and $G_n,H_n \overset{W_1}{\to} G_0$ 
such that 
\begin{equation}
\lim_{ n\to\infty} \frac{\sup_{\myf\in \myF} |G_n \myf-H_n \myf|}{\mbf_{2d_1-1}(G_n-\theta_0,H_n-\theta_0)} = 0. \label{eqn:limitzeronew}
\end{equation}

Since $\Theta$ is compact, by taking subsequence if necessary, we have that for each $n$: 1) $G_n\in \Ec_{m_1}(\Theta)$ and $H_n\in \Ec_{m'}(\Theta)$ with $m_1,m'\in [k_0,k]$ independent of $n$; 2) $G_n=\sum_{j\in [m_1]}p_{jn}\delta_{\theta_{jn}}$ and $H_n=\sum_{j \in [m']}\pi_{jn}\delta_{\eta_{jn}}$ with
\begin{align*}
    \sum_{j\in [m_1]} p_{jn} = 1, &  \sum_{j\in [m']} \pi_{jn} = 1, \\
    p_{jn}>0,\ \theta_{jn} \text{ all distinct}, & \quad  \pi_{jn}>0,\ \eta_{jn} \text{ all distinct},   \\
    \theta_{jn} \to \theta_j, & \quad \eta_{jn} \to \eta_j.
    \numberthis \label{eqn:supporpointslimitnew}
\end{align*}
For each $n$, set
$$
(\omega_{jn}, \nu_{jn}) = 
\begin{cases}  
(p_{jn},\theta_{jn}), & \text{ if } j\leq m_1, \\
(-\pi_{(j-m_1)n},\eta_{(j-m_1)n}), & \text{ if } m_1<j\leq m_1+m'.
\end{cases}
$$

\noindent
\textbf{Step 2}: (Decreasing rate of moment difference) 
We will reuse the same notation and definition of the Step 2 and Step 3 in the proof of Theorem \ref{thm:inversebound} \ref{thm:inversebounda}.

\begin{lem}[Characterization of the decreasing rate of moment difference] 
\label{lem:momentdifcharacterization}

If $k_0>1$, or equivalently $\epsilon_{J_r}(n)=1$,  we have 
\begin{align*}
\|\mbf_{2d_1-1}(G_n-\theta_0)- \mbf_{2d_1-1}(H_n-\theta_0)\|_\infty   \asymp & \max_{J\in \Child(J_r)} \max_{|\alpha|\leq |J|-1} \left|a(\alpha|J,\nu_{Jn}) \left(\epsilon_J(n)\right)^{|\alpha|} \right|\\
\asymp & \max_{J\in \Child(J_r)} \max_{|\alpha|\leq |J|-1} \left| m_\alpha(GH_{Jn}-\nu_{Jn}) \right|.
\end{align*}

If $k_0=1$, or equivalently $\epsilon_{J_r}(n)=o(1)$, we have  
\begin{align*}
\|\mbf_{2k-1}(G_n-\theta_0)- \mbf_{2k-1}(H_n-\theta_0)\|_\infty   \asymp &  \max_{|\alpha|\leq 2k-1} \left|a(\alpha|J_r,\nu_{J_r n}) \left(\epsilon_{J_r}(n)\right)^{|\alpha|} \right|\\
\asymp & \max_{|\alpha|\leq 2k-1}  \left| m_\alpha(GH_{J_r n}-\nu_{J_r n}) \right| \\
\asymp &   \|\mbf_{2k-1}(G_n-\nu_{J_r n}))- \mbf_{2k-1}(H_n-\nu_{J_r n}))\|_\infty .
\end{align*}
\end{lem}
\begin{proof}[Proof of Lemma \ref{lem:momentdifcharacterization}]\ \\
\noindent
\textit{Case 1}: Suppose $\epsilon_{J_r}(n)=1$ or equivalently $k_0>1$.
Apply Lemma \ref{lem:numeratortaylormul} to $\myf=\theta^\beta$  with $m=2d_1-1$ for each $J\in \Child(J_r)$, 
$$
m_\beta(G_n-\theta_0)- m_\beta(H_n-\theta_0) = \sum_{J\in \Child(J_r)} \sum_{\alpha\in\Ic_{2d_1-1}} a(\alpha|J,\nu_{J_n}) \left(\epsilon_J(n)\right)^{|\alpha|} \frac{\beta!}{(\beta-\alpha)!} (\nu_{Jn}-\theta_0)^{\beta-\alpha} 1_{\alpha\leq \beta}.
$$

Then 
\begin{align*}
|m_\beta(G_n-\theta_0)- m_\beta(H_n-\theta_0)| 
\leq & 
\sum_{J\in \Child(J_r)} \sum_{\alpha\in\Ic_{2d_1-1}} \left|a(\alpha|J,\nu_{Jn}) \left(\epsilon_J(n)\right)^{|\alpha|} \right| C(d_1,\Theta-\theta_0) \\
\preccurlyeq & \sum_{J\in \Child(J_r)} \max_{|\alpha|\leq |J|-1} \left|a(\alpha|J,\nu_{Jn}) \left(\epsilon_J(n)\right)^{|\alpha|} \right|,
\end{align*}
where the last step follows from Lemma \ref{lem:numeratortaylormul} \ref{itema:lem:numeratortaylormul}. Thus 
\begin{align*}
\|\mbf_{2d_1-1}(G_n-\theta_0)- \mbf_{2d_1-1}(H_n-\theta_0)\|_\infty 
\preccurlyeq & \sum_{J\in \Child(J_r)} \max_{|\alpha|\leq |J|-1} \left|a(\alpha|J,\nu_{Jn}) \left(\epsilon_J(n)\right)^{|\alpha|} \right|.
\end{align*}
By an argument similar to "Proof of \ref{itema:lem:numeratortaylormul} " in the proof of Theorem \ref{lem:numeratortaylormul}, we also have
$$
\|\mbf_{2d_1-1}(G_n-\theta_0)- \mbf_{2d_1-1}(H_n-\theta_0)\|_\infty  \succcurlyeq \max_{J\in \Child(J_r)} \max_{|\alpha|\leq |J|-1} \left|a(\alpha|J,\nu_{Jn}) \left(\epsilon_J(n)\right)^{|\alpha|} \right|.
$$

\noindent
\textit{Case 2}: Suppose $\epsilon_{J_r}(n)=o(1)$ or equivalently $k_0=1$. Apply Lemma \ref{lem:numeratortaylormul} to $\myf=\theta^\beta$  with $m=2d_1-1$ for each $J=J_r$,
$$
m_\beta(G_n-\theta_0)- m_\beta(H_n-\theta_0) = \sum_{\alpha\in\Ic_{2d_1-1}} a(\alpha|J_r,\nu_{J_n}) \left(\epsilon_{J_r}(n)\right)^{|\alpha|} \frac{\beta!}{(\beta-\alpha)!} (\nu_{J_r n}-\theta_0)^{\beta-\alpha} 1_{\alpha\leq \beta}.
$$
The remaining of the proof is similar to case 1 and is thus omitted.
\end{proof}

\noindent 
\textbf{Step 3}: (Deriving contradiction with that $\myF$ is a $(2d_1-1,k_0,k)$ linear independence domain). There are two cases: either
$\epsilon_{J_r}(n)=1$ or $\epsilon_{J_r}(n) = o(1)$. \\

\noindent
\textit{Case 1}: Suppose $\epsilon_{J_r}(n)=1$ or equivalently $k_0>1$. Notice that by Lemma \ref{lem:momentdifcharacterization},
\begin{align*}
    \|\mbf_{2d_1-1}(G_n-\theta_0)- \mbf_{2d_1-1}(H_n-\theta_0)\|_\infty
    \asymp   
    \underbrace{\max_{J\in \Child(J_r)}\ \max_{|\alpha|\leq |J|-1} \left|a(\alpha|J,\nu_{Jn}) \left(\epsilon_J(n)\right)^{|\alpha|}\right|}_{:=d_n}.
\end{align*}

Since $G_n\neq H_n$ and $\epsilon_{J_r}(n)=1$, $\Child(J_r)$ is not empty. 
Since $\epsilon_J(n)=o(1)$ for any $J\in \Child(J_r)$, by Lemma \ref{lem:numeratortaylormul} with $m=|J|-1$ for each $J$, 
\begin{align*}
    \frac{|G_n \myf-H_n \myf|}{\|\mbf_{2d_1-1}(G_n-\theta_0)- \mbf_{2d_1-1}(H_n-\theta_0)\|_\infty} \succcurlyeq & \frac{|G_n \myf-H_n \myf|}{d_n}\\
    = &
     \left|\sum_{J\in \Child(J_r)} \left(\sum_{|\alpha|\leq {|J|-1}}  D^\alpha \myf(\nu_{Jn}) \frac{a(\alpha|J,\nu_{Jn}) \left(\epsilon_J(n)\right)^{|\alpha|}}{d_n} + \frac{R(\myf,J,\nu_{Jn})}{d_n}\right)\right|. \numberthis \label{eqn:childexpansionnew}
\end{align*}
 It follows from Lemma \ref{lem:numeratortaylormul} \ref{itemc:lem:numeratortaylormul} and the condition $m= |J|-1$ that
\begin{align*}
\frac{R(\myf,J,\nu_{Jn})}{d_n}    = o(1). \numberthis \label{eqn:remainderlimit0new}
\end{align*}
By taking subsequence if necessary, we have that 
\begin{equation}
\frac{a(\alpha|J,\nu_{Jn}) \left(\epsilon_J(n)\right)^{|\alpha|}}{d_n} \to b_{J\alpha} \label{eqn:coefficientlimitnew}
\end{equation} 
for some $b_{J\alpha}\in [-1,1]$. Moreover, at least one of $\{b_{J\alpha}\}$ has magnitude $1$. We also have 
\begin{equation}
    \sum_{J\in \Child(J_r)} b_{J\bm{0}} = 0 \label{eqn:sum0new}
\end{equation}
since $\sum_{J\in \Child(J_r)} a(\bm{0}|J,\nu_{Jn})=\sum_{J\in \Child(J_r)}\bar{\omega}_{Jn}=\sum_{j\in [m_1]}p_{jn}-\sum_{j\in [m']}\pi_{jn}=0$. 

Then following \eqref{eqn:limitzeronew},
\begin{align*}
0= &\lim_{ n\to\infty} \frac{\sup_{\myf\in \myF} |G_n \myf-H_n \myf|}{ \|\mbf_{2d_1-1}(G_n-\theta_0)- \mbf_{2d_1-1}(H_n-\theta_0)\|_\infty}\\
\geq & \sup_{\myf\in \myF} \liminf_{ n\to\infty} \frac{ |G_n \myf-H_n \myf|}{ \|\mbf_{2d_1-1}(G_n-\theta_0)- \mbf_{2d_1-1}(H_n-\theta_0)\|_\infty}\\
\succcurlyeq & \sup_{\myf\in \myF} \left| \sum_{J\in \Child(J_r)} \ \sum_{|\alpha|\leq |J|-1} b_{J\alpha} D^\alpha \myf(\nu_J)  \right|, \numberthis \label{eqn:differentiatedomainnew}
\end{align*}
where the last step follows from \eqref{eqn:childexpansionnew}, \eqref{eqn:remainderlimit0new}, \eqref{eqn:coefficientlimitnew} and that $\nu_{Jn}\to \nu_J$, due to our choice of $\nu_{Jn}$ in Lemma \ref{lem:numeratortaylormul}, and the limit $\nu_J$ exists due to \eqref{eqn:supporpointslimitnew}.

Since $\epsilon_{J_r}(n)\asymp 1$, $\nu_J$ for different $J\in \Child(J_r)$ are all distinct. Moreover, by Lemma \ref{lem:treecount},
 $|\Child(J_r)| \in  [k_0,2k-k_0]$ and $\sum_{J\in \Child(J_r)}|J| \in [2k_0,2k]$.
That the equations \eqref{eqn:sum0new} and \eqref{eqn:differentiatedomainnew} hold with at least one $b_{J\alpha}$ nonzero contradicts with the hypothesis that $\myF$ is a $(2d_1-1,k_0,k)$ linear independence domain. \\

\noindent
\textit{Case 2}: $\epsilon_{J_r}(n)=o(1)$ or equivalently $k_0=1$. 
This implies that $G_0=\delta_{\theta}$ for some $\theta\in \Theta$ and $\nu_{jn}\to \theta$ for any $j\in J_r$. Notice that by Lemma \ref{lem:momentdifcharacterization},
\begin{align*}
    \|\mbf_{2k-1}(G_n-\theta_0)- \mbf_{2k-1}(H_n-\theta_0)\|_\infty
    \asymp \underbrace{\max_{|\alpha|\leq 2d_1-1} \left|a(\alpha|J_r,\nu_{J_r n}) \left(\epsilon_{J_r}(n)\right)^{|\alpha|}\right|}_{:=d'_n}   .
\end{align*}

By Lemma \ref{lem:numeratortaylormul}
\begin{align*}
    \frac{G_n \myf-H_n \myf}{\|\mbf_{2k-1}(G_n-\theta_0)- \mbf_{2k-1}(H_n-\theta_0)\|_\infty} \succcurlyeq \sum_{|\alpha|\leq 2d_1-1}  D^\alpha \myf(\nu_{J_r n}) \frac{a(\alpha|J_r,\nu_{J_rn}) \left(\epsilon_{J_r}(n)\right)^{|\alpha|}}{d'_n} + \frac{R(\myf,J_r,\nu_{J_rn})}{d'_n}.
\end{align*}
 The remainder of the proof for this case involves deriving a contradiction, which is done in the same manner as that of Case 1 above.

\end{proof}
}

\myred{
\begin{proof}[Proof of Theorem \ref{thm:inverseboundmoment} \ref{itemc:thm:inverseboundmoment}]
The proof is divided into the following steps.

\noindent
\textbf{Step 1}: (Proof by contradiction and subsequences) 
Suppose that \eqref{eqn:localinverseboundmomentupp} does not hold. Then there exists  $G_n\neq H_n\in \Gc_k(\Theta)$ and $G_n,H_n \overset{W_1}{\to} G_0$ 
such that 
\begin{equation}
\lim_{ n\to\infty} \frac{\sup_{\myf\in \myF} |G_n \myf-H_n \myf|}{\mbf_{2d_1-1}(G_n-\theta_0,H_n-\theta_0)} = \infty. \label{eqn:limitzeronewupp}
\end{equation}

Since $\Theta$ is compact, by taking subsequence if necessary, we have that for each $n$: 1) $G_n\in \Ec_{m_1}(\Theta)$ and $H_n\in \Ec_{m'}(\Theta)$ with $m_1,m'\in [k_0,k]$ independent of $n$; 2) $G_n=\sum_{j\in [m_1]}p_{jn}\delta_{\theta_{jn}}$ and $H_n=\sum_{j \in [m']}\pi_{jn}\delta_{\eta_{jn}}$ with
\begin{align*}
    \sum_{j\in [m_1]} p_{jn} = 1, &  \sum_{j\in [m']} \pi_{jn} = 1, \\
    p_{jn}>0,\ \theta_{jn} \text{ all distinct}, & \quad  \pi_{jn}>0,\ \eta_{jn} \text{ all distinct},   \\
    \theta_{jn} \to \theta_j, & \quad \eta_{jn} \to \eta_j.
\end{align*}
For each $n$, set
$$
(\omega_{jn}, \nu_{jn}) = 
\begin{cases}  
(p_{jn},\theta_{jn}), & \text{ if } j\leq m_1, \\
(-\pi_{(j-m_1)n},\eta_{(j-m_1)n}), & \text{ if } m_1<j\leq m_1+m'.
\end{cases}
$$

We will reuse the same notation and definition of the Step 2 and Step 3 in the proof of Theorem \ref{thm:inversebound} \ref{thm:inversebounda}, and Step 2 in the proof of Theorem \ref{thm:inverseboundmoment} \ref{itemc:thm:inverseboundmoment}. \\

\noindent 
\textbf{Step 2}: (Deriving contradiction). There are two cases: either
$\epsilon_{J_r}(n)=1$ or $\epsilon_{J_r}(n) = o(1)$. \\

\noindent
\textit{Case 1}: Suppose $\epsilon_{J_r}(n)=1$ or equivalently $k_0>1$.

Notice that by Lemma \ref{lem:momentdifcharacterization},
\begin{align*}
    \|\mbf_{2d_1-1}(G_n-\theta_0)- \mbf_{2d_1-1}(H_n-\theta_0)\|_\infty
    \asymp   
    \underbrace{\max_{J\in \Child(J_r)}\ \max_{|\alpha|\leq |J|-1} \left|a(\alpha|J,\nu_{Jn}) \left(\epsilon_J(n)\right)^{|\alpha|}\right|}_{:=d_n}.
\end{align*}

Since $G_n\neq H_n$ and $\epsilon_{J_r}(n)=1$, $\Child(J_r)$ is not empty. 
Since $\epsilon_J(n)=o(1)$ for any $J\in \Child(J_r)$, by Lemma \ref{lem:numeratortaylormul} with $m=|J|-1$ for each $J$, 
\begin{align*}
    \frac{|G_n \myf-H_n \myf|}{\|\mbf_{2d_1-1}(G_n-\theta_0)- \mbf_{2d_1-1}(H_n-\theta_0)\|_\infty} \preccurlyeq & \frac{|G_n \myf-H_n \myf|}{d_n}\\
    = &
     \left|\sum_{J\in \Child(J_r)} \left(\sum_{|\alpha|\leq {|J|-1}}  D^\alpha \myf(\nu_{Jn}) \frac{a(\alpha|J,\nu_{Jn}) \left(\epsilon_J(n)\right)^{|\alpha|}}{d_n} + \frac{R(\myf,J,\nu_{Jn})}{d_n}\right)\right|. 
\end{align*}
Thus 
\begin{align*}
 &    \frac{\sup_{\myf\in\myF}|G_n \myf-H_n \myf|}{\|\mbf_{2d_1-1}(G_n-\theta_0)- \mbf_{2d_1-1}(H_n-\theta_0)\|_\infty} \\
    \preccurlyeq & 
     \sum_{J\in \Child(J_r)} \left(\sum_{|\alpha|\leq {|J|-1}} \sup_{\myf\in\myF} \left| D^\alpha \myf(\nu_{Jn}) \right| \frac{|a(\alpha|J,\nu_{Jn})| \left(\epsilon_J(n)\right)^{|\alpha|}}{d_n} + \frac{\sup_{\myf\in\myF}|R(\myf,J,\nu_{Jn})|}{d_n}\right). \numberthis \label{eqn:childexpansionnewnew}
\end{align*}

 It follows from Lemma \ref{lem:numeratortaylormul} \ref{itemd:lem:numeratortaylormul} and the condition $m= |J|-1$ that
\begin{align*}
\frac{\sup_{\myf\in\myF}|R(\myf,J,\nu_{Jn})|}{d_n}    = o(1). 
\end{align*}
By taking subsequence if necessary, we have that 
\begin{equation*}
\frac{a(\alpha|J,\nu_{Jn}) \left(\epsilon_J(n)\right)^{|\alpha|}}{d_n} \to b_{J\alpha} 
\end{equation*} 
for some $b_{J\alpha}\in [-1,1]$. Plug the above two equations into \eqref{eqn:childexpansionnewnew},
\begin{align*}
  \lim_{n\to\infty}  \frac{\sup_{\myf\in\myF}|G_n \myf-H_n \myf|}{\|\mbf_{2d_1-1}(G_n-\theta_0)- \mbf_{2d_1-1}(H_n-\theta_0)\|_\infty} 
    \preccurlyeq 
    \sup_{|\alpha|\leq 2d_1-1} \sup_{\theta\in \Theta} \sup_{\myf\in\myF} \left| D^\alpha \myf(\theta) \right| \sum_{J\in \Child(J_r)} \sum_{|\alpha|\leq {|J|-1}}  |b_{J\alpha}| <\infty
\end{align*}
which contradicts with \eqref{eqn:limitzeronewupp}.
\\

\noindent
\textit{Case 2}: $\epsilon_{J_r}(n)=o(1)$ or equivalently $k_0=1$. 
This implies that $G_0=\delta_{\theta}$ for some $\theta\in \Theta$ and $\nu_{jn}\to \theta$ for any $j\in J_r$. Notice that by Lemma \ref{lem:momentdifcharacterization},
\begin{align*}
    \|\mbf_{2k-1}(G_n-\theta_0)- \mbf_{2k-1}(H_n-\theta_0)\|_\infty
    \asymp \underbrace{\max_{|\alpha|\leq 2d_1-1} \left|a(\alpha|J_r,\nu_{J_r n}) \left(\epsilon_{J_r}(n)\right)^{|\alpha|}\right|}_{:=d'_n}   .
\end{align*}

By Lemma \ref{lem:numeratortaylormul}
\begin{align*}
   & \frac{\sup_{\myf\in \myF}|G_n \myf-H_n \myf|}{\|\mbf_{2k-1}(G_n-\theta_0)- \mbf_{2k-1}(H_n-\theta_0)\|_\infty} \\
   \preccurlyeq & \sum_{|\alpha|\leq 2d_1-1}  \sup_{\myf\in \myF} |D^\alpha \myf(\nu_{J_r n})| \frac{|a(\alpha|J_r,\nu_{J_rn}) |\left(\epsilon_{J_r}(n)\right)^{|\alpha|}}{d'_n} + \frac{\sup_{\myf\in \myF}|R(\myf,J_r,\nu_{J_rn})|}{d'_n}.
\end{align*}
 The remainder of the proof for this case involves deriving a contradiction, which is done in the same manner as that of Case 1 above. 
\end{proof}
}

\section{Additional material and Proofs for Section \ref{sec:examples}}

\subsection{Additional material for Section \ref{sec:MMD}}

\begin{lem}
\label{lem:injectiveembeddingRd}
    Consider a measurable bounded kernel $\ker(\cdot,\cdot)$. 
    \begin{enumerate}[label=(\alph*)]
        \item 
    The map $\mu:\Mc_b(\Xf,\Xc)\to \Hc$ is injective if and only if 
    $$\int \int \ker(x,y) d\Pb(y) d\Pb(x)>0,\quad \forall \Pb\in \Mc_b(\Xf,\Xc)\setminus\{0\}.  $$
    \item 
    If $(\Xf,\Bc(\Xf))=(\Rb^d,\Bc(\Rb^d))$ and $\ker(\cdot,\cdot)  $ is translation invariant, i.e., $\ker(x,y)=\psi(x-y)$, where
    $\psi: \Rb^d\to \Rb$ is the Fourier transform of a finite nonnegative Borel measure $\Lambda$ on $\Rb^d$:
    \begin{equation*}
    \psi(x) = \int_{\Rb^d} e^{-\ibm x^\top \omega} d\Lambda(\omega).  
    \end{equation*}
    then the map $\mu: \Mc_b(\Rb^d,\Bc(\Rb^d))\to \Hc$ is injective if and only if $\supp(\Lambda)=\Rb^d$.
    \item 
    If $(\Xf,\Bc(\Xf))=(\Rb^d,\Bc(\Rb^d))$ and $\ker(\cdot,\cdot)  $ is a radial kernel, i.e.,  there is a finite nonnegative Borel measure $\nu$ on $[0,\infty)$ such that for all $x,y\in \Rb^d$,
    \begin{equation*}
    \ker(x,y) = \int_{[0,\infty)} e^{-t \|x-y\|_2^2} d\nu(t).  
    \end{equation*}
    then the map $\mu: \Mc_b(\Rb^d,\Bc(\Rb^d))\to \Hc$ is injective if and only if $\supp(\nu)\neq \{0\}$.
    \end{enumerate}
\end{lem}

\begin{proof}[Proof of Lemma \ref{lem:injectiveembeddingRd}]
    (a) Notice that 
    $$
    \|\mu(\Pb)\|_\Hc^2 = \langle \mu(\Pb),\mu(\Pb)  \rangle_\Hc =  \int \mu(\Pb)(x) d\Pb(x) = \int \int \ker(x,y) d\Pb(y) d\Pb(x).
    $$
    So $\mu$ is injective if and only if $\mu(\Pb)=0\in \Hc$ implies $\Pb=0\in \Mc_b(\Xf,\Xc)$, if and only if $$\int \int \ker(x,y) d\Pb(y) d\Pb(x)=0$$ implies $\Pb=0\in \Mc_b(\Xf,\Xc)$.\\
    (b) See \cite[Theorem 6 and Proposition 11]{sriperumbudur2011universality}. \\
    (c) See \cite[Theorem 6 and Proposition 16]{sriperumbudur2011universality}.
\end{proof}

\begin{lem} 
\label{eqn:locationmixture}
Let $f(x)$ be a function on $\Rb$ that is $m$-th order differentiable for every $x\in \Rb$ and that the $j$-th derivative $\frac{d^j}{dx^j}f(x)$ is Lebesgue integrable for any $j\in [m]$. Then the location mixture with density (w.r.t. Lebesgue measure) kernel $p(x \mid \theta)=f(x-\theta)$ is $m$-strongly identifiable. 
\end{lem}
The above lemma is a small improvement of \cite[Theorem 3]{chen1995optimal} or \cite[Theorem 2.4]{heinrich2018strong} in that we remove the assumption that $f$ and its derivatives vanish when $|x|$ approach infinity. 

\begin{proof}[Proof of Lemma \ref{eqn:locationmixture}]
Consider any distinct $\{\theta_i\}_{i\in [\ell] }\subset \Theta$. Assume 
\begin{align}
\sum_{i=1}^{\ell}\  \sum_{j=0}^m a_{ij} \frac{d^j p}{d\theta^j} (x \mid \theta_i )  = 0,  \quad a.e.\ x\in \Rb  \label{eqn:linindoned}
\end{align}
and we want to show that 
$$
a_{ij}=0, \quad \forall\ i\in [\ell], 0\leq j\leq m.
$$
Note that 
$$
\frac{d^j p}{d\theta^j} (x\mid \theta) = (-1)^j \frac{d^j f}{dx^j} (x-\theta).
$$
Plugging the above equation into \eqref{eqn:linindoned}, 
\begin{equation}
\sum_{i=1}^{\ell}\  \sum_{j=0}^m a_{ij} (-1)^j \frac{d^j f}{dx^j} (x-\theta_i)  = 0,  \quad a.e.\ x\in \Rb. \label{eqn:temp316}
\end{equation}
Denote the $\Ff$ to be the Fourier transform, i.e. 
$$
\Ff h (\xi) = \int_\Rb h(x) e^{-2\pi i \xi x} dx.
$$
Now taking Fourier transform on both sides of \eqref{eqn:temp316}, we obtain
\begin{align*}
 0 = &  \sum_{i=1}^{\ell}\  \sum_{j=0}^m a_{ij} (-1)^j  \Ff \frac{d^j f}{dx^j} (x-\theta_i) \\
  = &   \sum_{i=1}^{\ell}\  \sum_{j=0}^m a_{ij} (-1)^j e^{-2\pi i \xi \theta_i} \Ff \frac{d^j f}{dx^j}\\
  = &  \sum_{i=1}^{\ell}\  \sum_{j=0}^m a_{ij} (-1)^j e^{-2\pi i \xi \theta_i} (2\pi i \xi)^j \Ff f
\end{align*}
where the last step follows from that $\frac{d^j f}{dx^j}\in C_0$ for $j\in [m-1]$ and \cite[Theorem 8.22]{folland2013real}. Since $f$ is a probability density, then $\Ff f$ is continuous and $\Ff f (0)=1$. Thus $\Ff f (\xi)>0$ on a neighborhood of $0$, which then implies that on that neighborhood, 
\begin{align*}
 0 = \sum_{i=1}^{\ell}\  \sum_{j=0}^m a_{ij} (-1)^j e^{-2\pi i \xi \theta_i} (2\pi i \xi)^j .
\end{align*}
Since the right hand side (or consider its real and imaginary counterparts) is analytic function of $\xi$, we know the above equation must hold for every $\xi\in \Rb$. It then follows from a similar proof to that of \cite[Theorem 3]{chen1995optimal} that 
$$
a_{ij}=0  \quad \forall\ i\in [\ell], 0\leq j\leq m.
$$
\end{proof}

\subsection{Proofs for Section \ref{sec:MMD} }

\begin{lem}
\label{lem:separatingclass}
    If $\Xf$ is a metrizable, $\Pb,\Qb\in \Mc_b(\Xf,\Bc(\Xf))$ and $\int f d\Pb =\int f d\Qb$ for all $f\in C_b(\Xf)$, then $\Pb=\Qb$. If $\Xf$ is a metric space, then $C_b(\Xf)$ can be replaced by a smaller set of functions: the set of all bounded and Lipchitz functions on $\Xf$.
\end{lem}
\begin{proof}
    The second statement is with $\Pb,\Qb$ being probability measure is proved in \cite[Theorem 9.3.2]{dudley2018real}, but its proof indeed works for $\Pb,\Qb\in \Mc_b(\Xf)$ and the functions constructed in the proof are bounded Lipichitz functions. The first statement follows from the same proof, but for a topological space, Lipchitz functions are not meaningful so we use $C_b(\Xf)$ instead.
\end{proof}

\begin{proof}[Proof of Lemma \ref{lem:injectiveembedding}]
    Assume $\mu(\Pb)=\mu(\Qb)$ and we want to prove that $\Pb=\Qb$. By Lemma \ref{lem:separatingclass} it suffices to show that for  $\forall g\in C_b(\Xf)$, $\int g d\Pb = \int g d\Qb $. Since $C_b(\Xf)\subset \bar{\Hc}$ by assumption, for any $\epsilon>0$ there exists $h\in \Hc$ such that $\sup_{x\in\Xf}|g-h| \leq \epsilon$. Then  by the triangular inequality,
    \begin{align*}
        \left|\int g d\Pb - \int g d\Qb\right| \leq &  \left|\int g d\Pb - \int h d\Pb\right| +  \left|\int h d\Pb - \int h d\Qb\right| +  \left|\int h d\Qb - \int g d\Qb\right| \\
        \leq & 2\epsilon + \left|\int h d\Pb - \int h d\Qb\right| \\
        =&  2\epsilon + \left| \langle \mu(\Pb), h \rangle_{\Hc} - \langle \mu(\Qb), h \rangle_{\Hc}  \right| \\
        = & 2\epsilon,
    \end{align*}
    where the second inequality follows from our choice of $g$, and the first equality follows from the definition of map $\mu(\cdot)$. Since $\epsilon$ is arbitrary, we have $\int g d\Pb = \int g d\Qb $.
\end{proof}

\begin{proof}[Proof of Lemma \ref{lem:exchangealorder}]
    For any $\gamma\in \Ic_{m-1}$, by assumption, for any $0<\|\Delta'\|_2<\Delta'_\theta$:
$$
 \left|\frac{f_2(x) D^\gamma p(x|\theta+\Delta e_i)-f_2(x) D^\gamma p(x \mid \theta)}{\Delta}\right|\leq \|f_2\|_\infty \psi_\theta(x), \quad \lambda-a.e.\ x\in \Xf.
$$
It then follows by dominated convergence theorem and by induction that for any $\alpha\in \Ic_m$,
$$
D^\alpha \Psi(\theta) = D^\alpha \int_{\Xf} f_2(x) p(x \mid \theta) d\lambda =  \int_{\Xf} f_2(x) D^\alpha p(x \mid \theta) d\lambda.
$$

Next by assumption, for any $0<\|\Delta'\|_2<\Delta'_\theta$, and for any $\gamma \in \Ic_{m}\setminus \Ic_{m-1}$, 
$$
\left|f_2(x) D^\gamma p(x|\theta+\Delta' )- f_2(x) D^\gamma p(x \mid \theta)\right| \leq \|f_2\|_\infty \psi'_\theta(x), \quad \lambda-a.e.\ x\in \Xf.
$$ 
It then follows by the dominated convergence theorem that for any $\gamma\in \Ic_m\setminus \Ic_{m-1}$, $D^\gamma \Psi(\theta)$ is continuous.
\end{proof}

\begin{proof}[Proof of Lemma \ref{lem:invbouMMD}]

(a) In lieu of Lemma \ref{thm:inversebound} \ref{thm:inversebounda}, it suffices to show that $\myF_1$ is a $(2d_1-1,k_0,k)$ linear independent domain. 

Consider any member $\myf$ such that $\myf(\theta)=\int f_1(x) d\Pb_\theta(x)$ with $f_1 \in \Fc_1$. Note 
$$|f_1(x)| = |\langle f_1, \ker(\cdot,x)  \rangle_\Hc|\leq \|f_1\|_\Hc  \|\ker(\cdot,x)\|_\Hc \leq \|\ker\|_\infty,$$
so $\|f_1\|_\infty \leq \|\ker\|_\infty $. By Lemma \ref{lem:exchangealorder}, each member in $\myF_1$ is $2d_1-1$ continuously differentiable.

Let $m=2d_1-1$. Consider any integer $\ell\in [k_0,2k-k_0]$, and any vector $(m_1,m_2,\ldots,m_\ell)$ such that $1\leq m_i\leq m+1$ for $i\in [\ell]$ and $\sum_{i=1}^\ell m_i\in [2k_0, 2k]$.  
For any distinct $\{\theta_i\}_{i\in [\ell] }\subset \Theta$, we want to show that the following equations
\begin{subequations}
\begin{align}
\sum_{i=1}^{\ell}\  \sum_{|\alpha|\leq m_i-1} a_{i\alpha} D^\alpha|_{\theta=\theta_i} \int_{\Xf} f_1 d\Pb_\theta  = &0,  \quad \forall f_1\in \Fc_1 \label{eqn:linearinddomainMMDa}\\
\sum_{i\in [\ell]}   a_{i\bm{0}}   = &0, \label{eqn:linearinddomainMMDb}
\end{align}
\end{subequations}
imply that
$$
a_{i\alpha}=0, \quad \forall \ 0\leq |\alpha|< m_i, \ i\in [\ell].
$$
Note that \eqref{eqn:linearinddomainMMDa} is equivalent to 
$$
\sum_{i=1}^{\ell}\  \sum_{|\alpha|\leq m_i-1} a_{i\alpha} D^\alpha|_{\theta=\theta_i} \int_{\Xf} h d\Pb_\theta  = 0,  \quad \forall h \in \Hc,
$$
which implies that 
\begin{equation}
\sum_{i=1}^{\ell}\  \sum_{|\alpha|\leq m_i-1} a_{i\alpha} D^\alpha|_{\theta=\theta_i} \int_{\Xf} \ker(y,x) d\Pb_\theta(x)  = 0,  \quad \forall y \in \Xf.  \label{eqn:lineardomainkertemp}
\end{equation}

By Lemma \ref{lem:exchangealorder}, we have 
$$
D^\alpha|_{\theta=\theta_i} \int_{\Xf} \ker(y,x) d\Pb_\theta = D^\alpha|_{\theta=\theta_i} \int_{\Xf} \ker(y,x) p(x \mid \theta) d\lambda  = \int_{\Xf} \ker(y,x) D^\alpha p(x \mid \theta_i)  d\lambda.
$$
Plugging the above equation into \eqref{eqn:lineardomainkertemp}, one has
\begin{equation}
   \int_{\Xf} \ker(y,x)  \sum_{i=1}^{\ell}\  \sum_{|\alpha|\leq m_i-1} a_{i\alpha} D^\alpha p(x \mid \theta_i)  d\lambda  = 0,  \quad \forall y \in \Xf.  \label{eqn:MMDlineardomaintemp2}
\end{equation}
By assumption A2(m), $D^\alpha p(x \mid \theta_i)$ is integrable w.r.t. dominating measure $\lambda$, hence $\sum_{i=1}^{\ell}\  \sum_{|\alpha|\leq m_i-1} a_{i\alpha} D^\alpha p(x \mid \theta_i)$ is integrable w.r.t. $\lambda$. Consequently the measure $\Qb$ defined by $\frac{d\Qb}{d\lambda} = \sum_{i=1}^{\ell}\  \sum_{|\alpha|\leq m_i-1} a_{i\alpha} D^\alpha p(x \mid \theta_i)$ is a member of $\Mc_b(\Xf,\Xc)$. Then \eqref{eqn:MMDlineardomaintemp2} is the same as 
$$
0=\int_\Xf \ker(y,x) dQ(x) = \langle \ker(y,\cdot), \mu(Q) \rangle_\Hc, \quad \forall y\in \Xf,
$$
which implies that $\mu(\Qb)=0\in \Hc$. By injectivity of $\mu$ on $\Mc_b(\Xf,\Xc)$, $\Qb=0\in \Mc_b(\Xf,\Xc)$, or equivalently,
$$
\sum_{i=1}^{\ell}\  \sum_{|\alpha|\leq m_i-1} a_{i\alpha} D^\alpha p(x \mid \theta_i) = 0,  \quad \lambda-a.e.\ x\in \Xf.
$$
Since $\{p(x \mid \theta)\}_{x\in \Xf}$ is a $(2d_1-1,k_0,k)$ linear independent, 
$$
a_{i\alpha}=0, \quad \forall \ 0\leq |\alpha|< m_i, \ i\in [\ell].
$$

(b) By part \ref{lem:invbouMMDa}, we know that for any $k_0\in [k]$, \eqref{eqn:localinversebound} or \eqref{eqn:localinverseboundMMD} holds for any $G_0\in \Ec_{k_0}(\Theta)$. By Lemma \ref{lem:distinguishableMMD},  $\Gc_k(\Theta)$ is distinguishable by $\myF_1$. Then by Lemma  \ref{lem:localtoglobal}, \eqref{eqn:globalinversebound} or \eqref{eqn:globalinverseboundMMD} holds.

\end{proof}

\begin{proof}[Proof of Lemma \ref{lem:empiricalprocessMMD}]
    \begin{align*}
    \Eb \sup_{f_1 \in \Fc_1}\left|\int f_1 d\Pb_G-\frac{1}{n}\sum_{i\in [n]}f_1(X_i)\right|\leq  \frac{2}{n}\Eb \sup_{f_1 \in \Fc_1}\ \sum_{i\in [n]}\sigma_i f_1(X_i) \leq 2\Eb \sqrt{\frac{ \ker(X_1,X_1)}{n}} \leq \frac{2\|\ker\|_\infty}{\sqrt{n}},
    \end{align*}
    where the first inequality follows from the symmetrization method \cite[Lemma 7.5]{van2014probability} with $\sigma_i$ following $\iidtext$ from Rademacher distribution and independent of $\{X_i\}_{i\in [n]}$, and the second inequality follows from \cite[Lemma 22]{bartlett2002rademacher}.
\end{proof}

\begin{lem}
If a function $f$ on $\Rb$ is differentiable everywhere and $f, f'$ are Lebesgue integrable, then $f\in C_0$, i.e. $f$ is continuous and $\lim_{x\to\infty}f(x)=\lim_{x\to -\infty}f(x)=0$.  
\end{lem}
\begin{proof}
By \cite[Theorem 7.21]{rudin1986real}, for any $x_1<x_2$,
$$
|f(x_2)-f(x_1)| = \left|\int_{x_1}^{x_2} f' dx\right|\leq  \int_{x_1}^{x_2} \left|f' \right| dx \leq \int_{x_1}^{\infty} \left|f' \right| dx,
$$
which converges to $0$ when $x_1\to \infty$. By Cauchy's criteria, $\lim_{x\to\infty}f(x)$ exists. Now since $f$ is Lebesgue integrable, it must hold that  $\lim_{x\to\infty}f(x)=0$. Similarly, one also has $\lim_{x\to -\infty}f(x)=0$. 
\end{proof}

\subsection{Optimality of moment inverse bound}

\begin{lem} \label{lem:optimallocalmomentinversebou}
For any $G_0\in \Ec_{k_0}(\Theta)$ and any $\theta'\in \Rb^q$, there exists $G_n\neq H_n\in \Ec_k(\Theta)$ such that $G_n,H_n\overset{W_1}{\to} G_0$ and 
$
\|\mbf_{2k-2}(G_n-\theta')-\mbf_{2k-2}(H_n-\theta')\|_{\infty}=0.
$
Consequently, 
$$
\liminf_{ \substack{G,H\overset{W_1}{\to} G_0\\ G\neq H\in \Gc_k(\Theta) }}\frac{\|\mbf_{2k-2}(G-\theta')-\mbf_{2k-2}(H-\theta')\|_{\infty}}{W_{2d_1-1}^{2d_1-1}(G,H)} =0.
$$
\end{lem}
\begin{proof}
For any $\gamma\in \Ic_{2k-2}$, consider $\myf_\gamma(\theta)=\left(\theta-\theta'\right)^\gamma$. Then by Lemma \ref{lem:technicallemma}, there exist $G_n\neq H_n\in \Ec_k(\Theta)$ such that $G_n,H_n\overset{W_1}{\to} G_0$ and 
$$
m_\gamma(G_n-\theta') =  \int \myf_\gamma(\theta) dG_n =  \int \myf_\gamma(\theta) dH_n = m_\gamma(H_n-\theta'), \quad \forall \gamma\in \Ic_{2d_1-2},
$$
since $D^\alpha \myf_\gamma = 0$ for any $\alpha\in \Ic_{2d_1-1}$. Thus 
$$
\|\mbf_{2d_1-2}(G_n-\theta')-\mbf_{2d_1-2}(H_n-\theta')\|_{\infty}=0.
$$
\end{proof}

\subsection{Proofs for Section \ref{sec:momentdeviation}}


\begin{definition}[\cite{schudy2012concentration}]
	A random variable $Z$ is called moment bounded with parameter $L > 0$ if for any
	integer $i \geq 1$,
	$$
	\Eb |Z|^i \leq i L \Eb |Z|^{i-1}.
	$$
	A probability family $\{\Pb_{\theta}\}_{\theta\in \Theta}$ on $\Rb$ is uniformly moment bounded with parameter $L$ if  $Z$ is moment bounded with parameter $L$ for each $Z\sim P_\theta$.
\end{definition}

We sometimes write 
$p(x\mid G) :=p_G(x)=\int p(x \mid \theta) dG(\theta)$ for any $G$ a mixing measure on $\tilde{\Theta}$. 

\begin{lem} \label{lem:momentboundedmixture}
	
	\begin{enumerate}[label=(\alph*)]
		\item
		If $\{p(x \mid \theta)\}_{\theta\in\Theta}$ is uniformly moment bounded with parameter $L$, then the family of all mixtures  
		$$
		\left\{\left. p(x\mid G)=\int_\Theta p(x \mid \theta) dG(\theta) \right|G\in \Pc(\Theta) \right\}
		$$ 
		generated from $\{p(x \mid \theta)\}_{\theta\in\Theta}$ is also uniformly moment bounded with the parameter $L$.
		\item Let $p(x \mid \theta)$ be any one of the $6$ families in the NEF-QVF specified with the mean parameter $\theta\in \Theta=[M_1,M_2]$.  The family $\{p(x \mid \theta)\}_{\theta\in \Theta}$ is uniformly moment bounded with parameter $L(p(x \mid \theta),\Theta)$, where  $L(p(x \mid \theta),\Theta)$ is a constant depending on the family of probability kernels $p(x \mid \theta)$ and the constraint $\Theta$.
	\end{enumerate}
\end{lem}
\begin{proof} 
	(a) For any distribution $G$ on $\Theta$, consider $X\sim p(x\mid G)$. $X$ can be thought as being generated from the two steps: $\theta\sim G$ and then $X|\theta \sim p(x \mid \theta)$. Then
	$$
	\Eb_G |X|^j = \Eb_{\theta\sim G} \Eb[|X|^j|\theta] \leq jL \Eb_{\theta\sim G}  \Eb[|X|^{j-1}|\theta] = jL \Eb_G |X|^{j-1},
	$$
	where the inequality follows from that $\{p(x \mid \theta)\}_{\theta\in \Theta}$ is uniformly moment bounded. 
	\\
	(b)
	If $p(x \mid \theta)$ is gaussian, gamma or NEF-GHS, then $p(x \mid \theta)$ is log-concave for each fixed $\theta$. By \cite[Lemma 7.3]{schudy2012concentration}, the single distribution $p(x \mid \theta)$ is moment bounded with parameter $L_{1}=\Eb_\theta|Y|$ where $Y\sim p(x \mid \theta)$. \myred{In particular, for gaussian family, $L_1=\sqrt{\frac{2}{\pi}}\sigma$, independent of $\Theta$. }
	
	If $p(x \mid \theta)$ is Poisson, binomial or negative binomial, then the corresponding random variable is non-negative integer-valued log-concave for each fixed $\theta$.  
 By \cite[Lemma 7.6]{schudy2012concentration}, the single distribution $p(x \mid \theta)$ is moment bounded with parameter $L_{2}=1+\Eb_\theta|Y|$ where $Y\sim p(x \mid \theta)$. 
	
	Combining both cases we see that for a fixed $\theta$, $p(x \mid \theta)$ is moment bounded with $L=1+\Eb_\theta |Y|$  where $Y\sim p(x \mid \theta)$. It is not difficult to see that $L$, as a continuous function of $\theta$ on $\Theta=[M_1,M_2]$ for each given family $p(x \mid \theta)$ in NEF-QVF, has an upper bound $L(p(x \mid \theta),\Theta)$.
\end{proof}

\begin{lem} \label{lem:concentrationbou} 
	Consider any of the $6$ NEF-QVF families \eqref{eqn:NEFQVF} and let $t_j(\cdot|\theta_0)$ and $\bar{t}_j(\theta_0)$ be the same as in Lemma \ref{lem:concentrationboualljNEFQVF} for each specific family of probability kernels $p(x \mid \theta)$. Then for any $\epsilon>0$, and any $G\in \Pc(\Theta)$,
	\begin{align*}
		\Pb_G(|\bar{t}_j(\theta_0)-\Eb_G t_j(X|\theta_0)|\geq \epsilon) 
		\leq  e^2 \exp\left( -C(p(x \mid \theta), \Theta,j,\theta_0) \min \left\{ n\epsilon^2,  (n\epsilon)^{1/j} \right\}  \right),
	\end{align*}
	where the positive constant $C(p(x \mid \theta), \Theta,j,\theta_0)$  depends on the specific NEF-QVF family $p(x \mid \theta)$, the constraint $\Theta$, the polynomial degree $j$ and the choice of a reference point $\theta_0\in \tilde{\Theta}^\circ$.
\end{lem}

\begin{proof}
	Denote $\tilde{t}_j(x_1,\ldots,x_n|\theta) := \frac{1}{n}\sum_{i\in [n]} t_j(x_i|\theta)$. Then $\bar{t}_j(\theta_0)= \tilde{t}_j(X_1,\ldots,X_n|\theta_0)$. The proof follows by a general concentration inequality for independent random variables \cite[Theorem 1.4]{schudy2012concentration}. By Lemma \ref{lem:momentboundedmixture} $\{p(x\mid G)\}_{G\in \Gc_k(\Theta)}$ is uniformly moment bounded with parameter $L=L(p(x \mid \theta),\Theta)$.

	We now calculate the constants in the upper bound of \cite[Theorem 1.4]{schudy2012concentration}. For $\tilde{t}_j(x_1,\ldots,x_n|\theta_0)$, the total power is $q=j$ and the maximal variable power is $\Gamma = j$. To avoid notation conflict, we write $\nu_r$ for the $\mu_r$ in \cite[Theorem 1.4]{schudy2012concentration}. By \cite[(1.7)]{schudy2012concentration} we have 
	\begin{align*}
		\nu_r = & \frac{1}{n} |a_{jr}(\theta_0)|, \quad \forall r\in [j], \\  
		\nu_0 = & \frac{1}{n}\sum_{i\in [n]} \sum_{\ell\in [j] } |a_{j\ell}(\theta_0)| \Eb_G |X_i|^\ell= \sum_{\ell\in [j] } |a_{j\ell}(\theta_0)| \Eb_G |X_1|^\ell,
	\end{align*}
 where we recall $a_{ji}(\theta_0)$ are the coefficients of $t_j(x|\theta_0)$ defined after \eqref{eqn:gammamean}. 
	By \cite[Theorem 1.4]{schudy2012concentration},
	with the notion $w_r := n \nu_r=|a_{jr}(\theta_0)|$,
	\begin{equation}
		\Pb_G(|\bar{t}_j(\theta_0)-\Eb_G t_j(X|\theta_0)|\geq \epsilon) \leq e^2 \max\left\{ \max_{r\in [j]} \exp\left(-\frac{n\epsilon^2}{\nu_0 w_r L^r j^r R^j}\right), \max_{r\in [j]} \exp\left(-\left(\frac{n\epsilon}{w_r L^r j^r R^j}\right)^{1/r}\right)  \right\}, \label{eqn:momentboundedbou}
	\end{equation}
	where $R\geq 1$ is some absolute constant.
	
	Note that 
	\begin{align}
		\min_{r\in [j]} \frac{n\epsilon^2}{\nu_0 w_r L^r j^r R^j} =  \frac{n\epsilon^2}{\nu_0 R^j \max_{r\in [j]} w_r L^r j^r }, \label{eqn:momentboundedconcentrationfirstterm}
	\end{align}
	and
	\begin{align}
		\min_{r\in [j]}  \left(\frac{n\epsilon}{w_r L^r j^r R^j}\right)^{1/r} \geq \frac{1}{Lj   \max_{r\in[j]}(w_r R^j)^{1/r}} \ \min_{r\in[j]}(n\epsilon)^{1/r}. \label{eqn:momentboundedconcentrationsecondterm} 
	\end{align}
        \myred{
        By denoting
        $$
        A:= \max\left\{\nu_0 R^j \max_{r\in [j]} w_r L^r j^r , Lj   \max_{r\in[j]}(w_r R^j)^{1/r} \right\},
        $$
        the inequality \eqref{eqn:momentboundedbou}  becomes
        \begin{equation}
		\Pb_G(|\bar{t}_j(\theta_0)-\Eb_G t_j(X|\theta_0)|\geq \epsilon) \leq e^2  \exp\left(- \frac{1}{A }\min\left\{ n\epsilon^2,n\epsilon,  (n\epsilon)^{1/j} \right\}\right)  . \label{eqn:momentboundedbou2}
	\end{equation}

	Moreover, since $\{p(x\mid G)\}_{G\in \Gc_k(\Theta)}$ is uniformly moment bounded, 
	\begin{equation}
		\nu_0\leq |a_{j0}(\theta_0)|+   \sum_{\ell=1}^{j} |a_{j\ell}(\theta_0)|L^{\ell-1} \Eb_G |X_1| \leq \max\{\Eb_G |X_1|,1\} \left( |a_{j0}(\theta_0)|+   \sum_{\ell=1}^{j} |a_{j\ell}(\theta_0)|L^{\ell-1} \right). \label{eqn:nu0bou1} 
	\end{equation}
	Write  $Y_\theta \sim f(x|\theta)$. Then
	\begin{equation}
		\Eb_G |X_1| =  \int_{\Theta} \Eb |Y_\theta| dG \leq L
  \label{eqn:absfirstmomentbou}
	\end{equation}
 where the last step follows from the definition of $L$ in the proof of Lemma \ref{lem:momentboundedmixture}. 
 Combining \eqref{eqn:nu0bou1} and \eqref{eqn:absfirstmomentbou}, we obtain $\nu_0\leq C(p(x \mid \theta), L,j,\theta_0)$, where the dependence of $t_j(\cdot|\theta_0)$ is absorbed in the dependence on $p(x \mid \theta)$, $j$ and $\theta_0$. Therefore, $A\leq C(p(x \mid \theta), L,j,\theta_0)$, and hence \eqref{eqn:momentboundedbou2} becomes:
	\begin{align*}
		\Pb_G(|\bar{t}_j(\theta_0)-\Eb_G t_j(X|\theta_0)|\geq \epsilon)
		\leq & e^2 \exp\left( -C(p(x \mid \theta),L,j,\theta_0) \min\left\{ n\epsilon^2,n\epsilon,  (n\epsilon)^{1/j} \right\}  \right) \\
		\leq & e^2 \exp\left( -C(p(x \mid \theta),L,\Theta,j,\theta_0) \min \left\{ n\epsilon^2,  (n\epsilon)^{1/j} \right\}  \right).
	\end{align*}
}
\end{proof}

\begin{proof}[Proof of Lemma \ref{lem:concentrationboualljNEFQVF}]
\myred{
 By Lemma \ref{lem:concentrationbou},
	\begin{align*}
		\Pb_G\left(\max_{j\in [2k-1]}|\bar{t}_j(\theta_0)-\Eb_G t_j(X|\theta_0)|\geq \lambda\right) 
		\leq & \sum_{j\in [2k-1]} e^2 \exp\left( -C(p(x \mid \theta), L,j,\theta_0) \min \left\{ n\lambda^2,  (n\lambda)^{1/j} \right\}  \right)\\
		\leq &  (2k-1) e^2 \exp\left( -C(p(x \mid \theta), L,k,\theta_0) \min \left\{ n\lambda^2,  (n\lambda)^{\frac{1}{2k-1}} \right\}  \right) ,
	\end{align*}
	where in the second inequality $C(p(x \mid \theta), L,k,\theta_0)=\min_{j\in [2k-1]} C(p(x \mid \theta), L,j,\theta_0) $. 
 }
\end{proof}

\begin{lem} \label{lem:highdimnortalp}
Consider $X\sim \Nc(U,\Sigma)$ on $\Rb^d$. 
\begin{enumerate}[label=(\alph*)]
\item 
Then 
$$
\|M_\ell (X) \|_2 \leq \left(  \|  U\|_2+   \|\Sigma^{\frac{1}{2}}\|_2  C \sqrt{\ell} \right)^\ell,
$$
where $C$ is a universal constant. 
\item  
For any $\beta\in [d]^\ell$, 
$$
\mathrm{Var} \left( \prod_{j=1}^\ell X_{\beta_j} \right)  \leq \left(  \|  U\|_2+   \|\Sigma^{\frac{1}{2}}\|_2  C \sqrt{2\ell} \right)^{2\ell}
$$
where $C$ is a universal constant. 
\item  \label{lem:highdimnortalpc}
Consider $\alpha\in \Omega_{d,\ell}$. Then
$$
\mathrm{Var}(t_\alpha(X) )  \leq \left(  \|  U\|_2+    \| \Sigma^{\frac{1}{2}} \|_2 C \sqrt{\ell}  \right)^{2\ell}
$$
where $C$ is a universal constant.
\end{enumerate}
\end{lem}
\begin{proof}
(a)
Write $X= U + \Sigma^{\frac{1}{2}} Z$ where  $Z\sim \Nc(\bm{0}, I)$. 
\begin{align*}
\|M_\ell (X) \|_2 = & \sup_{\|v\|_2=1} \Eb |\langle X,  v \rangle|^\ell \\
= & \sup_{\|v\|_2=1}  \| \langle X,  v \rangle \|_{L^\ell}^\ell \\
\overset{(*)}{\leq} & \sup_{\|v\|_2=1}  \left(| \langle U,  v \rangle |+\| \langle \Sigma^{\frac{1}{2}} Z,  v \rangle \|_{L^\ell}\right)^\ell \\
 \leq &  \left(  \|  U\|_2+  \sup_{\|v\|_2=1} \|\Sigma^{\frac{1}{2}}\|_2  \|  \langle  Z,  v \rangle \|_{L^\ell}\right)^\ell \\
 \leq & \left(  \|  U\|_2+   \|\Sigma^{\frac{1}{2}}\|_2  C \sqrt{\ell} \right)^\ell,
\end{align*}
where (*) follows from triangular inequality of $\|\cdot\|_{L^\ell}$ and that $U$ is deterministic, and the last inequality follows that $\langle  Z,  v \rangle$ is standard normal and $C$ is an universal constant. \\

(b)
$$
\mathrm{Var} \left( \prod_{j=1}^\ell X_{\beta_j} \right)  \leq \Eb \prod_{j=1}^\ell X^2_{\beta_j} = \left(M_{2\ell}(X)  \right)_{(\beta,\beta)} \leq \| M_{2\ell}(X) \|_2
$$
where the last inequality follows that the spectrum norm of a tensor is larger than every entry. The proof is then completed by utilizing part (a). \\

(c)
 Choose any $\beta\in \pi_\ell^{-1}(\alpha)$. Then 
\begin{equation}
\mathrm{Var}(t_\alpha(X) ) =  \mathrm{Var}( (F_\ell(X))_\beta ). \label{eqn:normaltalp1}
\end{equation}
Since standard deviation is the $L^2$ norm of centered random variables, by the triangle inequality, 
\begin{equation}
\sqrt{\mathrm{Var}( (F_\ell(X))_\beta )} \leq  \sum_{j=0}^{\lfloor \ell/2 \rfloor } A_{\ell,j} \sqrt{\mathrm{Var} \left(\operatorname{sym}_{\beta}\left( X^{\otimes \ell-2j} \otimes \Sigma^{\otimes j} \right) \right)}, \label{eqn:normaltalp2}
\end{equation}
and, 
\begin{align*}
\sqrt{\mathrm{Var} \left(\operatorname{sym}_{\beta}\left( X^{\otimes \ell-2j} \otimes \Sigma^{\otimes j} \right) \right)} 
\leq & \frac{1}{\ell !} \sum_{\sigma\in S_\ell} \sqrt{\mathrm{Var} \left(  \left(\prod_{i=1}^{\ell-2j} X_{\beta_{\sigma(i)}}\right) \left( \prod_{i=1}^j \Sigma_{\beta_{ \sigma( \ell-2j + 2i-1)} \beta_{\sigma(\ell-2j + 2i)} } \right)  \right) } \\
\leq & \frac{1}{\ell !} \sum_{\sigma\in S_\ell} \left( \prod_{i=1}^j \Sigma_{\beta_{ \sigma( \ell-2j + 2i-1)} \beta_{\sigma(\ell-2j + 2i)} } \right) \sqrt{\mathrm{Var} \left(  \left(\prod_{i=1}^{\ell-2j} X_{\beta_{\sigma(i)}}\right)   \right) } \\
\leq & \frac{1}{\ell !} \sum_{\sigma\in S_\ell} \| \Sigma \|_2^{j} \sqrt{\mathrm{Var} \left(  \left(\prod_{i=1}^{\ell-2j} X_{\beta_{\sigma(i)}}\right)   \right) } \\
\leq &  \| \Sigma \|_2^{j} \left(  \|  U\|_2+   \|\Sigma^{\frac{1}{2}}\|_2  C \sqrt{2(\ell-2j)} \right)^{(\ell-2j)},  \numberthis \label{eqn:normaltalp3}
\end{align*}
where in the first inequality $S_\ell$ denotes the set of all permutations on $[\ell]$, and the last step follows from part (b). 

By combining \eqref{eqn:normaltalp1}, \eqref{eqn:normaltalp2} and \eqref{eqn:normaltalp3}, 
\begin{align*}
\sqrt{\mathrm{Var}(t_\alpha(X) )}
\leq &  \sum_{j=0}^{\lfloor \ell/2 \rfloor } A_{\ell,j}  \| \Sigma \|_2^{j} \left(  \|  U\|_2+   \|\Sigma^{\frac{1}{2}}\|_2  C \sqrt{2(\ell-2j)} \right)^{(\ell-2j)}  \\
\leq & \sum_{j=0}^{\lfloor \ell/2 \rfloor } A_{\ell,j}  \| \Sigma^{\frac{1}{2}} \|_2^{2j} \left(  \|  U\|_2+   \|\Sigma^{\frac{1}{2}}\|_2  C \sqrt{2\ell} \right)^{\ell-2j} \\
= & \Eb \left(  \|  U\|_2+   \|\Sigma^{\frac{1}{2}}\|_2  C \sqrt{2\ell} + \| \Sigma^{\frac{1}{2}} \|_2 Z_1  \right)^{\ell} \\
\leq & \left(  \|  U\|_2+   \|\Sigma^{\frac{1}{2}}\|_2  C \sqrt{2\ell} + \| \Sigma^{\frac{1}{2}} \|_2 C \sqrt{\ell}  \right)^{\ell}
\end{align*}
where the equality follows from the formula of moments of one-dimensional Gaussian distribution, and the last inequality follows from part (a) when dimension $d=1$. 

\end{proof}

\begin{proof}[Proof of Lemma \ref{lem:tensorcon}]
Write $X_i= U_i + \Sigma^{\frac{1}{2}} Z_i$ where $U_i\sim G=\sum_{i\in [k]} p_i \delta_{\theta_i}$ and $Z_i \sim \Nc(\bm{0}, I)$.  Denote $R=\sup_{\theta\in \Theta}\|\theta\|_2$. 

Consider $\alpha\in \Omega_{d,\ell}$.  Denote $\bar{t}_\alpha=\frac{1}{n} \sum_{i\in [n]}t_\alpha(X_i)$.  Then by independence,
$$
\mathrm{Var}(\bar{t}_\alpha|U_1,\ldots,U_n) =  \frac{1}{n^2} \sum_{i=1}^n \mathrm{Var}(t_\alpha(X_i)|U_i)\leq \frac{1}{n^2} \sum_{i=1}^n \left(  \|  U_i\|_2+    \| \Sigma^{\frac{1}{2}} \|_2 C \sqrt{\ell}  \right)^{2\ell}\leq \frac{1}{n}  \left(  R +    \| \Sigma^{\frac{1}{2}} \|_2 C \sqrt{\ell}  \right)^{2\ell},
$$
where the first inequality follows from Lemma \ref{lem:highdimnortalp} \ref{lem:highdimnortalpc}. By the Hypercontractivity  inequality  \cite[Theorem  1.9]{schudy2012concentration},  
$$
\Pb\left( \left| \bar{t}_\alpha - \Eb[ \bar{t}_\alpha|   U_1,\ldots,U_n ] \right| > \epsilon | U_1,\ldots,U_n   \right) \leq e^2 \exp \left( - \left( \frac{c n \epsilon^2}{\left(  R +    \| \Sigma^{\frac{1}{2}} \|_2 C \sqrt{\ell}  \right)^{2\ell}}  \right)^{\frac{1}{\ell}}  \right).
$$
By taking expectation on both sides, 
\begin{equation}
\Pb\left( \left| \bar{t}_\alpha - \Eb[ \bar{t}_\alpha|   U_1,\ldots,U_n ] \right| > \epsilon    \right) \leq e^2 \exp \left( - \left( \frac{c n \epsilon^2}{\left(  R +    \| \Sigma^{\frac{1}{2}} \|_2 C \sqrt{\ell}  \right)^{2\ell}}  \right)^{\frac{1}{\ell}}  \right). \label{eqn:polgautail}
\end{equation}

Choose any $\beta\in \pi_\ell^{-1}(\alpha)$. Then by Lemma \ref{lem:denoisedmomenttensor}, 
$$
\Eb[ \bar{t}_\alpha|   U_1,\ldots,U_n ] = \frac{1}{n} \sum_{i=1}^n \Eb(t_\alpha(X_i)|U_i) =  \frac{1}{n} \sum_{i=1}^n \Eb( (F_\ell(X_i))_\beta|U_i) = \frac{1}{n} \sum_{i=1}^n (U_i^{\otimes \ell} )_\beta = \frac{1}{n} \sum_{i=1}^n \prod_{j=1}^\ell U_{i\beta_j}. 
$$
Since 
$$
\left|\prod_{j=1}^\ell U_{i\beta_j} \right|\leq R^\ell,
$$
by Hoeffding's inequality, 
\begin{equation}
\Pb\left( \left|   \Eb[ \bar{t}_\alpha|   U_1,\ldots,U_n ] -\Eb \bar{t}_\alpha \right| > \epsilon   \right) \leq 2 \exp\left( -\frac{2n\epsilon^2}{ (2R^\ell)^2 }  \right) = 2 \exp\left( -\frac{n\epsilon^2}{ 2R^{2\ell} }  \right). \label{eqn:conhoe}
\end{equation}
Combining \eqref{eqn:polgautail} and \eqref{eqn:conhoe}, 
\begin{align*}
\Pb\left( \left| \bar{t}_\alpha -\Eb \bar{t}_\alpha \right| > 2\epsilon   \right) 
\leq & e^2 \exp \left( - \left( \frac{c n \epsilon^2}{\left(  R +    \| \Sigma^{\frac{1}{2}} \|_2 C \sqrt{\ell}  \right)^{2\ell}}  \right)^{\frac{1}{\ell}}  \right) + 2 \exp\left( -\frac{n\epsilon^2}{ 2R^{2\ell} }  \right) \\
\leq & e^2  \exp \left( -C(R,\|\Sigma^{\frac{1}{2}}\|_2, \ell ) \min\{ n\epsilon^2, (n\epsilon^2)^{\frac{1}{\ell}} \} \right).
\end{align*}

Thus 
\begin{align*}
&\Pb\left( \max_{\ell\in [2k-1]}  \left\|M_\ell(G) - \frac{1}{n} \sum_{i\in [n]}F_\ell(X_i)   \right\|_\infty > 2\epsilon \right) \\
\leq & \sum_{\ell=1}^{2k-1} |\Omega_{d,\ell}| e^2  \exp \left( -C(R,\|\Sigma^{\frac{1}{2}}\|_2, \ell ) \min\{ n\epsilon^2, (n\epsilon^2)^{\frac{1}{\ell}} \} \right) \\
\leq & C(d,k) \exp \left( -C(R,\|\Sigma^{\frac{1}{2}}\|_2, k ) \min\{ n\epsilon^2, (n\epsilon^2)^{\frac{1}{2k-1}} \} \right). 
\end{align*}
\end{proof}

\section{Proofs for Section \ref{sec:pointwiserate}}

\begin{proof}[Proof of Lemma \ref{lem:componentnumberestimation}]
    Define $A_G(a_n):= \{ \sup_{\myf\in\myF}\left|G  \myf-\frac{1}{n}\sum_{i\in [n]}t_\myf(X_i)\right| \leq a_n \}$. Then 
    on the event $A_G(a_n)$, we have $\hat{k}_n\leq k(G)$ by the definition of $\hat{k}_n$. We also have
    \begin{align*}
         \{\hat{k}_n< k(G)\}  
        =  \left\{\sup_{\myf\in\myF}\left|\hat{G}_n(k(G)-1)  \myf-\frac{1}{n}\sum_{i\in [n]}t_\myf(X_i)\right|\leq a_n\right\} , \numberthis \label{eqn:eventequality}
    \end{align*}
    since $\sup_{\myf\in\myF}\left|\hat{G}_n(\ell)  \myf-\frac{1}{n}\sum_{i\in [n]}t_\myf(X_i)\right|$ is decreasing w.r.t. $\ell$.

     Following the definition of $b_G$ and the triangle inequality, we have
    \begin{equation}
    b_G \leq \sup_{\myf\in\myF}\left|\hat{G}_n(k(G)-1)  \myf-\frac{1}{n}\sum_{i\in [n]}t_\myf(X_i)\right| + \sup_{\myf\in\myF}\left|G  \myf-\frac{1}{n}\sum_{i\in [n]}t_\myf(X_i)\right|. \label{eqn:bGuppbou}
    \end{equation}
    Combining \eqref{eqn:bGuppbou} and \eqref{eqn:eventequality},  
    $$
    \{\hat{k}_n\neq k(G)\} \bigcap A_G(a_n) \subset \left\{ \sup_{\myf\in\myF}\left|G  \myf-\frac{1}{n}\sum_{i\in [n]}t_\myf(X_i)\right| \geq b_G-a_n \right\} \bigcap A_G(a_n).
    $$
    Thus,
    \begin{align*}
    \{\hat{k}_n\neq k(G)\} \subset & \left(\{\hat{k}_n\neq k(G)\} \bigcap A_G(a_n)\right) \bigcup \left(A_G(a_n)\right)^c \\
    \subset & \left\{ \sup_{\myf\in\myF}\left|G  \myf-\frac{1}{n}\sum_{i\in [n]}t_\myf(X_i)\right| \geq \min\{a_n,b_G-a_n\} \right\}. \numberthis \label{eqn:mixturecomponenttemp1}
    \end{align*}
\end{proof}

\begin{proof}[Proof of Lemma \ref{lem:inverseboundfixoneargument}]

Suppose that \eqref{eqn:localinverseboundfixargument} does not hold. Then there exists a sequence of $G_n\in \Gc_k(\Theta)$ and $G_n \overset{W_1}{\to} G_0$ such that 
\begin{equation}
\lim_{ n\to\infty} \frac{\sup_{\myf\in \myF} |G_n \myf-G_0 \myf|}{W_{2}^{2}(G_n,G_0)} = 0. \label{eqn:limitzerofixargument}
\end{equation}

Write $G_0=\sum_{i\in [k_0]}p_i^0\delta_{\theta_i^0}\in \Ec_{k_0}(\Theta)$. Since $\Theta$ is compact, by taking subsequence if necessary, we have that: 
1) $G_n\in \Ec_{k_*}(\Theta)$ for some $k_*\in [k_0,k]$ independent of $n$; 2) $G_n=\sum_{i\in [k_1]}\sum_{j\in [s_i]} p_{i j n} \delta_{\theta_{i j n}}$, with $k_1,s_i$ all independent of $n$ and
\begin{align*}
    \sum_{i\in [k_1]}\sum_{j\in [s_i]} 1 = k_*, &  \quad
    p_{i j n}>0,\ \theta_{i j n} \text{ all distinct for} \  j\in [s_i],   i\in [k_1],  \\
    \theta_{i j n} \to \theta_i^0,   & \quad \sum_{j\in [s_i]}p_{i j n} \to p_i^0, \quad \forall \ i\in [k_0] \\
    \theta_{i j n} \to \theta_i,  & \quad p_{i j n}\to 0 , \quad \forall\ k_0 < i \leq k_1, 
\end{align*}
where $\{\theta_i\}_{i=k_0+1}^{k_1}$ are distinct elements in $\Theta \setminus \{\theta_i^0\}_{i\in [k_0]}$. 

Note that
\begin{align*}
G_n \myf-G_0 \myf =&\sum_{i\in [k_0]}\underbrace{\sum_{j\in [s_i]} p_{i j n} \left(\myf(\theta_{i j n}) -\myf(\theta_i^0)\right) }_{:=I_i}+ \sum_{i\in [k_0]} \left(\sum_{j\in [s_i]} p_{i j n}- p_i^0\right) \myf(\theta_i^0) + \sum_{i=k_0+1}^{k_1}\sum_{j\in [s_i]} p_{i j n} \myf(\theta_{i j n}).  
\end{align*}
Denote $\If_{>1}:=\{i\in [k_0]: s_i>1 \}$. By Taylor's theorem, for any $i\in \If_{>1}$,
\begin{align*}
I_i = & \sum_{j\in [s_i]} p_{i j n} \left(\sum_{1\leq |\alpha|\leq 2}\frac{1}{\alpha!} D^\alpha \myf (\theta_i^0) \left( \theta_{i j n}-\theta_i^0 \right)^\alpha + R_{i j n}  \right) \\
= &  \sum_{1\leq |\alpha|\leq 2}\frac{1}{\alpha!} D^\alpha \myf (\theta_i^0)  \sum_{j\in [s_i]} p_{i j n} \left( \theta_{i j n}-\theta_i^0 \right)^\alpha + \sum_{j\in [s_i]} p_{i j n} R_{i j n},
\end{align*}
where $R_{i j n} = o(\|\theta_{i j n}-\theta_i^0\|_2^2)$. 

Denote $\If_{1}:=\{i\in [k_0]: s_i=1 \}$. By Taylor's theorem, for any $i\in \If_{1}$,
\begin{align*}
I_i = &  p_{i 1 n} \left(\sum_{ |\alpha|= 1}\frac{1}{\alpha!} D^\alpha \myf (\theta_i^0) \left( \theta_{i j n}-\theta_i^0 \right)^\alpha + R_{i j n}  \right) \\
= &  \sum_{|\alpha|= 1}\frac{1}{\alpha!} D^\alpha \myf (\theta_i^0)   p_{i 1 n} \left( \theta_{i 1 n}-\theta_i^0 \right)^\alpha +  p_{i 1 n}R_{i 1 n},
\end{align*}
where $R_{i 1 n} = o(\|\theta_{i 1 n}-\theta_i^0\|_2)$. 

We also have
\begin{align*}
    W_2^2(G_n,G_0) 
\leq & \sum_{i\in [k_0]}\sum_{j\in [s_i]} p_{i j n} \|\theta_{i j n}-\theta_i^0\|_2^2 + \diam^2(\Theta) \left( \sum_{i\in [k_0]} \left|\sum_{j\in [s_i]} p_{i j n}- p_i^0\right|+ \sum_{i=k_0+1}^{k_1}\sum_{j\in [s_i]} p_{i j n} \right) \\
\leq & \sum_{i\in \If_{>1}}  \max\left\{\max_{|\alpha|=1}\left|\sum_{j\in [s_i]} p_{i j n} \left( \theta_{i j n}-\theta_i^0 \right)^\alpha\right|, \sum_{j\in [s_i]} p_{i j n} \|\theta_{i j n}-\theta_i^0\|_2^2\right\} +\\
& +\sum_{i\in \If_{1}} p_{i 1 n} \|\theta_{i 1 n}-\theta_i^0\|_2 + \diam^2(\Theta) \left( \sum_{i\in [k_0]} \left|\sum_{j\in [s_i]} p_{i j n}- p_i^0\right|+ \sum_{i=k_0+1}^{k_1}\sum_{j\in [s_i]} p_{i j n} \right)\\
 := & d_n.
\end{align*}
Then by combining the previous three equations, with $m_i=1$ for $i\in \If_1$ and $m_i=2$ for $i\in \If_{>1}$, we obtain
\begin{align*}
&\frac{|G_n \myf-G_0\myf|}{W_2^2(G_n,G_0)} \\
\geq & \left|\sum_{i\in [k_0]}\sum_{1\leq |\alpha|\leq m_i}\frac{D^\alpha \myf (\theta_i^0)}{\alpha!}   A_{i n}(\alpha) + \sum_{i\in [k_0]}\sum_{j\in [s_i]}\frac{ p_{i j n} R_{i j n}}{d_n}+ \sum_{i\in [k_0]} B_{i n} \myf(\theta_i^0) + \sum_{i=k_0+1}^{k_1}\sum_{j\in [s_i]} \frac{p_{i j n}}{d_n} \myf(\theta_{i j n})\right|, \numberthis \label{eqn:ratiolowerbou}
\end{align*}
where $A_{i n}(\alpha)= \frac{\sum_{j\in [s_i]} p_{i j n} \left( \theta_{i j n}-\theta_i^0 \right)^\alpha}{d_n} $ and $B_{i n}=\frac{\sum_{j\in [s_i]} p_{i j n}- p_i^0}{d_n}$ for $i\in [k_0]$. 

Note from above that for each pair $(i,j)$ where $i\in [k_1]$ and $j\in [s_i]$, we have $\frac{ p_{i j n} R_{i j n}}{d_n}\to 0$. Moreover, by taking subsequence if necessary, we also have
\begin{equation}
A_{i n}(\alpha) \to a_{i\alpha} \text{ and } B_{i n}\to b_i \text{ for } i\in [k_0], \quad \frac{p_{i j n}}{d_n} \to g_{ij} \text{ for } k_0+1\leq i\leq k_1,    \label{eqn:coefficientlimitfixargument}
\end{equation}
and at least one of elements in $\{a_{i\alpha}\}_{1\leq |\alpha|\leq m_i, i\in [k_0]}$ or $\{b_i\}_{i\in [k_0]}$ or $\{g_{ij}\}_{j\in [s_i], k_0+1\leq i\leq k_1}$ is not zero.
Denote $g_i=\sum_{j\in [s_i]}g_{ij}$ for $k_0+1\leq i\leq k_1$. Then at least one of elements in $\{a_{i\alpha}\}_{i\in [k_0],1\leq |\alpha|\leq 2}$ or $\{b_i\}_{i\in [k_0]}$ or $\{g_{i}\}_{k_0+1\leq i\leq k_1}$ is not zero since $g_{ij}\geq 0$. 
Since $B_{i n} + \sum_{i=k_0+1}^{k_1} \sum_{j\in [s_i]} \frac{p_{i j n}}{d_n}=0$ for all $n$, it follows that 
\begin{equation}
    \sum_{i\in [k_0]} b_i + \sum_{i=k_0+1}^{k_1} g_i =0. \label{eqn:coefficientsum0}
\end{equation}
Now, from \eqref{eqn:limitzerofixargument}, we have
\begin{align*}
    0=&\lim_{ n\to\infty} \frac{\sup_{\myf\in \myF} |G_n \myf-G_0 \myf|}{W_{2}^{2}(G_n,G_0)} \\
     \geq & \sup_{\myf\in \myF} \liminf_{ n\to\infty} \frac{ |G_n \myf-G_0 \myf|}{W_{2}^{2}(G_n,G_0)}\\
     \geq & \sup_{\myf\in \myF} \left| \sum_{i\in [k_0]}\ \sum_{1\leq |\alpha|\leq m_i}\frac{D^\alpha \myf (\theta_i^0)}{\alpha!}   a_{i\alpha} + \sum_{i\in [k_0]} b_i \myf(\theta_i^0) + \sum_{i=k_0+1}^{k_1}g_i  \myf(\theta_{i}) \right| \numberthis \label{eqn:differentiatedomainfixargument}
\end{align*}
where the last inequality follows from \eqref{eqn:ratiolowerbou} and \eqref{eqn:coefficientlimitfixargument}. That the equations \eqref{eqn:coefficientsum0} and \eqref{eqn:differentiatedomainfixargument} hold with at least one coefficient nonzero contradicts with the hypothesis that $\myF$ is a $(G_0,k)$ second-order linear independent domain.

\end{proof}

\begin{proof}[Proof of Lemma \ref{lem:inverseboundfixoneargumentexactfitted}]

Suppose that \eqref{eqn:localinverseboundfixargumentexactfitted} does not hold. Then there exists a sequence of $G_n\in \Gc_{k_0}(\Theta)$ such that $G_n \overset{W_1}{\to} G_0$ and
\begin{equation}
\lim_{ n\to\infty} \frac{\sup_{\myf\in \myF} |G_n \myf-G_0 \myf|}{W_{1}(G_n,G_0)} = 0. \label{eqn:limitzerofixargumentexafit}
\end{equation}
Write $G_0=\sum_{i\in [k_0]}p_i^0\delta_{\theta_i^0}\in \Ec_{k_0}(\Theta)$. By taking subsequence if necessary, we have that:  $G_n=\sum_{i\in [k_0]} p_{i n} \delta_{\theta_{i n}}$, with 
\begin{align*}
    \theta_{i n} \to \theta_i^0,   & \quad p_{i n} \to p_i^0, \quad \forall \ i\in [k_0] .
\end{align*}

Note that
\begin{align*}
G_n \myf-G_0 \myf =&\sum_{i\in [k_0]}\underbrace{ p_{i  n} \left(\myf(\theta_{i  n}) -\myf(\theta_i^0)\right) }_{:=I_i}+ \sum_{i\in [k_0]} \left( p_{i  n}- p_i^0\right) \myf(\theta_i^0) .  
\end{align*}
By Taylor's theorem, 
\begin{align*}
I_i = &  p_{i  n} \left(\sum_{ |\alpha|= 1}\frac{1}{\alpha!} D^\alpha \myf (\theta_i^0) \left( \theta_{i n}-\theta_i^0 \right)^\alpha + R_{i  n}  \right) ,
\end{align*}
where $R_{in} = o(\|\theta_{i n}-\theta_i^0\|_2)$. 

We also have
\begin{align*}
    W_1(G_n,G_0) 
\leq & \sum_{i\in [k_0]} p^0_{i} \|\theta_{i j n}-\theta_i^0\|_2 + 2\rho  \sum_{i\in [k_0]} \left| p_{i n}- p_i^0\right|  \\
 := & d_n,
\end{align*}
where $\rho:=\max_{1\leq i< j\leq k_0}\|\theta_i^0-\theta_j^0\|_2$. 
Then by combining the previous three equations,  we obtain 
\begin{align*}
&\frac{|G_n \myf-G_0\myf|}{W_1(G_n,G_0)} \\
\geq & \left|\sum_{i\in [k_0]}\sum_{ |\alpha|=1}\frac{D^\alpha \myf (\theta_i^0)}{\alpha!}   A_{i n}(\alpha) + \sum_{i\in [k_0]}\sum_{j\in [s_i]}\frac{ p_{i  n} R_{i  n}}{d_n}+ \sum_{i\in [k_0]} B_{i n} \myf(\theta_i^0)\right| , \numberthis \label{eqn:ratiolowerbouexafit}
\end{align*}
where $A_{i n}(\alpha)= \frac{ p_{i  n} \left( \theta_{i j n}-\theta_i^0 \right)^\alpha}{d_n} $ and $B_{i n}=\frac{  p_{i  n}- p_i^0}{d_n}$ for $i\in [k_0]$. 
Note that for any $i\in [k_0]$ we have $\frac{ p_{i  n} R_{i  n}}{d_n}\to 0$. Moreover, by taking subsequence if necessary, we also have
\begin{equation}
A_{i n}(\alpha) \to a_{i\alpha} \text{ and } B_{i n}\to b_i \text{ for } i\in [k_0], \label{eqn:coefficientlimitfixargumentexafit}
\end{equation}
and at least one of the elements in $\{a_{i\alpha}\}_{1\leq |\alpha|\leq m_i, i\in [k_0]}$ or $\{b_i\}_{i\in [k_0]}$  is not zero.
It also holds that 
\begin{equation}
    \sum_{i\in [k_0]} b_i  =0. \label{eqn:coefficientsum0exafit}
\end{equation}
since $B_{i n} =0$ for all $n$.

Now, from \eqref{eqn:limitzerofixargumentexafit},
\begin{align*}
    0=&\lim_{ n\to\infty} \frac{\sup_{\myf\in \myF} |G_n \myf-G_0 \myf|}{W_{1}(G_n,G_0)} \\
     \geq & \sup_{\myf\in \myF} \liminf_{ n\to\infty} \frac{ |G_n \myf-G_0 \myf|}{W_1(G_n,G_0)}\\
     \geq & \sup_{\myf\in \myF} \left| \sum_{i\in [k_0]}\ \sum_{ |\alpha|=1}\frac{D^\alpha \myf (\theta_i^0)}{\alpha!}   a_{i\alpha} + \sum_{i\in [k_0]} b_i \myf(\theta_i^0)    \right| \numberthis \label{eqn:differentiatedomainfixargumentexafit}
\end{align*}
where the last inequality follows from \eqref{eqn:ratiolowerbouexafit} and \eqref{eqn:coefficientlimitfixargumentexafit}. That the equations \eqref{eqn:coefficientsum0exafit} and \eqref{eqn:differentiatedomainfixargumentexafit} hold with at least one coefficient nonzero contradicts with the hypothesis that $\myF$ is a $(G_0,k_0)$ first-order linear independence domain.

\end{proof}

Lemma \ref{lem:compactness} is reproduced and proved below.

\begin{lem} 
Suppose that $\Theta=\Rb^q$ and the function class $\myF$ is uniformly bounded, i.e. $\sup_{\myf\in \myF} \sup_{\theta\in \Theta}|\myf(\theta)|<\infty$.  Consider $G_0\in \Ec_{k_0}(\Theta)$ and $k>k_0$. Then for any $r>0$, 
\begin{equation}
\liminf_{ \substack{G\overset{W_1}{\to} G_0\\ G\in \Gc_k(\Theta) }} \frac{\sup_{\myf\in \myF} |G \myf-G_0 \myf|}{W_{r}^{r}(G,G_0)} =0. 
\end{equation}
\end{lem}
\begin{proof}[Proof of Lemma \ref{lem:compactness}]
Write $G_0=  \sum_{i\in [k_0]} p_i^0 \delta_{\theta_i^0}   $. Without loss of generality assume that $\theta_{k_0}^0$ has the largest first coordinate. Consider $\theta_n=\theta_{k_0}^0+n^{\frac{1}{2r}} e_1$ where $e_1$ is the vector with $1$ on first coordinate and $0$ on other coordinates. 
Consider $G_n=\sum_{i\in [k_0-1]} p_i^0 \delta_{\theta_i^0}  +  (p_{k_0}^0 - \frac{1}{n}) \delta_{\theta_{k_0}^0} + \frac{1}{n} \delta_{\theta_n}  $. Since $\|\theta_n-\theta_{k_0}^0\|_2\leq \|\theta_n-\theta_{i}^0\|_2 $ for $i\in [k_0]$, we have $W_r^r(G_n,G_0)= \frac{1}{n}\|\theta_n- \theta_{k_0}^0\|_2^r = \frac{1}{\sqrt{n}}\to 0$. On the other hand, 
$$
G_n\myf - G_0 \myf = \frac{1}{n} \left(\myf(\theta_n) - \myf(\theta_{k_0}^0) \right)
$$
and thus 
$$
\frac{\sup_{\myf\in \myF} |G_n \myf-G_0 \myf|}{W_{r}^{r}(G_n,G_0)} = \frac{\sup_{\myf\in \myF} |\myf(\theta_n) - \myf(\theta_{k_0}^0)|}{\sqrt{n}} \to 0,
$$
where the last step follows from that $\myF$ is uniformly bounded.
\end{proof}

\begin{proof}[Proof of Theorem \ref{thm:fixedG_0estimator}]
(a)
    By \eqref{eqn:localinverseboundfixargument}, there exist $r,\epsilon>0$ such that for any $H\in B_{W_1}(G_0,r)$, the $W_1$-ball centering at $G_0$ of radius $r$ in $\Gc_{k_0}(\Theta)$, we have
\begin{equation}
\sup_{\myf\in \myF} |G_0 \myf-H \myf|\geq \epsilon W_{1}(G_0,H).  \label{eqn:inverseboundfixcon1}
\end{equation}

Define the constant
$$
z := \inf_{ \substack{ H \in \Gc_{k_0}(\Theta) \setminus B_{W_1}(G_0,r) } }\ \  \sup_{\myf\in \myF} |G_0 \myf-H \myf|.
$$
Since $\sup_{\myf\in \myF} |G_0 \myf-H \myf|$ is lower semicontinuous 
on the compact set $\Gc_{k_0}(\Theta) \setminus B_{W_1}(G_0,r)$, the infimum is attained. Since $G_{k_0}(\Theta)$ is distinguishable by $\myF$, we have $z>0$.

Set the event
$$A_{G_0}(z):= \left\{ \sup_{\myf\in\myF}\left|G_0  \myf-\frac{1}{n}\sum_{i\in [n]}t_\myf(X_i)\right| \leq \frac{1}{4}z \right\}.$$
Then on the event $A_{G_0}(z) \cap  \{\hat{k}_n = k_0\}$, by triangle inequality and the definition of $\hat{G}_n(\ell)$ we have
\begin{equation}
\sup_{\myf\in\myF}\left|G_0  \myf-\hat{G}_n(\hat{k}_n)\myf \right| \leq 2 \sup_{\myf\in\myF}\left|G_0  \myf-\frac{1}{n}\sum_{i\in [n]}t_\myf(X_i)\right| \leq \frac{1}{2}z, \label{eqn:constraintuppbou}
\end{equation}
which then implies that $\hat{G}_n(\hat{k}_n) \in B_{W_1}(G_0,r)$ by our choice of $z$. Thus on the event $A_{G_0}(z) \cap  \{\hat{k}_n=  k_0\}$, by \eqref{eqn:inverseboundfixcon1} and \eqref{eqn:constraintuppbou},
\begin{align*}
    W_{1}(G_0,\hat{G}_n(\hat{k}_n)) \leq \frac{1}{\epsilon} \sup_{\myf\in \myF} |G_0 \myf-\hat{G}_n(\hat{k}_n)) \myf| \leq \frac{2}{\epsilon} \sup_{\myf\in\myF}\left|G_0  \myf-\frac{1}{n}\sum_{i\in [n]}t_\myf(X_i)\right|.
\end{align*}

Denote
\begin{align*}
J_{G_0}(z) := \left(A_{G_0}(z) \cap  \{\hat{k}_n=  k_0\}\right)^c = &  \{\hat{k}_n\neq k_0\} \bigcup A_{G_0}^c(z)  . 
\end{align*}
Then we have 
\begin{equation}
W_{1}(G_0,\hat{G}_n(\hat{k}_n)) \leq  \frac{2}{\epsilon} \sup_{\myf\in\myF}\left|G_0  \myf-\frac{1}{n}\sum_{i\in [n]}t_\myf(X_i)\right|  + \diam(\Theta) 1_{J_{G_0}(z)}, \label{eqn:W1uppbou}
\end{equation}
where we use that $W_{1}(G_0,\hat{G}_n(\hat{k}_n))\leq \diam(\Theta) $. 
It follows that
\begin{align*}
    \{W_{1}(G_0,\hat{G}_n(\hat{k}_n)) \geq t \} 
    \subset & \left\{   \sup_{\myf\in\myF}\left|G_0  \myf-\frac{1}{n}\sum_{i\in [n]}t_\myf(X_i)\right| \geq \frac{\epsilon}{2} t \right\} \bigcup  J_{G_0}(z) \\
    = & \left\{ \sup_{\myf\in\myF}\left|G_0  \myf-\frac{1}{n}\sum_{i\in [n]}t_\myf(X_i)\right| \geq \min\left\{\frac{\epsilon}{2} t, \frac{1}{4}z\right\} \right\} \bigcup \{\hat{k}_n\neq k_0\} \\
    \subset & \left\{ \sup_{\myf\in\myF}\left|G_0  \myf-\frac{1}{n}\sum_{i\in [n]}t_\myf(X_i)\right| \geq \min\left\{\epsilon_1 t,\epsilon_1\right\} \right\} \bigcup \{\hat{k}_n\neq k_0\}, \numberthis \label{eqn:setinclusion}
\end{align*}
where in the last step $\epsilon_1:=\min\left\{\frac{\epsilon}{2},\frac{1}{4} z\right\}$.

(b)
Apply Lemma \ref{lem:componentnumberestimation} with $G=G_0$, we get
 \begin{equation}
    \{\hat{k}_n\neq k_0\}   \subset \left\{ \sup_{\myf\in\myF}\left|G  \myf-\frac{1}{n}\sum_{i\in [n]}t_\myf(X_i)\right| \geq \min\{a_n,b_{G_0}-a_n\} \right\}. \label{eqn:notequalsign}
\end{equation}
Then
\begin{align*}
J_{G_0}(z_n) 
\subset & \left\{ \sup_{\myf\in\myF}\left|G_0  \myf-\frac{1}{n}\sum_{i\in [n]}t_\myf(X_i)\right| \geq \min\{\frac{1}{4}z, a_n,b_{G_0}-a_n\} \right\} \\
\subset & \left\{ \sup_{\myf\in\myF}\left|G_0  \myf-\frac{1}{n}\sum_{i\in [n]}t_\myf(X_i)\right| \geq  \min\{a_n,\epsilon'_0-a_n\} \right\} , \numberthis \label{eqn:fixedG_0temp2}
\end{align*}
where the second set inclusion is achieved by setting $\epsilon'_0=\frac{1}{4}z\wedge b_{G_0}$. Combining the previous equation with \eqref{eqn:W1uppbou}, the conclusion on $\Eb_{G^*}W_1(G^*,\hat{G}_n(\hat{k}_n))$ is completed.

By \eqref{eqn:setinclusion} and \eqref{eqn:notequalsign},
\begin{align*}
    \{W_{1}(G_0,\hat{G}_n(\hat{k}_n)) \geq t \} 
    \subset & \left\{   \sup_{\myf\in\myF}\left|G_0  \myf-\frac{1}{n}\sum_{i\in [n]}t_\myf(X_i)\right| \geq \min\left\{\epsilon_1 t,\epsilon_1,a_n,b_{G_0}-a_n\right\} \right\} \\
    \subset & \left\{ \sup_{\myf\in\myF}\left|G_0  \myf-\frac{1}{n}\sum_{i\in [n]}t_\myf(X_i)\right| \geq \min\left\{\epsilon_0 t, a_n,\epsilon_0-a_n\right\} \right\},
\end{align*}
where in the last step $\epsilon_0=\min\left\{b_{G_0},\epsilon_1\right\}$.
\end{proof}

\begin{proof}[Proof of Lemma \ref{lem:empiricalprocessMMDconcentration}]
    Note that 
    $$
    \MMD(\Pb, \hat{\Pb}_n) = \sup_{f_1 \in \Fc_1}\left|\int f_1 d\Pb-\frac{1}{n}\sum_{i\in [n]}f_1(X_i)\right|:=g(X_1,\ldots,X_n)
    $$ 
    with $\Fc_1$ to be the unit ball in the associate RKHS $\Hc$. For any $f_1\in \Fc_1$,
    $$
    |f_1(x)|=|\langle f_1, \ker(x,\cdot) \rangle| \leq \|f_1\|_{\Hc} \|\ker(x,\cdot)\|_{\Hc} = \|\ker\|_\infty, 
    $$
    so $\|f_1\|_\infty\leq \|\ker\|_\infty$. It is then easy to see that for any $i$
    $$
    |g(X_1,\ldots,X_n) - g(X_1,\ldots,X_{i-1},Y_i,X_{i+1},\ldots,X_n) | \leq \frac{2\sup_{f_1\in \Fc_1}\|f_1\|_\infty}{n}\leq  \frac{2\|\ker\|_\infty}{n}.
    $$
    By McDiarmid's bounded difference inequality, we then have 
    $$ \Pb  \left( \MMD(\Pb, \hat{\Pb}_n) \geq \Eb \MMD(\Pb, \hat{\Pb}_n) + \epsilon \right)
    \leq 2\exp\left(-\frac{n\epsilon^2}{2 \|\ker\|^2_\infty}\right). $$
   The proof is then completed by combining the above inequality and Lemma \ref{lem:empiricalprocessMMD}.
\end{proof}

\end{document}